\documentclass[pdflatex,sn-apa]{sn-jnl}% Basic Springer Nature Reference Style/Chemistry Reference Style
\usepackage{graphicx,color}
\usepackage{amsmath,amsthm,amssymb,amsfonts,mathtools}
\usepackage{epsfig}
\usepackage{booktabs}
\usepackage{enumitem}
\usepackage{bm}
\usepackage{setspace}
\usepackage{threeparttable}
\usepackage{array}
\usepackage{xr}
\shortcites{cuesta2019,adler1998,GHK,meintanis2023,pfister2018,lindsay2008,HJM,giacomini2013,chen2022,nolan2001,byczkowski1993,cadirci2022,kokoszka2017,rice2020,bugni2009,nieto2014,barcenas2017,jiang2019,gorecki2020}

\jyear{2021}%

\theoremstyle{thmstyleone}%
\newtheorem{thm}{Theorem}[section]

\newtheorem{prop}[thm]{Proposition}

\newcommand{\be}{\begin{eqnarray}}
\newcommand{\ee}{\end{eqnarray}}
\newcommand{\bq}{\begin{eqnarray*}}
\newcommand{\eq}{\end{eqnarray*}}

\DeclareMathOperator*{\argmax}{arg\,max}

\newcommand{\verk}{\stackrel{{\cal D}}{\longrightarrow}}

\renewcommand{\hat}{\widehat}
\renewcommand{\tilde}{\widetilde}

\DeclareMathOperator{\Prob}{P}

\DeclareMathOperator{\sign}{sign}

\newcommand{\stk}{\stackrel{\Prob}{\longrightarrow}}

\renewcommand{\vec}[1]{\bm{#1}}
\newcommand{\mat}[1]{\bf{#1}}	

\DeclareRobustCommand{\S}{\ifmmode\mathsection\else\textsection\fi}

\newcolumntype{H}{>{\setbox0=\hbox\bgroup}c<{\egroup}@{}} 

\begin{document}

\title[Specification procedures for multivariate stable-Paretian laws]{Specification procedures for multivariate stable-Paretian laws for independent and for conditionally heteroskedastic data}


\author*[1,2]{\fnm{Simos G.} \sur{Meintanis}}\email{simosmei@econ.uoa.gr}

\author[3]{\fnm{John P.} \sur{Nolan}}

\author[4]{\fnm{Charl} \sur{Pretorius}}

%\author[5]{\fnm{Zhou} \sur{Zhou}}

\affil*[1]{\orgdiv{Department of Economics}, \orgname{National and Kapodistrian University of Athens}, \orgaddress{\city{Athens}, \country{Greece}}}

\affil[2]{\orgdiv{Pure and Applied Analytics}, \orgname{North-West University}, \orgaddress{\city{Potchefstroom}, \country{South Africa}}}

\affil[3]{\orgdiv{Department of Statistics}, \orgname{American University}, \orgaddress{\city{Washington DC}, \country{USA}}}

\affil[4]{\orgdiv{Centre for Business Mathematics and Informatics}, \orgname{North-West University}, \orgaddress{\city{Potchefstroom}, \country{South Africa}}}

%\affil[5]{\orgdiv{Department of Statistical Sciences}, \orgname{University of Toronto}, \orgaddress{\state{Ontario}, \country{Canada}}}


\abstract{We consider goodness-of-fit  methods for multivariate symmetric and asymmetric stable Paretian random vectors in arbitrary dimension. The methods are based on the empirical characteristic function and are implemented both in the i.i.d.  context as well as for innovations in GARCH models. Asymptotic properties of the proposed procedures are discussed, while the finite-sample properties are illustrated by means of an extensive Monte Carlo study. The procedures are also applied to real data from the financial markets.}

\keywords{Empirical characteristic function, Goodness-of-fit, Heavy-tailed distribution}



\maketitle

\section{Introduction}\label{sec_1}

Stable-Paretian (SP) distributions are extremely important from the theoretical point of view   as they are  closed under convolution, and also they are the only possible limit laws for normalized sums of i.i.d.\ random variables. In fact, this last feature makes SP laws particularly appealing for financial applications since stock returns and other financial assets often come in the form of sums of a large numbers of independent terms. Moreover the empirical density of such assets is leptokurtic  and in many cases skewed, thus making the family of SP laws particularly suited for related applications; see for instance the celebrated  stability hypothesis that goes back at least to  \citet{mandelbrot1963} and \citet{fama1965}.

The above findings prompted further research on the stochastic properties and inference of SP laws and on their application potential. The reader is referred to  \citet{samarodnitsky1994}, \cite{adler1998}, \cite{uchaikin1998chance}, \citet{rachev2000} and \citet{nolan2020} for an overview on the stochastic theory, statistical inference and applications of SP laws. %The recent volume by \citet{nolan2020} contains many results on theory, estimation and testing techniques and are also good sources for recent contributions on inferential procedures and applications of SP laws.  

In the aforementioned works, statistical inference and applications are mostly restricted to univariate SP laws. The respective topics for SP random vectors are much less explored. In this connection, certain distributional aspects of SP random vectors are discussed in \citet{press1972multivariate}, while \citet{press1972estimation} defines moment-type estimators of the parameters of a multivariate SP distribution utilizing the characteristic function (CF).  Maximum-likelihood based estimation is discussed by \citet{ogata2013}, and \citet{nolan2013}, Bayesian methods in \citet{tsionas2016}, \citet{lombardi2009} use indirect estimation methods, whereas \cite{koutrouvelis1980}, \citet{nolan2013} and \citet{sathe2019} consider CF-regression methods. As far as testing is concerned, the only available formal method seems to be the CF-based goodness-of-fit test of \citet{MNT}. % while for multivariate GARCH modeling with SP innovations we refer to \citet{bonato2012}.     

In this article we propose CF-based goodness-of-fit procedures for SP random vectors in the elliptically symmetric case. Moreover we also consider the asymmetric case for which to the best of our knowledge goodness-of-fit  have not been considered before. The remainder of the article unfolds as follows. In Section~\ref{sec_2} we introduce a general  goodness-of-fit test for  elliptically symmetric SP laws. In Section \ref{sec_3} we particularize this test in terms of  computation.  Section~\ref{sec_4} addresses the problem of estimation of the SP parameters and  the study of asymptotics of the proposed test, while in Section \ref{sec_5} we extend the test to a multivariate GARCH model with SP innovations.  The results of an extensive simulation study illustrating the finite-sample properties of the method are presented in Section \ref{sec_6}. In Section \ref{sec_7} the case of asymmetric SP law  is considered, with known as well as unknown characteristic exponent,  by means of a different testing procedure that avoids integration over a complicated  CF and computation of the corresponding density.   The impact of high dimension on the test statistic is also discussed in this section. Applications are given in Section \ref{sec_8}, and we conclude in Section \ref{sec_9} with a discussion. 
Some technical material is deferred to an online supplement along with additional Monte Carlo results.

\section{Goodness-of-fit tests} \label{sec_2}
Let $\vec X$  be a random vector in general dimension $p\geq 1$, with CF $\varphi(\vec t)=\mathbb E(e^{{\rm{i}} \vec t^\top \vec X}), \ \vec t \in \mathbb R^p, \ {\rm{i}}=\sqrt{-1}$. Here we consider  goodness-of-fit tests for the elliptically symmetric SP law. In this connection we note that SP random vectors are parameterized by the triplet $(\alpha,{\vec {\delta}}, {\mat{Q}})$, where $\alpha\in(0,2]$ denotes the characteristic exponent, and   ${\vec {\delta}} \in \mathbb R^p$ and ${\mat{Q}} \in {\cal {M}}^p$, are location vector and dispersion matrix, respectively,  with  ${\cal {M}}^p$ being the set of all symmetric positive definite $(p\times p)$ matrices. On the basis of i.i.d. copies  $(\bm X_j, \ j=1,...,n)$ of $\bm X$ we wish to test the null hypothesis 
\begin{equation}\label{null} {\cal{H}}_{0}: \vec X \sim {\cal{S}}_\alpha(\vec \delta,\vec Q),   \text{ for some }  (\bm{\delta},{\bf{Q}}) \in  \mathbb R^p\times {\cal {M}}^p,  \end{equation}  
where we write $\vec X\sim {\cal{S}}_\alpha(\vec \delta,\vec Q)$ to denote that $\vec X$ follows a SP law with the indicated parameters. 
For subsequent use we mention that if $\vec X \sim {\cal{S}}_\alpha(\vec \delta,\vec Q)$, then it admits the stochastic representation
\begin{equation} \label{sr} 
\vec X=\vec \delta +A^{1/2}\bm N,
\end{equation}
where $A$ is a totally skewed to the right SP random variable with characteristic exponent $\alpha/2$ ($\alpha<2$)  and $\bm N$ is a zero-mean Gaussian vector with covariance matrix $\vec Q$, independent of $A$; see \citet[\S2.5]{samarodnitsky1994}.

Our test will make use of the fact that if $\vec X \sim {\cal{S}}_\alpha(\vec \delta,\vec Q)$, then the CF of $\vec X$  is given by  
$\varphi_\alpha({\bm{t}};{\bm{\delta}},{\bf{Q}})=e^{{\rm{i}} {\bm {t}}^\top {\bm{\delta}}-({\bm{t}}^\top {\bf{Q}} {\bm{t}})^{\alpha/2}} $. The cases $\alpha=2$ and $\alpha=1$, respectively, correspond to the most prominent members of the SP family, i.e. the Gaussian and Cauchy distributions, while  
 $\phi_{\alpha}(\bm t)=\exp({-\|\bm{t}\|^{\alpha}})$, will denote the CF of a spherical SP law, i.e. an SP law  with location $\bm \delta=\vec 0$ and dispersion matrix $\bf Q$ set equal to the identity matrix. 

% \begin{equation} \label{cfnull} 
% \varphi(\bm t)=\varphi_\alpha({\bm{t}};{\bm{\delta}},{\bf{Q}}), \ \bm t \in \mathbb R^p, 
% \end{equation}
% where
% \begin{equation}\label{cfnull1}
% \varphi_\alpha({\bm{t}};{\bm{\delta}},{\bf{Q}})=e^{{\rm{i}} {\bm {t}}^\top {\bm{\delta}}-({\bm{t}}^\top {\bf{Q}} {\bm{t}})^{\alpha/2}}
% \end{equation}

As it is already implicit in \eqref{null} the parameters $(\bm \delta,\bf Q)$ are considered unknown, and hence they will be replaced by estimators  $(\widehat {\bm  \delta}_n,\widehat {\bf Q}_n)$ and the test procedure will be applied on the standardized data  
\begin{equation}\label{stand}
{\bm{Y}}_j=\widehat {\bf{Q}}^{-1/2}_n({\bm {X}}_j-\widehat {\bm{\delta}}_n), \quad j=1,...,n.
\end{equation} %involving  estimators  $(\widehat {\bm  \delta}_n,\widehat {\bf Q}_n)$ of $(\bm \delta, \bf Q)$ that are obtained from $(\bm X_j, \ j=1,...,n)$.  

Specifically for a non-negative weight function $w(\cdot)$  we propose the test criterion 
\begin{equation} \label{TS}
T_{n,w}=n \int_{\mathbb R^p} \left\lvert\phi_n(\bm t)-\phi_{\alpha_0}(\bm t)\right\rvert^2 w(\bm t) {\rm{d}} \bm t, 
\end{equation}
where $\phi_n(\vec t)=n^{-1} \sum_{j=1}^n \exp({{\rm{i}} \vec t ^\top \vec Y_j})$
is the empirical CF computed from  $(\bm Y_j, \ j=1,...,n)$. Here the characteristic exponent $\alpha$ is considered fixed at value $\alpha=\alpha_0$, while the case of unknown $\alpha$ will be considered in Section \ref{sec_7}.

\section{Computational aspects} \label{sec_3}
 It is already transparent that the test statistic $T_{n,w}$ depends on the weight function $w(\cdot)$, the choice of which we consider in this section. Specifically, from \eqref{TS} we have by simple algebra
\begin{equation}\label{TS1}
T_{n,w}=\frac{1}{n} \sum_{j,k=1}^n I_w{(\vec Y_j-\vec Y_k)}+n J_w-2 \sum_{j=1}^n K_w(\vec Y_j), 
\end{equation}
where 
% \begin{align}
% I_w(\bm x)&=\int_{\mathbb R^p} \cos(\bm t^\top \bm x)w(\bm t) {\rm{d}} \bm t, \label{int1} \\
% K_w(\bm x)&=\int_{\mathbb R^p} \cos(\bm t^\top \bm x)\phi_{\alpha_0}(\bm t)w(\bm t) {\rm{d}} \bm t, \label{int2} \\
% J_w&=\int_{\mathbb R^p} \phi^2_{\alpha_0}(\bm t)w(\bm t) {\rm{d}} \bm t. \label{int3}
% \end{align}
\begin{equation}\begin{split} \label{int1}
I_w(\bm x)&=\int_{\mathbb R^p} \cos(\bm t^\top \bm x)w(\bm t) {\rm{d}} \bm t, \quad
J_w=\int_{\mathbb R^p} \phi^2_{\alpha_0}(\bm t)w(\bm t) {\rm{d}} \bm t, \\
K_w(\bm x)&=\int_{\mathbb R^p} \cos(\bm t^\top \bm x)\phi_{\alpha_0}(\bm t)w(\bm t) {\rm{d}} \bm t.
\end{split}\end{equation}
\subsection{Using the CF of the Kotz-type distribution} \label{sec_31}
We now discuss the computation of the test statistic figuring in \eqref{TS1}. In doing so we will make use of an appropriate weight function $w(\cdot)$ that facilitates explicit representations of the integrals in \eqref{int1}. Specifically we choose  $w(\bm t)=(\|\bm t \|^2)^{\nu}e^{-(\|\bm{t}\|^2)^{\alpha_0/2}}$ as weight function, which for $(r,s)=(1,\alpha_0/2)$ is proportional to the density of the spherical Kotz-type distribution $c(\|\bm x\|^2)^{\nu} e^{-r (\|\bm x\|^2)^s}$,  with $c$ a normalizing constant; see \citet{nadarajah2003}.  With this weight function the integrals in   \eqref{int1} may be derived as special cases of the integral
\begin{equation} \label{int4}
I_{\nu,r}(\bm x;s)=\int_{\mathbb R^p} \cos(\bm t^\top \bm x) \left(\|\bm t\|^2\right)^{\nu}
e^{-r (\|\bm t\|^2)^s} {\rm{d}} \bm t, \ r>0, \ \nu=0,1,\ldots, 
\end{equation}
with the cases of interest being $0<s(=\alpha_0/2)\leq 1$. In turn this integral can be computed  by making use of the CF of the Kotz-type distribution in a form of an absolutely convergent series of the type (see \citealp{streit1991}; \citealp{Iyengar1989}; \citealp{kotz1994})
$
I_{\nu,r}(\bm x;s)=\sum_{k=0}^\infty \psi_k(\|\bm x\|^2), 
$
which for selected values of $s$ reduces to a finite sum. For more details we refer to Section~S1 of the online supplement. 
%as $I_w(\bm x)=I_1(\bm x;\alpha_0/2)$, $K_w(\bm x;\alpha_0)=I_2(\bm x;\alpha_0/2)$, and $J_w(\alpha_0)=I_3(\bm 0;\alpha_0/2)$, respectively. 

Given $I_{\nu,r}(\cdot;\cdot)$, the test statistic figuring in \eqref{TS1} may be written as 
\begin{equation} T_{n,w}=\frac{1}{n} \sum_{j,k=1}^n I_{\nu,r}\left(\bm Y_j-\bm Y_k;\frac{\alpha_0}{2}\right)+n I_{\nu,r+2}\left(\bm 0;\frac{\alpha_0}{2}\right)-2 \sum_{j=1}^n I_{\nu,r+1}\left(\bm Y_j;\frac{\alpha_0}{2}\right). 
\label{TS2}
\end{equation}
%where $I_{\nu,r}(\bm x;\frac{\alpha_0}{2})$, is given in the onby \eqref{int5}--\eqref{int11}
% depending on the arguments $(r,\bm x,\alpha_0)$ and the dimension $p$.

\subsection{Using the inversion theorem}
In this section we will use the inversion theorem for CFs in order to compute the integrals defined  by  \eqref{int1}. Specifically for an absolutely integrable CF $\varphi(t)$, the inversion theorem renders the density $f(\cdot)$ corresponding to $\varphi(\cdot)$  as 
\begin{equation}\label{eq:fourier.inv}
f(\bm x)=\frac{1}{(2\pi)^p} \int_{\mathbb R^p} e^{- {\rm{i}} \bm t^\top \bm x}\varphi(\bm t)
{\rm{d}} \bm t.  
\end{equation}
%see \cite[Theorem~1.8.5]{ushakov}.
%by means of which the density $f(\cdot)$ of a given random variable in $\mathbb R^p$ may be computed from  the corresponding CF $\varphi(\cdot)$.  
In this connection, we start from the expression of the test statistic in \eqref{TS1}, and adopt the weight function $w(\bm t)=e^{-r \|\bm{t}\|^{\alpha_0}}$. This choice amounts to taking $\nu=0$ in the Kotz-type density, which is the same as if we incorporate the CF $\phi_{\alpha_0}(\cdot)$ of the SP law under test in the weight function.

With this weight function, the statistic figuring in \eqref{TS1}, say $T_{n,r}$, becomes
% With this weight function, the test statistic figuring in \eqref{TS1}, say $T_{n,r}$, is rendered as
\begin{equation}\label{STS}
T_{n,r}=\frac{1}{n} \sum_{j,k=1}^n \Lambda_r(\bm Y_j-\bm Y_k;\alpha_0) +n\Lambda_{r+2}(\bm 0;\alpha_0) -2 \sum_{j=1}^n  \Lambda_{r+1}(\bm Y_j;\alpha_0),
\end{equation}
where  by making use of the inversion theorem in \eqref{eq:fourier.inv},
\begin{equation} \label{lambda}
\Lambda_r(\bm x;\alpha_0)=\int_{\mathbb R^p} \cos(\bm t^\top \bm x) e^{-r \|\bm{t}\|^{\alpha_0}}
{\rm{d}} \bm t= \frac{(2\pi)^p}{r^{\frac{p}{\alpha_0}}} f_{\alpha_0}\left(\frac{\bm x}{r^{\frac{1}{\alpha_0}}}\right), 
\end{equation}
 with  $f_\alpha(\cdot)$ being the density of the  spherical SP law with CF  $\phi_{\alpha}(\cdot)$. %For details on the computation of the density of the spherical SP law we refer to the on--line supplement.  
% Notice that clearly from \eqref{int4} we have $\Lambda_r(\cdot;\alpha_0)=I_{0,r}(\cdot;\alpha_0/2)$.  

% {\color{magenta}From \eqref{lambda} for $\alpha_0=2$ we obtain 
% \[
% \Lambda_r(\bm x;2)= \frac{(2\pi)^p}{r^{\frac{p}{2}}} f_{2}\left(\frac{\bm x}{r^{\frac{1}{2}}}\right), 
% \]
%  with  
%  \[f_{2}(\bm z)=\frac{1}{(2\pi)^{p/2}} e^{-\|\bm z\|^2/2}
%  \]
%  being the density of the  p-variate standard normal distribution, so that
%  \[
% \Lambda_r(\bm x;2)= \left(\frac{2\pi}{r}\right)^{\frac{p}{2}} e^{-\|\bm x\|^2/(2r)}, 
% \]
% Now consider the test statistic with no standardization, i.e. replace in \eqref{STS} $\bm Y_j$ by $\bm X_j, \ j=1,...,n$,
%  \begin{eqnarray}\label{extra3} \nonumber
% T_{n,r}&=&\frac{1}{n} \left (\sum_{j=1}^n \Lambda_r(\bm 0;2)+\sum_{j\neq k}^n \Lambda_r(\bm X_j-\bm X_k;2)\right) +n\Lambda_{r+2}(\bm 0;2) -2 \sum_{j=1}^n  \Lambda_{r+1}(\bm X_j;2) \\ \nonumber 
% &=& \frac{1}{n} \left (\sum_{j=1}^n\left(\frac{2\pi}{r}\right)^{\frac{p}{2}} +\sum_{j\neq k}^n \left(\frac{2\pi}{r}\right)^{\frac{p}{2}}  e^{-\|\bm X_j-\bm X_k\|^2/(2r)}\right)+n \left(\frac{2\pi}{r+2}\right)^{\frac{p}{2}}\\ \nonumber &-& 2 \sum_{j=1}^n   \left(\frac{2\pi}{r+1}\right)^{\frac{p}{2}} e^{-\|\bm X_j\|^2/(2(r+1))}.
% \end{eqnarray}
% Now the quantity $e^{-\|x\|^2}$ behaves asymptotically as $e^{-p^2}$, as $p\to\infty$, so that    
% \[
% T_{n,r} \approx \left(\frac{2\pi}{r}\right)^{\frac{p}{2}} +n \left(\frac{2\pi}{r+2}\right)^{\frac{p}{2}},
% \]
% from which it follows that if $r<2\pi$, $r=2\pi$, $r>2\pi$, we have  
% \[
% T_{n,r} \to \infty, 1, 0, \ p \to \infty, 
% \]
% respectively. 
%  } %For

\section{On estimation of parameters and limit properties of the test} \label{sec_4}

\subsection{Estimation of parameters} \label{sec_4.1}    

The parameters $\bm\delta$ and $\bf Q$ in \eqref{null} are assumed unknown and need to be estimated.
Seeing that reliable procedures for calculation of stable densities are available,
we use maximum likelihood estimation.
To avoid searching over the space of positive definite matrices,
% we follow \citet[see Section~2.3]{nolan2013} and use the estimators
% $\hat{\bm\delta}_n$ and $\hat{\bf Q}_n=\hat{\bm L}_n \hat{\bm L}_n^\top$
% given by
we use the estimators
$\hat{\bm\delta}_n$ and $\hat{\bf Q}_n=\hat{\bm L}_n \hat{\bm L}_n^\top$
given by
\begin{equation} \label{MLE}
	\left(\hat{\bm\delta}_n, \hat{\bm L}_n\right)
	= \argmax_{\bm\delta\in\mathbb{R}^p, \bm L \in\mathbb{L}^p} \;
			\log\left[\lvert\det \bm L\rvert^{-n} \prod_{j=1}^n f_{\alpha_0}\left(\bm L^{-1}(\bm X_j - \bm\delta)\right)\right],
\end{equation}
where $\mathbb{L}^p$ denotes the space of lower triangular $p\times p$ matrices,
%\[
	%\operatorname{lik}(\bm\delta, \bm L)
	%= |\det \bm L|^{-n} \prod_{j=1}^n f_{\alpha_0}\left(\bm L^{-1}(\bm X_j - \bm\delta)\right),
%\]
and, as before, $f_{\alpha}(\cdot)$ denotes the density of a $p$-dimensional stable
distribution with CF $\phi_{\alpha_0}(\cdot)$.
Initial values for the optimization procedure are obtained using projection estimators of
$\bm\delta$ and $\bf Q$ as outlined in \citet[Section~2.3]{nolan2013}.
%We first estimate the $p$ location parameters
%
%The procedure is based on the fact that, for any $\bm t\in\mathbb{R}^p$,
%the linear combination $\vec t^\top\vec X$ is univariate $\alpha$-stable
%with location $\vec t^\top\vec\delta$.
%
 %if $\vec X$ has CF $e^{-(\vec t^\top\mat Q\vec t)^{\alpha/2}}$,
%$\vec t\in\mathbb{R}^p$,
%then $\vec t^\top\vec X$ is univariate $\alpha$ stable
%with (squared) scale
%\[
	%\left[-\ln\Expec e^{i\vec t^\top\vec X}\right]^{2/\alpha}
	%= \vec t^\top\mat Q\vec t.
%\]

\subsection{Limit null distribution and consistency} \label{subsec_4.2}    
We present here the main elements involved in the limit behavior of the test statistic $T_{n,w}$. In this connection and despite the fact that, as already emphasized, it is computationally convenient to use a weight function that is proportional to the density of a spherical Kotz-type distribution,  our limit results apply under a general weight function satisfying certain assumptions and, under given regularity conditions, with arbitrary estimators of the distributional parameters. 

Specifically, we assume that the weight function satisfies $w(\bm t)>0$ (apart from a set of measure zero),   $w(\bm t)=w(-\bm t)$ and $\int_{\mathbb R^p} w(\bm t) {\rm{d}}\bm t<\infty$.

We also suppose that the estimators figuring in the standardization defined by \eqref{stand}  admit the Bahadur-type asymptotic representations
\[
\sqrt{n} \left(\widehat {\bm \delta}_n-\bm \delta_0\right)=\frac{1}{\sqrt{n}} \sum_{j=1}^n \bm \ell_\delta(\bm X_j)+{\rm{o}}_{\mathbb P}(1)
\]
and
\[
\sqrt{n} \left(\widehat {\bf Q}_n-\bf Q_0\right)=\frac{1}{\sqrt{n}} \sum_{j=1}^n \bm \ell_{\bf Q}(\bm X_j)+{\rm{o}}_{\mathbb P}(1),
\]
where $(\bm \delta_0,\bf Q_0)$ denote true values, with $\bm X_1\sim {\cal{S}}_{\alpha_0}(\bm \delta_0,\bf Q_0)$ and $\left({\bm {\ell}}_{\bf \delta}({\bm {X}}_1),{\bf {\ell}}_{\bf Q}({\bm {X}}_1)\right)$ satisfying  $\mathbb E\left(\bm \ell_{\bm \delta},\bm \ell_{\bf Q}\right)=(\bf 0,\bf 0)$ and $\mathbb E\left\|\bm \ell_{\bm \delta}+\bm \ell_{\bf Q}\right\|^2<\infty$.

Then we may write from \eqref{TS} 
\[
T_{n,w}=\int_{\mathbb R^p} \Xi^2_n(\bm t) w(\bm t) {\rm{d}} \bm t, 
\quad
\Xi_n(\bm t)=\frac{1}{\sqrt{n}} \sum_{j=1}^n \left(\cos \bm t^\top \bm Y_j+\sin \bm t^\top \bm Y_j-\phi_{\alpha_0}(\bm t)\right).
\]

It also follows that $\|\Xi_{n}-\Xi_{n,0}\|^2  \stk 0$, where  the approximating process $\Xi_{n,0}(\cdot)$ admits the i.i.d.\ representation 
\[
\Xi_{n,0}(\bm t)=\frac{1}{\sqrt{n}} \sum_{j=1}^n L_{{\bm \ell_{\bm {\delta}}},{\bm {\ell_{\bf {Q}}}}}(\bm X_j),
\]
upon which the central limit theorem applies, and which together with a subsequent application of the continuous mapping theorem entails   
\[
T_{n,w} \verk \int_{\mathbb R^p} \Xi^2_0(\bm t) w(\bm t) {\rm{d}} \bm t:=T_w, 
\]
where $\Xi_0(\cdot)$ is a zero-mean Gaussian process with covariance kernel, say, $K(\bm s,\bm t)$. In turn the law of $T_w$ is that of $\sum_{j=1}^\infty \lambda_j {\cal{N}}^2_j$, where $({\cal{N}}_j, \ j\geq 1)$ are i.i.d. standard normal random variables. The covariance kernel $K(\cdot,\cdot)$ of the limit process $\Xi_0(\cdot)$ enters the distribution of $T_w$ via the eigenvalues $\lambda_1\leq \lambda_2\leq \cdots $ and corresponding eigenfunctions $f_1,f_2,\cdots$, through the integral equation
\[
\lambda f(\bm s)=\int_{\mathbb R^p} K(\bm t,\bm s) f(\bm t) w(t){\rm{d}}\bm t.  
\]
In this connection the maximum likelihood estimators defined  by \eqref{MLE}  satisfy   certain equivariance/invariance properties (refer to Section~S2 of the online supplement),  and as a consequence the resulting test statistic is affine invariant. In this case we may set $(\bm \delta_0,\bf Q_0)$ equal to the zero vector and identity matrix, respectively,  thus rendering the limit null distribution free of true parameter values; see \cite{EH20} and \cite{MNT}.

Moreover the standing  assumptions imply the strong consistency of the new test under fixed alternatives.   

\begin{prop} \label{cons}
Suppose that under the given law of $\bm X$ the estimators of the parameters $\bm \delta$ and $\bf Q$ satisfy,  $(\widehat {\bm  \delta}_n,\widehat {\bf Q}_n) \to (\bm \delta_X, {\bf {Q}}_{\bm X}) \in \mathbb R^p\times {\cal{M}}^p$, a.s. as $n\to\infty$. Then we have \begin{equation} \label{limstat}
\frac{ T_{n,w}}{n}  \rightarrow \int_{\mathbb R^p} \left\lvert\varphi\left({\bf Q}_{\bm X}^{-1/2} \bm t\right) -e^{{\rm{i}} \bm t^\top {\bf Q}_{\bm X}^{-1/2} {\bm \delta}_{\bm X} }\phi_{\alpha_0}(\bm t)\right\rvert^2w(\bm t){\rm{d}}\bm t:=  {\cal{T}}_{w}, 
\end{equation}
a.s. as  $n\rightarrow \infty$. 
\end{prop}
\begin{proof}
Recall from \eqref{TS} that
\begin{equation} \label{tsn} 
\frac{ T_{n,w}}{n} =\int_{\mathbb R^p} \lvert\phi_n(\bm t)-\phi_{\alpha_0}(\bm t)\rvert^2 w(\bm t){\rm{d}}\bm t. \end{equation}
Now the strong uniform consistency of the empirical CF in bounded intervals \citep[see][]{csorgHo1981multivariate} entails
\[
\varphi_n(t) \rightarrow e^{-{\rm{i}}\bm t^\top {\bf {Q}}_{\bm X}^{-1/2}{\bm \delta}_{\bm X}} \varphi\left({\bf {Q}}_{\bm X}^{-1/2}\bm t\right),
\]
a.s. as  $n\rightarrow \infty$. Consequently, since  $\lvert\phi_n(\bm t)-\phi_{\alpha_0}(\bm t)\rvert^2\leq 4$, an application of Lebesgue's theorem of dominated convergence on \eqref{tsn} yields \eqref{limstat}.
\end{proof}

Since $w>0$, ${\cal{T}}_{w}$ is positive unless 
\[
\varphi\left( {\bf {Q}}_{\bm X}^{-1/2}\bm t\right)= e^{{\rm{i}} {\bf t}^\top {\bf{Q}}_{\bm X}^{-1/2}{\bm \delta}_{\bm X}} \phi_{\alpha_0}(\bm t),
\]
identically in $\bm t$, which is equivalent to  
% \begin{equation} \label{cf1a} 
$\varphi(\bm t)=\varphi_{\alpha_0}({\bm{t}};{\bm{\delta}}_{\bm X},{\bf{Q}}_{\bm X})$, $\bm t \in \mathbb R^p$,
% \end{equation}
i.e.\ ${\cal{T}}_{w}>0$ unless the CF of $\bm X$ coincides with the CF of a SP law with $\alpha=\alpha_0$ and $({\bm {\delta}}, {\bf {Q}})=({\bm {\delta}}_{\bm X},{\bf {Q}}_{\bm X})$, and thus by the uniqueness of CFs, the test which rejects the null hypothesis   ${\cal{H}}_0$ in \eqref{null} for large values of $T_{n,w}$ is consistent against each fixed  alternative distribution.      

The above limit results have been developed in a series of papers, both in the current setting as well as in related settings,  and with varying conditions on the weight function and the family of distributions under test; see for instance \cite{HW97}, \cite{GHK}, \cite{MNT},  \cite{HRA}, and \cite{EH20}. In this regard, the solution of the above integral equation, and thus the approximation of the the limit null distribution of $T_{n,w}$, is extremely complicated, and in fact  constitutes a research problem in itself.  This line of research has been followed by a few authors. We refer to \cite{matsui2008}, \cite{HRM}, \cite{meintanis2023} and \cite{ebner2023}. In these works, the infinite sum distribution of $T_w$ is approximated by a corresponding finite sum employing numerically computed eigenvalues and then large-sample critical points for $T_{n,w}$ are found by Monte Carlo. It should be noted that such approximation is specific in several problem-parameters: the distribution under test, the type of  estimators of the distributional parameters, and the weight function employed, and thus have to be performed on a case-to-case basis. A different more heuristic approach is moment-matching between the first few moments of $T_w$ (computed numerically) and a known distribution,  like the gamma distribution or one from a Pearson family of distributions; see \cite{henze1990,henze1997} and \cite{pfister2018}, while yet another, Satterthwaite-type, approximation is studied by \cite{lindsay2008}.

%Matsui and Takemura (2008), TEST volume 17, pages 546–566
%Henze (1997) Statistics & Probability Letters 35,  203-213
%Henze (1990) Metrika 37:7-18
%Lindsay et al. (2008), Ann.Statist. 36, 983-1006
%Ebner and Henze, Statistical Papers (2023) 64:739–752
%Pfister et al. JRS-B: Statistical Methodology, Volume 80, Issue 1, January 2018, Pages 5–31
%Meintanis et al. (2023). Metrika 86, 723-751 

The validity and usefulness of the above approximation methods notwithstanding, we hereby favor Monte Carlo simulation and  bootstrap resampling for the computation of critical points and for test implementation, not only because these are large-sample approximations performed mostly in the univariate setting and thus inappropriate for small sample size $n$ and/or dimension $p>1$, but more importantly because, in the case of GARCH-type observations considered in the next section, the finite-sample counterpart of this distribution may even involve true parameter values; see for instance \cite{HJM}.  %In view of these considerations here we emphasize  the pragmatic aspects of the test which rely on bootstrap resampling for the computation of critical points and for actual test implementation.  

\section{The case of the stable-GARCH model} \label{sec_5}
Assume that observations $(\bm X_{j}, \ j=1,...,n)$ arise from a  multivariate GARCH model defined by
\begin{equation}
\label{garch}
\bm X_{j}= {\bf Q}^{1/2}_{j} \bm \varepsilon_{j}, \quad j=1,\ldots,n,
\end{equation}
where $(\bm \varepsilon_{j}, j=1,...,n)$ is a sequence of i.i.d.\ $p$-dimensional random variables with mean zero and identity dispersion matrix, and ${\bf Q}_{j}:={\bf Q}(\bm X_j\vert\mathbb I_{j-1})$,
with $\mathbb I_{j}$ denoting the information available at time $j$, is a symmetric positive definite matrix of dimension $(p\times p)$. We wish to test the null hypothesis stated in \eqref{null} for the innovations $\bm \varepsilon_j$ figuring in model \eqref{garch}. Note that in view of \eqref{sr},  ${\bf Q}_j$ may be interpreted  as the conditional covariance matrix of the corresponding latent Gaussian vector $\bm N$. 

As GARCH model we adopt the constant conditional correlation (CCC) specification \citep[see, for instance,][]{silvennoinen2009} %,francq2019}
\begin{equation}\begin{split}\label{CCC}
{\bf Q}_j(\bm \vartheta, \kappa_x,\kappa_q)&=  {\mat D}_j  {\mat R}  {\mat D}_j, \\
\bm q_{j}& =  (q^2_{1,j},\ldots,q^2_{p,j})^\top= \bm \mu+\sum_{k=1}^{\kappa_x}  {\mat{A}}_k  \bm X^{(2)}_{j-k}+\sum_{k=1}^{\kappa_q}  {\mat{B}}_k  \bm q_{j-k},
\end{split}\end{equation}
where $\bm X^{(2)}_{j}=(X^2_{1,j},\ldots,X^2_{p,j})^\top$, ${\mat{D}}_j={\rm{diag}}(q_{1,j},\ldots, q_{p,j})$, $\bm \mu$ is an $(p\times 1)$ vector of positive elements, $\mat R$ is a correlation matrix, and ${\mat{A}}_k$ and ${\mat{B}}_k$ are $(p\times p)$ matrices with nonnegative elements.
Note that the vector $\bm \vartheta$ which incorporates all parameters figuring in (\ref{CCC}) is of dimension $p+p^2(\kappa_x+\kappa_q)+(p(p-1)/2)$.

Now let $\bm  \Upsilon_{j,\kappa_x,\kappa_q}=(\bm X^\top_{j-1},\ldots,\bm X^\top_{j-\kappa_x},\bm q^\top_{j-1},\ldots, \bm q^\top_{j-\kappa_q})^\top$ be a vector of dimension $p\times (\kappa_x+\kappa_q)$ associated with model (\ref{garch})--(\ref{CCC}). Since the innovations are unobserved, any test  should be applied on the corresponding residuals,
 \begin{equation} \label{residual}
 \bm {\widetilde \varepsilon}_{j}={\bf {\widetilde Q}}_j^{-1/2}(\bm {\widehat \vartheta}_n) \bm X_{j}, \quad j=1,...,n,
 \end{equation}
where  $\bm {\widehat\vartheta}_n$ is a consistent estimator of $\bm \vartheta$ that uses $\bm  \Upsilon_{j,\kappa_x,\kappa_q}$ as input data. We also employ initial values in order to start the estimation process. As an estimator of ${\vec\vartheta}$ we use the equation-by-equation (EbE) estimators
proposed by \citet{francq2014}. The reader is referred to Section~S3 of the online supplement for more details on the EbE estimator. 

%Note however that  ${\bf Q}_j(\bm \vartheta)$ depends on $\{\bm X_k, \, -\infty<k \leq j-1\}$, whereas  we only observe $\bm X_1, \ldots, \bm X_n$. So for calculating the residuals we consider ${\bf {\widetilde Q}}_j (\bm {\widehat\vartheta}_n)$ defined  as ${\bf Q}_j (\bm {\widehat\vartheta}_n)$, and also employ initial values  $\bm {\widetilde{\Upsilon}}_1:=(\bm {\widetilde X}^\top_0,\ldots,\bm {\widetilde X}^\top_{1-\kappa_x},\bm {\widetilde q}^\top_{0},\ldots,\bm {\widetilde q}^\top_{1-\kappa_q})^\top$, in order to start the estimation process. 

%\textcolor{blue}{(Simos, I've added the following paragraph on EbE estimation.
%Maybe we could move it to the supplement, or just mention that we use EbE with a reference?)}

% Under the specification in \eqref{CCC},

% Assume that the $k$th component of $D_j$ is parametrised by ${\vec\theta}^{(k)}_0$,
% so that $q_{k,j}({\vec\theta}^{(k)}_0)=q_k(\vec X_{j-1},\vec X_{j-2},\ldots;{\vec\theta}^{(k)}_0)$
% for some positive function $q_k$, e.g.\
% \[
%     q_{k,j}^2 = \mu_k + a_{k,k} X_{k,j-1}^2 + b_{k,k} q_{k,j-1}^2.
% \]

% Each ${\vec\theta}^{(k)}_0$ is estimated using MLE:
% \[
%     \hat{\vec\theta}^{(k)}
%     = \argmax_{\bm\theta^{(k)}}
%             \; \prod_{j=1}^n \frac{1}{q_{kj}\left(\bm\theta^{(k)}\right)}
%                         f_{\alpha_0}\left(\frac{X_{kj}}{q_{kj}\left(\bm\theta^{(k)}\right)}\right),
%     \qquad
%     k=1,\ldots,p.
% \]

% Francq and Zakoian (2014) studied asymptotic and finite sample properties
% of EbE estimators.

%\section{Asymptotics} \label{sec_6}

\section{Numerical study} \label{sec_6}
We now turn to a simulation study to demonstrate the performance of the new tests
under various alternatives, specifically:
\begin{enumerate}[label=(\roman*)]
	\item elliptically symmetric SP distributions, denoted by $S_{\alpha}$;
	\item Student $t$-distributions with $\nu$ degrees of freedom, denoted by $t_\nu$;
	\item skew normal distributions with skewness parameter $\nu$, denoted by $SN_\nu$;
	\item skew Cauchy distributions with skewness parameter $\nu$, denoted by $SC_\nu$;
\end{enumerate}
The cases $p=4,6$ are considered for the dimension of the distribution.
Throughout, we choose ${\vec\delta}={\vec 0}$ and ${\mat Q}={\mat I}_p$,
where ${\mat I}_p$ denotes the $p\times p$ identity matrix. Recall in this connection that our  method requires the computation of the density of a spherical SP  law. Details on this computation are included in Section~S4 of the online supplement.

\subsection{Monte Carlo results for i.i.d.\ data}
We focus on our test criterion in \eqref{STS} that employs $\phi_{\alpha_0}(\cdot)$ in the weight function, say $T_r$ for simplicity,
and we consider several values of the tuning parameter $r>0$.
Corresponding results for the test in \eqref{TS2} that employs the Kotz-type weight function are included in Tables S8 and S9 of the online supplement. 

For comparison, we include results obtained when using the test of
\citet*{MNT} using the Gaussian weight function $\exp(-\|\bm t\|^2)$.
% As in the original paper, we use the Gaussian weight function
% $\exp(-\gamma\|\bm t\|^2)$ and, based on the numerical results reported
% in that paper, we choose $\gamma=1$ throughout.
The test, denoted by $M_a$ in the tables, depends on a tuning parameter denoted by
$a$ for which we consider the choices $a=4,6,10$ and $a=15$.
As pointed out by the authors, the number of operations required to compute
their test statistic is in the order of $n^{2a}$ (for integer-valued $a$) and becomes
time-consuming for larger values of $n$ and $a$. We therefore use the Monte Carlo
approach suggested by the authors to approximate the value of the test statistic
(using 1\,000 replications for each approximation); see p.~180 of \citet{MNT}.
When testing for multivariate normality (i.e.\ $\mathcal{H}_0$ with
$\alpha=2$), we consider the test of \citet*{HJM}, denoted by HJM in the
tables, with weight function $\exp(-1.5\|\bm t\|^2)$, which yielded good results
in the original paper.

The rejection percentages of the tests are shown in
Tables~\ref{tab:alpha2_vs_stable_estimate}--\ref{tab:alpha1_vs_stable_estimate}.
All simulations results are based on 1,000 independent Monte Carlo iterations and
a significance level of 10\% is used throughout.

We first consider the case where we test $\mathcal H_0$ with $\alpha=2$,
i.e.\ when testing for departures from multivariate normality.
Table~\ref{tab:alpha2_vs_stable_estimate}
% (which contains selected results
% available in Table~\ref{tab:alpha2_vs_stable_estimateFULL} of the online supplement)
shows that, when heavier-tailed symmetric Paretian alternatives are considered, the
newly proposed tests based on $T_r$ are more powerful than the existing tests $M_a$ of
\citet{MNT}, and have power slightly lower but comparable to the test HJM of \citet{HJM}.
In light of the above-mentioned computational complexity of the existing tests,
this gain in power makes the new tests attractive for implementation in practice.
Moreover, the favorable power is visible for all the considered choices of the tuning parameter $r$,
the choice of which, as opposed to $a$ in $M_a$, has no significant impact on computational complexity.
Finally, as is expected, the power of the new tests increase as the sample size
is increased.
Similar conclusions can also be made in the case of elliptically symmetric $t$ alternatives
(see Table~S2 of the online supplement).

\begin{table}[tbp]
	\centering\footnotesize
	\begin{threeparttable}
	\caption{Percentage of rejection of $\mathcal{H}_0$ with $\alpha=2$
					 against stable alternatives.
					 Tests done at the 10\% significance level.
                      Additional cases are shown Table~S1
                      of the online supplement.}
	\label{tab:alpha2_vs_stable_estimate}
	\begin{tabular}{cccHHHHHHrrrrHHH|HrrrrHHr}
		\toprule
		$p$ & $n$ & $\mathcal{H}_1$ & $T_{0.001}$ & $T_{0.005}$ & $T_{0.01}$ & $T_{0.05}$ & $T_{0.1}$
							& $T_{0.5}$ & $T_{1}$ & $T_{2}$ & $T_{5}$ & $T_{10}$ & $T_{20}$ & $T_{50}$ & $T_{100}$
						& $M_{2}$ & $M_{4}$ & $M_{6}$ & $M_{10}$ & $M_{15}$ & $M_{20}$ & $M_{25}$ & HJM \\
		\midrule
\multicolumn{1}{r}{4} & 50    &  $S_{1.9}$  & 3.1   & 13.9  & 15.2  & 17.8  & 21.3  & 34.1  & 41.7  & 52.5  & 59.9  & 60.9  & 61.8  & 62.0  & 62.0  & 6.9   & 12.3  & 15.8  & 15.1  & 13.4  & 14.2  & 15.0  & 62.6 \\
      &       &  $S_{1.8}$  & 3.0   & 17.6  & 21.8  & 32.4  & 37.9  & 57.8  & 68.1  & 76.6  & 82.9  & 82.8  & 83.2  & 83.2  & 83.2  & 5.4   & 20.0  & 24.2  & 27.7  & 26.1  & 24.5  & 21.7  & 84.0 \\
      &       &  $S_{1.7}$  & 4.9   & 26.2  & 30.1  & 48.7  & 60.0  & 82.0  & 87.9  & 92.1  & 93.0  & 93.2  & 93.5  & 93.5  & 93.4  & 9.0   & 37.1  & 43.8  & 44.3  & 41.6  & 35.0  & 29.4  & 94.6 \\
      % &       &  $S_{1.5}$  & 13.1  & 42.8  & 52.6  & 75.9  & 88.0  & 98.5  & 99.1  & 99.4  & 99.6  & 99.4  & 99.4  & 99.3  & 99.3  & 28.4  & 74.8  & 80.3  & 81.8  & 77.0  & 69.3  & 62.6  & 99.8 \\
% \cmidrule{2-24}      & 75    &  $S_{1.9}$  & 4.2   & 15.1  & 14.9  & 22.8  & 28.0  & 42.0  & 50.7  & 61.4  & 68.8  & 71.0  & 71.5  & 71.4  & 71.4  & 5.6   & 12.0  & 15.9  & 15.2  & 16.6  & 14.2  & 12.6  &  \\
%       &       &  $S_{1.8}$  & 5.7   & 20.8  & 22.9  & 43.2  & 52.5  & 74.0  & 82.3  & 87.7  & 90.8  & 91.0  & 91.2  & 91.2  & 91.2  & 4.3   & 22.4  & 32.0  & 32.3  & 32.4  & 31.3  & 23.8  &  \\
%       &       &  $S_{1.7}$  & 8.4   & 31.0  & 35.8  & 64.4  & 77.2  & 92.3  & 95.8  & 97.4  & 98.2  & 98.1  & 98.2  & 97.9  & 97.8  & 9.4   & 44.6  & 55.5  & 54.0  & 50.8  & 47.3  & 39.5  &  \\
%       &       &  $S_{1.5}$  & 22.4  & 60.5  & 69.2  & 94.9  & 98.1  & 99.8  & 99.8  & 99.8  & 99.8  & 99.8  & 99.8  & 99.8  & 99.8  & 47.5  & 90.0  & 94.0  & 90.9  & 88.0  & 82.2  & 80.0  &  \\
% \cmidrule{2-24}   
& 100   &  $S_{1.9}$  & 6.2   & 14.9  & 17.1  & 28.9  & 33.6  & 50.6  & 60.1  & 68.7  & 76.7  & 78.4  & 78.7  & 79.1  & 79.2  & 5.1   & 13.9  & 17.0  & 18.6  & 17.0  & 15.8  & 16.8  & 81.4 \\
      &       &  $S_{1.8}$  & 8.2   & 22.9  & 28.1  & 48.8  & 59.4  & 82.2  & 88.7  & 92.3  & 94.6  & 95.1  & 94.8  & 94.9  & 94.9  & 6.0   & 25.1  & 35.7  & 37.2  & 34.6  & 27.4  & 28.6  & 96.7 \\
      &       &  $S_{1.7}$  & 11.0  & 35.8  & 41.6  & 76.0  & 85.9  & 97.4  & 98.4  & 99.1  & 99.3  & 99.6  & 99.6  & 99.6  & 99.6  & 12.9  & 56.9  & 67.3  & 66.0  & 61.4  & 54.5  & 49.9  & 99.7 \\
\midrule
\multicolumn{1}{r}{6} & 50    &  $S_{1.9}$  & 1.9   & 2.4   & 4.0   & 14.6  & 15.3  & 28.3  & 36.9  & 49.4  & 61.5  & 63.7  & 65.0  & 65.8  & 66.0  & 3.4   & 9.4   & 13.5  & 9.1   & 10.5  & 12.1  & 11.0  & 66.9 \\
      &       &  $S_{1.8}$  & 0.3   & 2.7   & 7.3   & 27.0  & 31.7  & 61.0  & 74.6  & 83.8  & 88.7  & 90.0  & 90.1  & 90.6  & 90.8  & 2.7   & 14.8  & 16.7  & 15.1  & 13.1  & 11.4  & 14.0  & 92.3 \\
      &       &  $S_{1.7}$  & 0.3   & 3.9   & 9.9   & 41.8  & 54.4  & 82.8  & 90.8  & 95.4  & 97.2  & 97.4  & 97.5  & 97.6  & 97.6  & 2.9   & 24.8  & 27.4  & 19.9  & 16.4  & 14.1  & 11.1  & 97.0 \\
      % &       &  $S_{1.5}$  & 1.9   & 19.2  & 36.1  & 77.8  & 88.2  & 99.2  & 99.8  & 99.9  & 100.0 & 100.0 & 100.0 & 100.0 & 100.0 & 17.4  & 65.9  & 68.6  & 47.8  & 38.6  & 32.4  & 27.8  & 99.8 \\
% \cmidrule{2-24}      & 75    &  $S_{1.9}$  & 0.7   & 2.3   & 6.0   & 18.4  & 22.3  & 40.0  & 50.7  & 65.1  & 75.8  & 78.3  & 79.2  & 79.7  & 79.9  & 3.5   & 12.1  & 12.2  & 8.6   & 9.5   & 12.7  & 11.2  &  \\
%       &       &  $S_{1.8}$  & 0.6   & 3.6   & 10.4  & 32.6  & 44.3  & 74.4  & 84.0  & 91.6  & 95.4  & 96.0  & 96.4  & 96.3  & 96.3  & 3.5   & 18.7  & 20.4  & 14.8  & 14.0  & 14.4  & 13.4  &  \\
%       &       &  $S_{1.7}$  & 0.9   & 6.6   & 17.9  & 53.2  & 68.0  & 93.4  & 97.2  & 98.8  & 99.3  & 99.6  & 99.7  & 99.7  & 99.7  & 7.0   & 35.9  & 37.0  & 22.7  & 18.0  & 17.1  & 14.2  &  \\
%       &       &  $S_{1.5}$  & 2.9   & 24.1  & 49.7  & 90.0  & 96.9  & 100.0 & 100.0 & 100.0 & 100.0 & 100.0 & 100.0 & 100.0 & 100.0 & 34.7  & 83.4  & 83.8  & 66.0  & 51.1  & 40.4  & 31.4  &  \\
% \cmidrule{2-24}     
& 100   &  $S_{1.9}$  & 0.7   & 3.0   & 8.3   & 19.9  & 25.1  & 51.3  & 64.4  & 75.7  & 82.9  & 86.0  & 86.8  & 87.3  & 87.5  & 2.4   & 10.1  & 12.1  & 11.6  & 12.1  & 12.3  & 11.7  & 90.6 \\
      &       &  $S_{1.8}$  & 0.3   & 4.6   & 14.5  & 38.8  & 55.5  & 86.8  & 92.8  & 96.9  & 98.8  & 99.4  & 99.5  & 99.5  & 99.4  & 1.3   & 17.4  & 23.9  & 16.8  & 13.6  & 12.0  & 13.1  & 99.6 \\
      &       &  $S_{1.7}$  & 0.7   & 11.1  & 26.6  & 62.2  & 80.1  & 98.1  & 99.7  & 99.9  & 100.0 & 100.0 & 100.0 & 100.0 & 100.0 & 5.3   & 36.8  & 43.0  & 32.3  & 22.5  & 19.8  & 16.8  & 99.9 \\
		\bottomrule
	\end{tabular}
	\end{threeparttable}
\end{table}

Considering skew normal alternatives, even more favorable behavior can be observed
in the results shown in Table~\ref{tab:alpha2_vs_skew_estimate}.
% (full results in Table~\ref{tab:alpha2_vs_skew_estimateFULL} of the supplement).
Notice that the tests $M_a$ and HJM seem to have very low power against
skew normal alternatives, which is not the case with the newly proposed tests.

\begin{table}[tbp]
	\centering\footnotesize
	\begin{threeparttable}
	\caption{Percentage of rejection of $\mathcal{H}_0$ with $\alpha=2$
			 against skew normal alternatives.
			 Tests done at the 10\% significance level.
              Also see Table~S3 of the online supplement.}
	\label{tab:alpha2_vs_skew_estimate}
	\begin{tabular}{cccHHHHHHrrrrHHH|HrrrrHHr}
		\toprule
		$p$ & $n$ & $\mathcal{H}_1$ & $T_{0.001}$ & $T_{0.005}$ & $T_{0.01}$ & $T_{0.05}$ & $T_{0.1}$
							& $T_{0.5}$ & $T_{1}$ & $T_{2}$ & $T_{5}$ & $T_{10}$ & $T_{20}$ & $T_{50}$ & $T_{100}$
						& $M_{2}$ & $M_{4}$ & $M_{6}$ & $M_{10}$ & $M_{15}$ & $M_{20}$ & $M_{25}$ & HJM \\
         \midrule
\multicolumn{1}{r}{4} & \multicolumn{1}{r}{50} &   $SN_{1.5}$  & 10.3  & 10.1  & 9.7   & 11.0  & 12.7  & 21.8  & 29.1  & 32.0  & 30.6  & 29.5  & 28.7  & 28.2  & 28.0  & 12.9  & 12.9  & 12.1  & 13.5  & 9.6   & 12.3  & 11.0  & 16.3 \\
          &       &  $SN_{1.0}$  & 10.3  & 11.2  & 11.3  & 10.6  & 12.7  & 16.9  & 20.2  & 22.0  & 24.1  & 23.3  & 23.3  & 23.4  & 23.3  & 11.6  & 11.5  & 11.9  & 11.8  & 9.2   & 12.7  & 9.1   & 14.4 \\
          &       &  $SN_{0.5}$  & 9.8   & 9.8   & 9.6   & 9.7   & 10.3  & 11.4  & 11.9  & 13.6  & 15.0  & 14.1  & 14.3  & 13.9  & 14.1  & 9.5   & 7.8   & 7.1   & 10.2  & 12.1  & 10.5  & 9.0   & 12.6 \\
% \cmidrule{2-24}		    & \multicolumn{1}{r}{75} &  $SN_{2.0}$  & 9.2   & 10.0  & 10.3  & 10.7  & 14.6  & 27.1  & 35.2  & 42.9  & 47.0  & 46.9  & 46.0  & 44.7  & 43.6  & 12.5  & 13.0  & 11.4  & 12.5  & 11.7  & 9.2   & 11.6  &  \\
% 		    &       &  $SN_{1.5}$  & 8.7   & 9.8   & 10.2  & 8.8   & 14.0  & 21.7  & 27.6  & 32.8  & 37.8  & 37.5  & 36.9  & 36.0  & 35.4  & 11.5  & 10.9  & 10.0  & 12.6  & 15.0  & 13.9  & 10.0  &  \\
% 		    &       &  $SN_{1.0}$  & 9.9   & 9.3   & 9.3   & 8.4   & 11.8  & 16.1  & 17.5  & 21.2  & 23.4  & 23.9  & 24.4  & 23.9  & 23.6  & 11.2  & 10.6  & 9.5   & 10.0  & 9.2   & 8.8   & 13.9  &  \\
% 		    &       &  $SN_{0.5}$  & 10.7  & 10.8  & 11.0  & 10.3  & 11.9  & 11.9  & 11.3  & 12.7  & 13.7  & 13.6  & 14.1  & 13.3  & 13.1  & 9.9   & 10.7  & 9.8   & 12.1  & 11.4  & 9.1   & 12.2  &  \\
% \cmidrule{2-24}   
& \multicolumn{1}{r}{100} &  $SN_{1.5}$  & 7.5   & 10.4  & 10.7  & 15.2  & 19.8  & 36.3  & 46.1  & 57.0  & 60.4  & 59.2  & 58.3  & 56.7  & 56.0  & 14.1  & 15.5  & 14.7  & 12.4  & 13.4  & 11.5  & 13.9  & 25.7 \\
      &       &  $SN_{1.0}$  & 8.9   & 9.5   & 9.2   & 11.0  & 14.1  & 24.1  & 29.6  & 37.0  & 41.0  & 41.9  & 41.2  & 41.4  & 41.4  & 9.6   & 14.2  & 11.5  & 13.7  & 12.5  & 10.5  & 6.4   & 20.4 \\
      &       &  $SN_{0.5}$  & 9.7   & 11.3  & 10.1  & 9.0   & 8.8   & 10.7  & 12.8  & 14.8  & 18.5  & 17.9  & 17.9  & 17.6  & 17.4  & 9.7   & 10.1  & 8.3   & 10.4  & 10.7  & 10.1  & 7.7   & 14.0 \\
\midrule
% \multicolumn{1}{r}{6} 
& \multicolumn{1}{r}{50} &   $SN_{1.5}$  & 6.9   & 6.8   & 7.7   & 14.4  & 13.9  & 15.4  & 18.5  & 21.8  & 20.7  & 19.8  & 18.5  & 18.8  & 18.6  & 11.0  & 11.4  & 9.8   & 7.9   & 11.9  & 12.8  & 8.7   & 12.7 \\
          &       &  $SN_{1.0}$  & 10.1  & 10.0  & 9.6   & 10.5  & 10.3  & 12.9  & 13.2  & 18.2  & 19.7  & 19.0  & 19.2  & 19.0  & 19.4  & 10.4  & 10.3  & 7.8   & 9.6   & 10.1  & 9.7   & 11.2  & 13.7 \\
          &       &  $SN_{0.5}$  & 9.6   & 9.2   & 10.8  & 9.7   & 9.4   & 9.2   & 8.9   & 9.6   & 9.9   & 11.1  & 10.6  & 10.6  & 9.9   & 11.1  & 10.5  & 9.9   & 9.1   & 12.0  & 10.8  & 10.9  & 10.1 \\
% \cmidrule{2-24}		    & \multicolumn{1}{r}{75} &  $SN_{2.0}$  & 8.9   & 9.0   & 10.1  & 11.4  & 10.9  & 17.8  & 23.2  & 29.0  & 30.0  & 29.5  & 29.9  & 29.1  & 28.9  & 12.2  & 9.9   & 10.0  & 9.8   & 11.0  & 8.6   & 8.4   &  \\
% 		    &       &  $SN_{1.5}$  & 8.8   & 8.8   & 10.2  & 10.1  & 10.1  & 15.7  & 19.9  & 23.6  & 25.0  & 25.2  & 25.2  & 24.7  & 24.3  & 11.4  & 9.9   & 9.7   & 11.1  & 9.5   & 6.5   & 8.6   &  \\
% 		    &       &  $SN_{1.0}$  & 8.8   & 8.7   & 9.2   & 9.2   & 9.3   & 12.3  & 15.1  & 17.7  & 18.2  & 18.0  & 17.8  & 17.2  & 17.1  & 10.6  & 10.6  & 11.0  & 12.0  & 10.1  & 9.9   & 10.6  &  \\
% 		    &       &  $SN_{0.5}$  & 9.6   & 8.9   & 8.8   & 8.5   & 8.2   & 10.0  & 10.8  & 11.2  & 11.1  & 11.5  & 12.3  & 11.6  & 11.6  & 11.4  & 10.0  & 9.8   & 11.1  & 11.1  & 11.2  & 10.1  &  \\
% \cmidrule{2-24}     
& \multicolumn{1}{r}{100} &  $SN_{1.5}$  & 7.7   & 8.9   & 9.8   & 11.8  & 11.6  & 19.3  & 27.4  & 32.2  & 37.2  & 37.2  & 36.5  & 35.2  & 34.7  & 13.6  & 13.6  & 12.0  & 7.8   & 9.8   & 11.7  & 12.0  & 16.3 \\
      &       &  $SN_{1.0}$  & 8.4   & 7.9   & 8.5   & 9.4   & 10.1  & 12.9  & 16.2  & 21.1  & 22.6  & 23.3  & 24.2  & 23.8  & 24.6  & 8.2   & 12.1  & 10.2  & 11.2  & 10.2  & 9.2   & 8.8   & 10.9 \\
      &       &  $SN_{0.5}$  & 7.7   & 8.2   & 7.9   & 9.7   & 10.8  & 11.4  & 13.0  & 14.0  & 13.4  & 14.3  & 13.1  & 13.5  & 13.1  & 9.6   & 9.7   & 9.3   & 10.5  & 12.3  & 11.3  & 11.0  & 10.5 \\
		\bottomrule
	\end{tabular}
	\end{threeparttable}
\end{table}

We now shift our attention to the case of testing $\mathcal H_0$ with $\alpha=1.8$.
Table~\ref{tab:alpha1.8_vs_stable_estimate}
% (full results in Table~\ref{tab:alpha1.8_vs_stable_estimate_FULL})
shows that, compared to the existing tests,
the new tests are quite powerful against heavier-tailed alternatives, i.e.\ alternatives
with stability index less than $1.8$. Despite the evident dependence of the performance of the
new tests on the tuning parameter $r$, we note that the power is very competitive
to the existing tests in most cases, and significantly outperforms the existing
tests for most choices of $r$.
For lighter-tailed alternative distributions, the new tests exhibit some
under-rejection for certain choices of $r$. Nevertheless, the problem
seems to disappear as the sample size is increased.

\begin{table}[tbp]
	\centering\footnotesize
	\begin{threeparttable}
	\caption{Percentage of rejection of $\mathcal{H}_0$ with $\alpha=1.8$
					 against stable alternatives.
					 Tests done at the 10\% significance level.
            Also see Table~S4 of the
            online supplement.}
	\label{tab:alpha1.8_vs_stable_estimate}
	\begin{tabular}{cccHHHHHHrrrrHHH|HrrrrHH}
		\toprule
		$p$ & $n$ & $\mathcal{H}_1$ & $T_{0.001}$ & $T_{0.005}$ & $T_{0.01}$ & $T_{0.05}$ & $T_{0.1}$
							& $T_{0.5}$ & $T_{1}$ & $T_{2}$ & $T_{5}$ & $T_{10}$ & $T_{20}$ & $T_{50}$ & $T_{100}$
						& $M_{2}$ & $M_{4}$ & $M_{6}$ & $M_{10}$ & $M_{15}$ & $M_{20}$ & $M_{25}$ \\
		\midrule
  \multicolumn{1}{r}{4} & \multicolumn{1}{r}{50} &  $S_{2.0}$  & 31.4  & 19.3  & 10.6  & 10.4  & 10.6  & 11.8  & 9.1   & 3.1   & 0.1   & 0.0   & 0.0   & 0.0   & 0.0   & 14.7  & 12.3  & 12.2  & 12.1  & 12.5  & 12.4  & 10.6 \\
      &       &  $S_{1.9}$  & 19.5  & 13.2  & 9.8   & 10.3  & 9.6   & 9.6   & 9.1   & 6.0   & 2.0   & 1.8   & 1.9   & 2.0   & 3.0   & 10.8  & 10.0  & 9.3   & 11.1  & 8.3   & 10.1  & 9.4 \\
      &       &  $S_{1.7}$  & 6.4   & 8.4   & 10.1  & 10.9  & 10.5  & 11.3  & 13.2  & 18.8  & 31.2  & 32.4  & 31.1  & 27.5  & 24.6  & 11.9  & 12.7  & 11.8  & 11.5  & 9.4   & 10.9  & 11.9 \\
      &       &  $S_{1.5}$  & 3.5   & 9.2   & 12.6  & 13.4  & 13.4  & 22.9  & 33.2  & 54.7  & 79.0  & 78.8  & 75.8  & 69.2  & 62.8  & 29.1  & 30.0  & 26.4  & 16.7  & 12.7  & 11.1  & 11.3 \\
& \multicolumn{1}{r}{100} &  $S_{2.0}$  & 36.8  & 10.8  & 10.5  & 10.4  & 10.3  & 13.4  & 12.6  & 5.3   & 1.0   & 30.3  & 12.1  & 0.0   & 0.0   & 18.0  & 16.6  & 15.8  & 15.8  & 12.0  & 10.3  & 12.3 \\
      &       &  $S_{1.9}$  & 22.7  & 9.8   & 9.6   & 9.1   & 9.8   & 10.4  & 9.9   & 5.8   & 1.4   & 3.0   & 1.0   & 0.7   & 0.9   & 12.8  & 11.1  & 10.2  & 11.4  & 11.5  & 10.2  & 10.5 \\
      &       &  $S_{1.7}$  & 5.8   & 10.5  & 11.0  & 11.2  & 10.9  & 13.4  & 17.6  & 23.1  & 36.6  & 40.0  & 38.5  & 33.2  & 30.0  & 13.0  & 13.7  & 13.7  & 12.4  & 11.7  & 10.5  & 11.8 \\
      &       &  $S_{1.5}$  & 4.0   & 11.9  & 13.0  & 15.1  & 19.6  & 42.3  & 58.2  & 76.6  & 91.8  & 93.6  & 91.3  & 85.9  & 79.4  & 44.4  & 48.5  & 47.0  & 33.4  & 21.6  & 19.4  & 16.2 \\
\midrule
\multicolumn{1}{r}{6} & \multicolumn{1}{r}{50} &  $S_{2.0}$  & 43.8  & 42.9  & 41.0  & 10.4  & 10.1  & 10.8  & 8.2   & 3.2   & 0.0   & 0.0   & 0.0   & 0.0   & 0.0   & 14.8  & 9.9   & 8.7   & 8.7   & 10.5  & 9.7   & 10.9 \\
      &       &  $S_{1.9}$  & 24.0  & 22.4  & 21.1  & 9.8   & 9.2   & 9.8   & 8.9   & 5.5   & 1.3   & 0.8   & 1.2   & 1.6   & 2.1   & 11.9  & 9.0   & 10.5  & 13.1  & 10.1  & 8.6   & 11.0 \\
      &       &  $S_{1.7}$  & 4.8   & 4.8   & 5.5   & 12.1  & 12.2  & 12.4  & 16.0  & 22.1  & 35.8  & 36.0  & 34.4  & 30.6  & 27.6  & 8.9   & 10.2  & 10.2  & 10.7  & 11.8  & 10.3  & 11.3 \\
      &       &  $S_{1.5}$  & 1.2   & 2.4   & 3.3   & 15.3  & 17.0  & 31.2  & 44.0  & 62.2  & 82.2  & 84.2  & 82.0  & 76.8  & 72.0  & 14.0  & 16.0  & 14.2  & 11.6  & 10.0  & 9.1   & 10.7 \\
  & \multicolumn{1}{r}{100} &  $S_{2.0}$  & 59.0  & 51.8  & 44.1  & 8.9   & 8.6   & 11.3  & 12.1  & 5.0   & 0.2   & 49.9  & 80.7  & 0.0   & 0.0   & 19.4  & 12.1  & 10.8  & 11.0  & 13.3  & 7.9   & 9.4 \\
      &       &  $S_{1.9}$  & 33.6  & 30.3  & 25.4  & 9.2   & 8.9   & 10.8  & 9.4   & 5.6   & 0.9   & 3.9   & 8.3   & 0.5   & 1.2   & 13.1  & 10.7  & 9.9   & 9.3   & 9.3   & 10.7  & 9.6 \\
      &       &  $S_{1.7}$  & 2.0   & 3.2   & 5.4   & 12.2  & 12.0  & 15.4  & 20.2  & 27.7  & 41.9  & 46.2  & 44.1  & 40.2  & 37.5  & 9.0   & 11.3  & 10.4  & 7.7   & 10.0  & 12.0  & 9.6 \\
      &       &  $S_{1.5}$  & 0.6   & 2.2   & 5.0   & 16.0  & 18.9  & 50.4  & 68.1  & 82.6  & 95.6  & 97.0  & 96.0  & 93.5  & 88.9  & 21.9  & 26.2  & 16.4  & 12.0  & 10.2  & 9.8   & 9.8 \\
\bottomrule
	\end{tabular}
	\end{threeparttable}
\end{table}

Another advantageous feature of tests based on $T_r$ is that the high power
is also visible in the higher-dimensional setting where $p=6$. In fact,
in this case the power of the tests seem to increase as the dimension is increased
from $p=4$ to $p=6$. In contrast, the tests $M_a$ exhibit a clear decrease
in power as the dimension $p$ is increased.

We finally consider testing $\mathcal H_0$ with $\alpha=1$, i.e.\ when testing for departures
from multivariate Cauchy. Overall, in agreement with previous observations,
Table~\ref{tab:alpha1_vs_stable_estimate}
% (also see Table~\ref{tab:alpha1_vs_stable_estimate_FULL} in the supplement)
shows that the new tests seem to be competitive in terms of power,
with the existing tests having some advantage when the true data-generating
distribution has a stability index greater than $1$.
However, this advantage in power disappears when considering the
higher-dimensional case $p=6$.
Similar behavior can be seen in the case of $t$ and skew Cauchy alternatives
(see Tables~S6 and S7
in the supplement).

\begin{table}[tbp]
	\centering\footnotesize
	\begin{threeparttable}
	\caption{Percentage of rejection of $\mathcal{H}_0$ with $\alpha=1$
					 against stable alternatives.
					 Tests done at the 10\% significance level.
            Also see Table~S5 of the
             supplement.}
	\label{tab:alpha1_vs_stable_estimate}
	\begin{tabular}{cccHHHHHHrrrrHHH|HrrrrHH}
		\toprule
		$p$ & $n$ & $\mathcal{H}_1$ & $T_{0.001}$ & $T_{0.005}$ & $T_{0.01}$ & $T_{0.05}$ & $T_{0.1}$
							& $T_{0.5}$ & $T_{1}$ & $T_{2}$ & $T_{5}$ & $T_{10}$ & $T_{20}$ & $T_{50}$ & $T_{100}$
						& $M_{2}$ & $M_{4}$ & $M_{6}$ & $M_{10}$ & $M_{15}$ & $M_{20}$ & $M_{25}$ \\
		\midrule
\multicolumn{1}{r}{4} & \multicolumn{1}{r}{50} &   $S_{1.5}$  & 83.3  & 80.6  & 77.6  & 40.1  & 12.3  & 26.9  & 48.8  & 48.5  & 77.7  & 92.0  & 90.5  & 54.5  & 0.0   & 69.2  & 88.2  & 77.2  & 50.8  & 26.9  & 18.3  & 15.4 \\
      &       &  $S_{1.2}$  & 40.2  & 36.8  & 33.6  & 14.6  & 8.3   & 12.8  & 17.6  & 13.6  & 15.6  & 26.3  & 24.2  & 8.0   & 0.4   & 17.8  & 31.2  & 22.9  & 16.9  & 13.3  & 14.2  & 9.6 \\
      &       &  $S_{0.8}$  & 1.6   & 3.0   & 4.0   & 11.7  & 14.9  & 20.6  & 25.0  & 36.7  & 57.5  & 61.9  & 63.3  & 64.2  & 60.0  & 32.7  & 18.0  & 13.4  & 12.7  & 9.5   & 10.0  & 10.7 \\
 & \multicolumn{1}{r}{100} &  $S_{1.5}$  & 98.3  & 97.4  & 95.6  & 13.8  & 8.5   & 57.4  & 84.6  & 94.6  & 99.8  & 100.0 & 100.0 & 98.2  & 74.0  & 97.2  & 98.8  & 96.2  & 71.7  & 38.2  & 25.7  & 20.6 \\
      &       &  $S_{1.2}$  & 63.5  & 58.3  & 49.1  & 8.7   & 7.4   & 17.1  & 26.2  & 26.9  & 46.6  & 60.5  & 59.4  & 37.1  & 13.1  & 36.7  & 49.1  & 36.7  & 20.6  & 12.9  & 12.7  & 10.7 \\
      &       &  $S_{0.8}$  & 1.1   & 3.2   & 5.8   & 14.4  & 17.1  & 38.2  & 50.8  & 67.0  & 77.8  & 80.8  & 83.3  & 82.4  & 78.2  & 55.9  & 40.0  & 25.9  & 12.1  & 10.8  & 11.7  & 8.1 \\
\midrule
\multicolumn{1}{r}{6} & \multicolumn{1}{r}{50} &  $S_{1.5}$  & 0.0   & 73.3  & 73.6  & 67.8  & 55.9  & 23.3  & 41.8  & 32.2  & 18.6  & 11.2  & 31.4  & 24.8  & 0.6   & 39.1  & 22.6  & 14.8  & 14.6  & 10.0  & 11.4  & 8.9 \\
      &       &  $S_{1.2}$  & 0.0   & 31.0  & 31.5  & 27.2  & 22.0  & 11.6  & 17.8  & 17.4  & 13.1  & 10.2  & 8.0   & 3.6   & 1.7   & 14.8  & 12.2  & 12.8  & 10.8  & 8.6   & 9.5   & 8.1 \\
      &       &  $S_{0.8}$  & 0.0   & 6.2   & 7.0   & 12.4  & 17.6  & 35.0  & 33.2  & 16.5  & 10.9  & 16.8  & 36.4  & 54.4  & 56.3  & 30.4  & 18.0  & 13.4  & 8.6   & 12.1  & 7.7   & 10.5 \\
   & \multicolumn{1}{r}{100} & $S_{1.5}$  & 0.0   & 98.7  & 98.7  & 97.6  & 95.0  & 49.8  & 96.2  & 98.9  & 96.9  & 99.7  & 100.0 & 99.3  & 91.8  & 87.4  & 52.4  & 19.0  & 13.6  & 14.2  & 13.5  & 11.8 \\
      &       &  $S_{1.2}$  & 0.0   & 67.8  & 66.5  & 53.9  & 41.4  & 9.4   & 34.6  & 40.7  & 19.8  & 52.3  & 63.5  & 45.3  & 23.9  & 23.9  & 18.0  & 11.9  & 9.4   & 8.5   & 11.1  & 10.9 \\
      &       &  $S_{0.8}$  & 0.0   & 1.9   & 2.8   & 10.7  & 19.6  & 50.8  & 67.6  & 72.1  & 59.5  & 69.6  & 81.4  & 81.5  & 80.3  & 52.2  & 30.4  & 14.4  & 9.4   & 10.4  & 11.7  & 11.5 \\
\bottomrule
	\end{tabular}
	\end{threeparttable}
\end{table}

\subsection{Monte Carlo results for GARCH data}\label{sec:GARCHcase}
We consider a CCC-GARCH$(1,1)$ model as defined in \eqref{CCC} with $\kappa_x=\kappa_q=1$. % and ${\mat{B}}_1$ assumed diagonal.
As parameters, we take ${\mat A}_1=0.2{\mat I}_p$, ${\mat B}_1=0.3{\mat I}_p$,
and correlation matrix $\mat{R}$ with all off-diagonal entries set to 0.5.

To test whether the innovations ${\vec\varepsilon}_j$
were generated by an elliptically symmetric SP law, we use the statistic
in \eqref{STS} applied to the residuals $\tilde{\vec\varepsilon}_j$
defined in \eqref{residual}. Denote the value of the test statistic
by $T_{r}:=T_{r}(\tilde{\vec\varepsilon}_1,\ldots,\tilde{\vec\varepsilon}_n)$.

We use the following bootstrap scheme to determine critical values:
\begin{enumerate}
	\item Independently generate innovations
                ${\vec\varepsilon}_1^*,\ldots,{\vec\varepsilon}_n^*$
				from the SP law $S_{\alpha_0}({\vec 0},{\mat I}_p)$.
        \item Construct a bootstrap sample ${\vec X}_1^*,\ldots,{\vec X}_n^*$
              using the recursive relation ${\vec X}_j^*=\tilde{\mat Q}_j^{1/2}{\vec\varepsilon}_j^*$, $j=1,\ldots,n$, with $\tilde{\mat Q}_j$
              as defined in \eqref{residual}.
	\item Based on the sample ${\vec X}_1^*,\ldots,{\vec X}_n^*$, fit the
                CCC-GARCH model in \eqref{CCC} to obtain estimates
                $\tilde{\mat Q}_j^*$, $j=1,\ldots,n$, and recover the bootstrap
                residuals $\tilde{\vec\varepsilon}_j^*=(\tilde{\mat Q}_j^{*})^{-1/2}{\vec X}_j^*$, $j=1,\ldots,n$
	\item A bootstrap version of the statistic is given by
				$T_{r}^{*}=T_{r}(\tilde{\vec\varepsilon}_1^*,\ldots,\tilde{\vec\varepsilon}_n^*)$.
\end{enumerate}
Steps 1--4 are repeated many times, say $B$, to obtain realisations
$T_{r}^{*}(1),\ldots,T_{r}^{*}(B)$ of
the bootstrap statistic.
The null hypothesis is rejected at significance level $\xi$ whenever
$T_{r}$ exceeds the
($1-\xi$)-level empirical quantile of the bootstrap realisations
$\{T_{r}^{*}(b)\}_{b=1}^B$.
In the Monte Carlo simulations that follow, instead of drawing $B$ bootstrap,
we employ the warp speed method of \citet{giacomini2013} which involves
drawing only one bootstrap sample for each Monte Carlo iteration.

Table~\ref{tab:alpha2_vs_stable_garch} (power against symmetric SP laws) exhibits similar favorable power properties
of the newly proposed test as was observed in the i.i.d.\ case.
Notice that the tests all have empirical size close to the nominal
level of 10\% when using critical values obtained by means of the bootstrap scheme
given above.
Also see Tables~S13 and S14
of the online supplement.

Table~S12 (power against skew normal alternatives)
shows that the new tests are especially
competitive in terms of power when the true innovation distribution
is not elliptically symmetric.
Finally we mention that although the tests of \citet{MNT} are competitive
when the innovations have an elliptically symmetric Student $t$-distribution
(see Table~S11), the slight
advantage in power disappears rapidly as the dimension $p$  increases.

\begin{table}[tbp]
	\centering\footnotesize
	\begin{threeparttable}
	\caption{Percentage of rejection of $\mathcal{H}_0$ with $\alpha=2$
             against stable alternatives in the CCC-GARCH$(1,1)$ case.
             % Tests done at the 10\% significance level.
             Also see Table~S10 of the online supplement.}
	\label{tab:alpha2_vs_stable_garch}
	\begin{tabular}{cccHHHHHHrrrrHHH|HrrrrHHr}
		\toprule
		$p$ & $n$ & $\mathcal{H}_1$ & $T_{0.001}$ & $T_{0.005}$ & $T_{0.01}$ & $T_{0.05}$ & $T_{0.1}$
							& $T_{0.5}$ & $T_{1}$ & $T_{2}$ & $T_{5}$ & $T_{10}$ & $T_{20}$ & $T_{50}$ & $T_{100}$
						& $M_{2}$ & $M_{4}$ & $M_{6}$
						& $M_{10}$ & $M_{15}$ & $M_{20}$
                        & $M_{25}$ & HJM \\
		\midrule
 \multicolumn{1}{r}{4} & \multicolumn{1}{r}{100} &  $S_{2.0}$  & 10.9  & 10.4  & 9.4   & 8.1   & 8.2   & 8.9   & 9.5   & 9.4   & 8.5   & 7.9   & 7.6   & 7.9   & 8.1   & 9.2   & 9.6   & 8.7   & 10.0  & 9.3   & 10.1  & 11.1  & 10.2 \\
      &       &  $S_{1.9}$  & 10.7  & 11.3  & 16.3  & 18.0  & 23.0  & 34.0  & 33.8  & 29.3  & 19.2  & 11.2  & 7.5   & 7.2   & 7.2   & 5.5   & 13.2  & 15.6  & 18.9  & 17.6  & 14.8  & 14.8  & 73.9 \\
      &       &  $S_{1.8}$  & 23.0  & 25.0  & 31.3  & 38.7  & 51.1  & 68.9  & 69.7  & 62.9  & 37.0  & 17.0  & 8.2   & 7.2   & 8.0   & 5.1   & 28.5  & 36.7  & 32.4  & 32.5  & 29.6  & 25.7  & 96.0 \\
      &       &  $S_{1.7}$  & 26.9  & 31.6  & 50.0  & 65.0  & 76.7  & 90.8  & 92.0  & 85.9  & 52.8  & 21.8  & 8.0   & 7.6   & 8.2   & 9.0   & 53.9  & 61.8  & 58.5  & 52.2  & 46.9  & 38.8  & 99.9 \\
   & \multicolumn{1}{r}{150} &  $S_{2.0}$  & 10.4  & 10.8  & 12.3  & 10.4  & 11.1  & 9.1   & 8.8   & 7.4   & 7.3   & 7.6   & 7.3   & 7.6   & 7.0   & 15.0  & 13.1  & 11.3  & 11.2  & 11.3  & 9.8   & 11.1  & 12.5 \\
      &       &  $S_{1.9}$  & 13.8  & 16.7  & 23.9  & 29.3  & 36.0  & 48.0  & 47.7  & 42.4  & 27.6  & 14.8  & 8.1   & 6.3   & 7.4   & 4.3   & 15.6  & 20.0  & 20.5  & 18.3  & 19.8  & 14.6  & 89.0 \\
      &       &  $S_{1.8}$  & 23.6  & 28.2  & 45.6  & 57.1  & 67.3  & 84.8  & 85.5  & 82.3  & 59.2  & 30.3  & 10.9  & 6.2   & 6.7   & 4.8   & 37.0  & 45.3  & 44.6  & 39.1  & 34.6  & 32.4  & 98.6 \\
      &       &  $S_{1.7}$  & 40.7  & 50.3  & 68.7  & 80.4  & 90.0  & 99.3  & 99.3  & 97.2  & 80.1  & 49.0  & 16.5  & 8.2   & 8.8   & 18.8  & 74.3  & 78.4  & 79.7  & 75.1  & 67.4  & 58.1  & 100.0 \\
\midrule
\multicolumn{1}{r}{6} & \multicolumn{1}{r}{100} &  $S_{2.0}$  & 12.1  & 13.2  & 11.0  & 11.5  & 11.5  & 11.7  & 10.6  & 9.1   & 8.9   & 7.9   & 8.1   & 8.4   & 8.5   & 10.6  & 10.3  & 11.3  & 8.5   & 9.9   & 8.2   & 9.0   & 8.1 \\
      &       &  $S_{1.9}$  & 3.5   & 8.1   & 18.6  & 21.1  & 26.3  & 41.3  & 40.8  & 35.4  & 24.4  & 14.4  & 9.4   & 7.7   & 7.8   & 2.5   & 11.3  & 12.6  & 12.9  & 10.2  & 9.2   & 9.3   & 81.0 \\
      &       &  $S_{1.8}$  & 5.8   & 13.0  & 28.3  & 34.8  & 47.0  & 75.9  & 76.3  & 69.8  & 41.8  & 22.2  & 12.0  & 9.9   & 9.9   & 1.4   & 18.1  & 22.6  & 15.7  & 14.9  & 12.2  & 12.1  & 95.8 \\
      &       &  $S_{1.7}$  & 9.4   & 22.7  & 46.3  & 56.5  & 72.8  & 92.3  & 92.3  & 85.1  & 59.6  & 27.1  & 7.6   & 5.5   & 6.1   & 2.3   & 36.5  & 37.4  & 31.6  & 23.3  & 19.3  & 18.3  & 99.1 \\
  & \multicolumn{1}{r}{150} &  $S_{2.0}$  & 8.9   & 11.4  & 12.1  & 12.3  & 12.2  & 8.4   & 10.8  & 9.8   & 7.7   & 7.3   & 8.1   & 8.8   & 9.6   & 7.2   & 9.2   & 7.4   & 10.3  & 10.1  & 10.2  & 11.0  & 10.0 \\
      &       &  $S_{1.9}$  & 5.4   & 14.2  & 21.1  & 22.4  & 31.4  & 53.3  & 54.3  & 48.9  & 34.9  & 21.8  & 14.4  & 11.2  & 11.2  & 2.1   & 13.6  & 13.6  & 11.1  & 11.2  & 12.1  & 10.8  & 93.1 \\
      &       &  $S_{1.8}$  & 7.9   & 20.0  & 31.9  & 44.8  & 65.0  & 89.8  & 88.8  & 83.7  & 60.6  & 32.5  & 15.8  & 9.0   & 9.4   & 1.0   & 20.4  & 21.7  & 18.6  & 17.8  & 11.2  & 11.8  & 99.2 \\
      &       &  $S_{1.7}$  & 13.4  & 31.2  & 51.9  & 71.6  & 88.4  & 99.7  & 99.1  & 97.5  & 79.1  & 50.6  & 13.8  & 7.5   & 7.5   & 5.3   & 51.2  & 54.4  & 41.8  & 32.6  & 25.0  & 18.6  & 100.0 \\
\bottomrule
	\end{tabular}
	\end{threeparttable}
\end{table}

\section{Testing asymmetric SP laws} \label{sec_7}

We now shift our focus to the more general case of testing whether multivariate observations from $\vec X$
originate from a SP  law, which need not necessarily be elliptically symmetric. In this connection note that the general multivariate SP law depends  on a location vector $\vec\delta\in\mathbb{R}^p$ and a spectral measure
$\Gamma(\cdot)$ on the unit sphere $S^p$. Accordingly we wish to test the null hypothesis 
\[
	\tilde{\cal H}_{0}:
	\vec X\sim \tilde {\cal{S}}_{\alpha}(\vec\delta,\Gamma), 
	 \text{ for some } \vec\delta\in\mathbb{R}^p \ \text{and some spectral measure } \Gamma \text{ on } S^p,
\]
where we write $\vec X \sim \tilde {\cal{S}}_\alpha(\vec\delta,\Gamma)$, when $\vec X$ follows a skew SP law with the indicated parameters. We will be considering the case of $\tilde{\cal H}_{0}$ for fixed $\alpha=\alpha_0$ as well as the case of testing the null hypothesis $\widetilde {\cal{H}}_0$ with an unspecified $\alpha$.

In this connection note that if $\vec X \sim \tilde {\cal{S}}_\alpha(\vec\delta, \Gamma)$, then the CF of $\vec X$ is given by 
% \begin{equation}\label{eq:mvstablelawcf}
% \varphi(\vec t)=\varphi_\alpha(\vec t;\vec \delta,\Gamma), \ \vec t \in \mathbb R^p, \end{equation}
% where 
% \begin{equation}\label{cfskew}
% \varphi_\alpha(\vec t;\vec \delta,\Gamma)= \exp\left(-\int_{\mathcal S^p} \psi_\alpha\left(\vec t^\top\vec s\right) \Gamma({\rm{d}}\vec s) + {\rm{i}} \:\vec t^\top \vec \delta \right),
% 	\quad\vec t\in\mathbb{R}^p,
% \end{equation}
\begin{equation}\label{eq:mvstablelawcf}
\varphi(\vec t)=\varphi_\alpha(\vec t;\vec \delta,\Gamma)
:= \exp\left(-\int_{\mathcal S^p} \psi_\alpha\left(\vec t^\top\vec s\right) \Gamma({\rm{d}}\vec s) + {\rm{i}} \:\vec t^\top \vec \delta \right),
	\quad\vec t\in\mathbb{R}^p,
\end{equation}
with
\[
	\psi_\alpha(u) =
	\begin{cases}
		\lvert u\rvert^\alpha\left(1 - {\rm{i}}\sign(u)\tan\left(\frac{1}{2}\pi\alpha\right)\right), & \text{for } \alpha\ne1, \\
		\lvert u\rvert\left(1 + {\rm{i}}\frac{2}{\pi}\sign(u)\log\lvert u\rvert\right), & \text{for } \alpha=1.
	\end{cases}
\]
See, e.g., \citet{nolan2001}.

In testing the null hypothesis $\tilde{\cal H}_{0}$, we consider an entirely different idea for the test statistic first put forward by \cite{chen2022}. Specifically, we 
consider a test statistic formulated as a two-sample test between the original
data $\mathcal{X}_n=(\vec X_j, \ j=1,...,n)$ and artificial data  $\mathcal{X}_{0n}=(\vec X_{0j}, \ j=1,...,n)$
generated under the null hypothesis $\tilde{\cal{H}}_0$.
More precisely, we propose the test criterion 
\begin{equation} \label{STS1}
\tilde T_{n,w}
:= \tilde T_{n,w}(\mathcal{X}_n;\mathcal{X}_{0n})
= n \int_{\vec t\in \mathbb R^p} \left\lvert\varphi_n(\vec t)-\varphi_{0n}(\vec t)\right\rvert^2 w(\vec t) {\rm{d}} \vec t, 
\end{equation}
where $\varphi_n(\vec t)=n^{-1} \sum_{j=1}^n e^{{\rm{i}} \vec t ^\top \vec X_j}$ is the empirical CF of the data at hand,
while $\varphi_{0n}(\vec t)=n^{-1} \sum_{j=1}^n e^{{\rm{i}} \vec t ^\top \vec X_{0j}}$
is the  empirical CF computed from the artificial data set $\mathcal{X}_{0n}$ generated under
the null hypothesis $\tilde{\cal H}_{0}$ with $\Gamma$ and $\vec\delta$ estimated
from the original observations $\mathcal{X}_{n}$.

By straightforward computations, we obtain
 \begin{equation}\label{STS2}
\tilde T_{n,w}(\mathcal{X}_n;\mathcal{X}_{0n})
=\frac{1}{n}\sum_{j,k=1}^n  \Big(I_w\left(\vec X_j-\vec X_k\right)+I_w\left(\vec X_{0j}-\vec X_{0k}\right)-2 I_w\left(\vec X_j-\vec X_{0k}\right)\Big), 
\end{equation}
where $I_w(\vec x)=\int_{\mathbb R^p} \cos(\vec t^\top \vec x)w(\vec t){\rm{d}}\vec t$, as also defined in \eqref{int1}.
Clearly then the numerical approaches of Section \ref{sec_3} are no longer required and the simplicity of this test  lies in the fact that in \eqref{STS2}  only the computation of $I_w(\cdot)$ is needed. Specifically the need for tailor-made weight functions such as those employed in Section \ref{sec_3} is circumvented. Furthermore we no longer need to compute the density of the underlying SP law as in \eqref{lambda}. In this connection, suppose that the weight function $w(\vec x)$  figuring in $I_w(\vec x)$ above is chosen as the  density of any spherical distribution
in $\mathbb  R^p$. Then the integral $I_w(\vec x)$  gives  the CF corresponding to this spherical distribution at the point $\vec x$.
Furthermore it is well known that this CF  may be  written as $\Psi(\|\vec x \|^2)$,
 where $\Psi(\cdot)$ is called the ``kernel" of the specific family of spherical distributions. Thus the test statistic figuring in \eqref{STS2} may be written as  
 \begin{multline}\label{STS3}
\tilde T_{n,\Psi}(\mathcal{X}_n;\mathcal{X}_{0n}) 
= \frac{1}{n}\sum_{j,k=1}^n  \Big(\Psi\left(\|\vec X_j-\vec X_k\|^2\right)\\+\Psi\left(\|\vec X_{0j}-\vec X_{0k}\|^2\right)-2 \Psi\left(\|\vec X_j-\vec X_{0k}\|^2\right)\Big), 
\end{multline}
where the kernel $\Psi(\cdot)$ can be chosen by the practitioner so that the resulting expression
in \eqref{STS3}  is as simple as possible. In this connection, as already clear from the preceding paragraphs,  a simple kernel is the
kernel of the spherical SP family of distributions with
$\Psi(\xi)=e^{-r \xi^{\alpha/2}}, \ r>0, \ \alpha \in (0,2]$. %and the kernel
%$\Psi(\xi)=(1+r \xi)^{-\alpha}, \ r,\alpha>0$, corresponding to the generalized Laplace family;
%see \citet{nolan2013} and \citet{kozubowski2013}, respectively.

Implementation of the test however relies on estimation of the spectral measure $\Gamma(\cdot)$
appearing in \eqref{eq:mvstablelawcf}.
Motivated by a result of \citet[Theorem~1]{byczkowski1993}, we assume that $\Gamma(\cdot)$ can be
approximated by the discrete spectral measure
$
	\widetilde\Gamma(\cdot) = \sum_{k=1}^K \gamma_k {\rm{I}}_{\vec s_k}(\cdot),
$
with weights $\gamma_k$ corresponding to mass points $\vec s_k\in S^p$, $k=1,\ldots,K$, and ${\rm{I}}_{\vec s_k}(\cdot)$ being the indicator index.
So in order to apply the test, we use the stochastic representation of $\vec X_0\sim \widetilde {\cal{S}}_{\alpha_0}(\vec\delta,\widetilde\Gamma)$  as
%A random vector satisfying $\tilde{\mathcal H}_0$ with stability index $\alpha_0$ is given by
\begin{equation}\label{eq:stochasticrepr}
	\vec X_{0} =
	\begin{cases}
		{\vec\delta} + \sum_{k=1}^K \gamma_k^{1/\alpha_0} A_k\vec{s}_k & \text{if } \alpha_0\ne1, \\
		{\vec\delta} + \sum_{k=1}^K \gamma_k\left(A_k + \frac{2}{\pi}\log(\gamma_k)\right)\vec{s}_k & \text{if } \alpha_0=1,
	\end{cases}
\end{equation}
where $(A_k, \ k=1,...,K)$ are i.i.d.  (univariate) SP variates following a totally skewed to the right SP law  with $\alpha=\alpha_0$. In turn this representation is used in order to generate observations $(\vec X_{0j}, \ j=1,...,n)$ under the null hypothesis $\widetilde {\cal{H}}_0$ with $\vec\delta$ and $(\gamma_k, \ k=1,...,K)$ replaced by appropriate estimates $\hat{\vec\delta}$ and $(\hat\gamma_k, \ k=1,...,K)$, respectively. The estimates 
$\hat{\vec\delta}$ and $(\hat\gamma_k, \ k=1,...,K)$ are obtained as shown in \citet{nolan2001}
and outlined in Section~S5 of the online supplement.
 
\subsection{Monte Carlo results}\label{sec_7:MCresults}
We now turn to a simulation study to demonstrate the performance of the test based on \eqref{STS2}, say 
$\tilde T_{\Psi}$ for simplicity, in the bivariate case. We specifically consider the following alternative
distributions:
\begin{enumerate}[label=(A\arabic*)]
	\item \label{skewalt1}
            the asymmetric SP law $\tilde {\cal{S}}_\alpha(\vec\delta,\Gamma)$ with CF defined in \eqref{eq:mvstablelawcf},
            where we take $\vec\delta=\vec 0$ and $\Gamma(\cdot)=\sum_{k=1}^L \gamma_k {\rm{I}}_{\vec s_k}(\cdot)$,
            with $L=5$, $(\gamma_1,\ldots,\gamma_L) = (0.1,0.3,0.1,0.4,0.1)$ and
            \begin{equation}\label{eq:masspoints}
            \vec s_k = \left[\cos\left(\frac{2\pi(k-1)}{L}\right), \sin\left(\frac{2\pi(k-1)}{L}\right)\right]^\top,
            \qquad
             k=1,\ldots,L;
            \end{equation}
	\item \label{skewalt2}
            spherically symmetric SP distributions, denoted by $S_{\alpha}$.
\end{enumerate}

The statistic in $\tilde T_{\Psi}(\mathcal{X}_{n};\mathcal{X}_{0n})$ in \eqref{STS3} is subject to randomness
introduced by the artificial data $\mathcal{X}_{0n}$. To address this randomness
in practical implementation, we follow \citet{chen2022} and base our test on the statistic
\begin{equation}\label{eq:avgstat}
	\tilde T_{\Psi}^{(m)} = \frac{1}{m}\sum_{r=1}^m \tilde T_{n,\Psi}(\mathcal{X}_{n};\mathcal{X}_{0n}^r),
\end{equation}
where, for each $(r=1,\ldots,m)$, the set $\mathcal{X}_{0n}^r$ is a random sample of observations
generated from the SP law $\tilde {\cal{S}}_{\alpha_0}(\hat{\vec\delta},\hat\Gamma)$, i.e.\ a random
sample satisfying $\tilde{\mathcal{H}}_0$.
% The estimates $\hat{\vec\delta}$ and $\hat\Gamma$ are calculated as described
% in \citet{nolan2001}, where we take
% $\hat\Gamma(\cdot) = \sum_{k=1}^K \hat\gamma_k {\rm{I}}_{\vec s_k}(\cdot)$ with
% $K=10$ mass points at ${\vec s_k}$ as defined in \eqref{eq:masspoints}.
				
Critical values of the test can be obtained using the following parametric bootstrap scheme:
\begin{enumerate}
	\item Generate a bootstrap sample $\mathcal{X}_{n}^*$
				from the SP law $\tilde {\cal{S}}_{\alpha_0}(\hat{\vec\delta},\hat\Gamma)$.
	\item Calculate bootstrap estimates $\hat{\vec\delta}^*$ and $\hat\Gamma^*$ from $\mathcal{X}_{n}^*$ and
				generate $m$ random samples $(\mathcal{X}_{0,n}^{*r}$, $r=1,\ldots,m)$,
				from the SP law $\tilde {\cal{S}}_{\alpha_0}(\hat{\vec\delta}^*,\hat\Gamma^*)$.
	\item A bootstrap version of the statistic is
				$
					\tilde T_{\Psi}^{*(m)} = m^{-1}\sum_{r=1}^m \tilde T_{\Psi}(\mathcal{X}_{n}^*;\mathcal{X}_{0,n}^{*r}).
				$
				% \[
				% 	\tilde T_{\Psi}^{*(m)} = \frac{1}{m}\sum_{r=1}^m \tilde T_{\Psi}(\mathcal{X}_{n}^*;\mathcal{X}_{0,n}^{*r}).
				% \]
\end{enumerate}
Steps 1--3 are repeated many times, say $B$, to obtain realizations
$\tilde T_{\Psi}^{*(m)}(1),\ldots,\tilde T_{\Psi}^{*(m)}(B)$ of
the bootstrap statistic.
The null hypothesis is rejected at significance level $\xi$ whenever
the statistic $\tilde T_{\Psi}^{(m)}$ in \eqref{eq:avgstat} exceeds the
($1-\xi$)-level empirical quantile of the sequence of bootstrap realizations
$\{\tilde T_{\Psi}^{*(m)}(b)\}_{b=1}^B$.

Table~\ref{tab:skew_alpha1.5} shows the empirical rejection percentages
of tests based on $\tilde T_{\Psi}$ with $\Psi(\xi)=e^{-r \xi^{\alpha_0/2}}$, where
$\alpha_0$ denotes the hypothesized stability index in $\tilde{\mathcal{H}}_0$.
We write $\tilde T_r$, $r>0$, for this test. In the case of unknown $\alpha$ in Table \ref{tab:skew_iid_estimate_alpha} we use the same weight function with $\alpha_0=2$. 

\begin{table}[bp]
	\centering\footnotesize
	\begin{threeparttable}
	\caption{Percentage of rejection of $\tilde{\mathcal{H}}_0$ with $\alpha=1.5$
					 against bivariate stable alternatives.
					 Tests done at the 10\% significance level.
            Also see Table~S15 of the
            online supplement.}
	\label{tab:skew_alpha1.5}
	\begin{tabular}{HccHHrrrrHHHHHHHHHrrrrHHHHHHHHrrrrr}
		\toprule
            &&&\multicolumn{13}{c}{Skew SP distributions}&\multicolumn{13}{c}{Symmetric SP distributions} \\
            \cmidrule(lr){4-16}\cmidrule(lr){17-29}
		$p$ & $n$ & $\mathcal{H}_1$ &
            $\tilde T_{0.001}$ & $\tilde T_{0.005}$ & $\tilde T_{0.01}$ & $\tilde T_{0.05}$ & $\tilde T_{0.1}$
							& $\tilde T_{0.5}$ & $\tilde T_{1}$ & $\tilde T_{2}$ & $\tilde T_{5}$ & $\tilde T_{10}$ & $\tilde T_{20}$ & $\tilde T_{50}$ & $\tilde T_{100}$
             & $\tilde T_{0.001}$ & $\tilde T_{0.005}$ & $\tilde T_{0.01}$ & $\tilde T_{0.05}$ & $\tilde T_{0.1}$
							& $\tilde T_{0.5}$ & $\tilde T_{1}$ & $\tilde T_{2}$ & $\tilde T_{5}$ & $\tilde T_{10}$ & $\tilde T_{20}$ & $\tilde T_{50}$ & $\tilde T_{100}$ \\
		\midrule
2     & 100   & $\tilde{\cal{S}}_{2.0}$   & 33.7  & 52.5  & 59.0  & 73.7  & 71.8  & 59.6  & 51.2  & 46.4  & 34.4  & 23.6  & 17.3  & 12.6  & 33.9  & 57.2  & 66.3  & 73.6  & 68.1  & 56.2  & 50.0  & 40.6  & 24.1  & 16.7  & 11.9  & 9.9 \\
      &       & $\tilde{\cal{S}}_{1.7}$   & 7.0   & 10.7  & 13.4  & 17.8  & 19.4  & 17.2  & 15.6  & 15.3  & 12.2  & 9.1   & 8.6   & 8.9   & 5.5   & 9.7   & 13.2  & 18.7  & 18.5  & 16.8  & 14.6  & 12.4  & 9.3   & 9.1   & 9.7   & 9.7 \\
      &       & $\tilde{\cal{S}}_{1.5}$   & 15.0  & 14.2  & 13.3  & 9.8   & 8.9   & 7.4   & 8.0   & 8.4   & 9.9   & 8.9   & 8.0   & 9.4   & 9.2   & 8.4   & 9.0   & 8.5   & 8.4   & 8.3   & 7.9   & 7.1   & 9.4   & 8.2   & 7.9   & 8.8 \\
      &       & $\tilde{\cal{S}}_{1.4}$  & 23.6  & 21.8  & 20.6  & 16.5  & 15.2  & 13.9  & 15.6  & 15.4  & 15.4  & 14.0  & 13.5  & 13.5  & 18.6  & 17.8  & 17.1  & 13.8  & 11.2  & 11.4  & 12.0  & 12.6  & 12.1  & 11.4  & 11.6  & 11.3 \\
      &       & $\tilde{\cal{S}}_{1.3}$   & 40.4  & 38.4  & 37.5  & 30.3  & 27.8  & 28.7  & 29.6  & 30.5  & 30.3  & 28.3  & 25.1  & 22.1  & 34.6  & 31.8  & 31.2  & 24.3  & 21.1  & 24.1  & 26.7  & 27.1  & 23.4  & 20.7  & 18.8  & 16.5 \\
      % &       & $\tilde{\cal{S}}_{1.2}$   & 55.5  & 58.0  & 55.1  & 48.4  & 44.5  & 44.9  & 49.0  & 51.5  & 51.7  & 48.3  & 43.6  & 35.4  & 54.8  & 55.2  & 53.1  & 47.0  & 43.2  & 48.9  & 54.0  & 52.7  & 47.0  & 43.4  & 39.3  & 29.4 \\
      & 150   & $\tilde{\cal{S}}_{2.0}$   & 40.8  & 70.5  & 81.2  & 91.7  & 90.5  & 78.6  & 71.4  & 61.8  & 43.9  & 29.1  & 20.0  & 11.8  & 46.3  & 79.0  & 85.7  & 91.2  & 84.5  & 73.2  & 70.8  & 54.1  & 35.4  & 23.1  & 17.0  & 13.2 \\
      &       & $\tilde{\cal{S}}_{1.7}$  & 5.8   & 11.3  & 17.3  & 22.2  & 24.4  & 19.3  & 18.9  & 15.4  & 13.1  & 12.4  & 9.7   & 7.8   & 6.6   & 13.1  & 19.3  & 24.3  & 23.3  & 22.6  & 18.0  & 14.3  & 11.8  & 11.0  & 9.8   & 11.1 \\
      &       & $\tilde{\cal{S}}_{1.5}$   & 11.0  & 7.7   & 8.6   & 6.3   & 6.6   & 6.2   & 8.1   & 7.5   & 6.9   & 7.1   & 8.1   & 9.5   & 8.0   & 7.4   & 8.2   & 8.4   & 7.1   & 5.9   & 7.5   & 8.4   & 7.7   & 9.6   & 11.9  & 8.7 \\
      &       & $\tilde{\cal{S}}_{1.4}$   & 22.5  & 19.7  & 20.2  & 16.3  & 13.8  & 12.3  & 14.6  & 15.7  & 15.1  & 13.4  & 12.3  & 12.6  &    18.1 &  23.2 &  22.4 &  13.8 &  10.6 &   7.1 &   9.1 &   9.6 &   9.5 &  10.7 &  14.0 &  13.1   \\
      &       & $\tilde{\cal{S}}_{1.3}$   & 40.2  & 40.0  & 39.8  & 33.8  & 29.8  & 25.5  & 31.7  & 32.1  & 29.6  & 26.6  & 22.5  & 18.5  &   41.0 &  43.3 &  42.6 &  33.6 &  27.5 &  24.0 &  25.3 &  25.3 &  23.3 &  20.3 &  16.6 &  13.2   \\
      % &       & $\tilde{\cal{S}}_{1.2}$   & 63.1  & 67.9  & 67.4  & 64.5  & 59.2  & 58.1  & 62.3  & 64.8  & 58.0  & 48.7  & 41.9  & 37.8  & 64.6  & 68.7  & 70.7  & 57.6  & 47.2  & 49.2  & 56.3  & 55.7  & 47.8  & 42.6  & 37.8  & 27.1 \\
		\bottomrule
	\end{tabular}
	\end{threeparttable}
\end{table}

The left-hand side of Table~\ref{tab:skew_alpha1.5} shows the rejection percentages
when observations
are sampled from an asymmetric SP distribution, whereas the right-hand side show
the results when observations are sampled from a symmetric SP distribution.
In all cases, the proposed procedure seems to respect the nominal size of the test,
although it seems to be somewhat conservative especially for smaller sample sizes.
% In agreement with earlier observations, there is also some under-rejection
% in the case of lighter-tailed alternatives, but this seems to disappear as the sample
% size is increased.
Nevertheless, the tests have good power against alternatives,
which increases with the extent of violation of the null hypothesis.
% Results for sample size $n=200$ are shown in Table~\ref{tab:skew_alpha1.5_FULL}
% of the supplement, and
Corresponding results when testing $\tilde{\mathcal{H}}_0$ with $\alpha=1.8$
are given in Table~S16.

\subsubsection*{The case of unknown $\alpha$}
In practice, the true value of stability index $\alpha$ will typically be unknown
and need to be estimated from data. Suppose the hypothesis of interest is
\[
	\tilde{\cal H}_{0}':
	\vec X\sim \tilde {\cal{S}}_{\alpha}(\vec\delta,\Gamma), 
	 \text{ for some } \alpha\in(0,2), \vec\delta\in\mathbb{R}^p \text{ and some measure } \Gamma \text{ on } S^p,
\]
To test this hypothesis, we again use the test based on the statistic in \eqref{STS3},
but now generate the artificial data ${\mathcal{X}}_{0n}$ from the SP law
${\mathcal{S}}_{\hat\alpha}(\hat{\vec \delta},\hat\Gamma)$, where $\hat\alpha$,
$\hat{\vec \delta}$ and $\hat\Gamma$ are projection estimates obtained as outlined in
Section~S5 of the supplement. Note that data can be generated from
${\mathcal{S}}_{\hat\alpha}(\hat{\vec \delta},\hat\Gamma)$ using the stocahstic
representation in \eqref{eq:stochasticrepr} with $\alpha_0$ replaced by $\hat\alpha$.
The bootstrap procedure for obtaining critical values follows similarly to the procedure
described in Section~\ref{sec_7:MCresults}.

The empirical rejection percentages (for the bivariate case, $p=2$) are shown in Table~\ref{tab:skew_iid_estimate_alpha}
for several distributions. Note that when observations are sampled from a $\mathcal{S}_{1.8}$
distribution (symmetric SP law with stability index 1.8), the rejection percentages
are close to the nominal level, indicating that the test has reasonable empirical size.

As alternatives, we consider the skew normal ($SN_{\nu}$),
symmetric Laplace and generalized Gaussian ($GG_{\nu})$ distributions. An observation $X$
from the $GG_{\nu}$ distribution (refer to \citealp{cadirci2022})
is generated according to $X=UV^{1/\nu}$, where $U$ is uniform on $S^{p-1}$ and
$V\sim\text{Gamma}(p/\nu,2)$.
Table~\ref{tab:skew_iid_estimate_alpha} shows that the proposed test procedure
has power against non-SP alternatives, and noting the increase in power associated with
increasing sample size, the results suggest that the procedure is consistent
against non-SP alternatives.

\begin{table}[tb]
	\centering\footnotesize
	\begin{threeparttable}
	\caption{Percentage of rejection of $\tilde{\mathcal{H}}_0'$ using tests based on
                  the statistic in \eqref{STS3} with $\Psi(\xi)=e^{-r\xi}$.
				 Tests done at the 10\% significance level.
            Also see Table~S17 of the
            online supplement.}
	\label{tab:skew_iid_estimate_alpha}
	\begin{tabular}{HlHrHHHHHrrrrrHHHHHrr}
		\toprule
		$p$ & $\mathcal{H}_1$ & & $n$ & $T_{0.001}$ & $T_{0.005}$ & $T_{0.01}$ & $T_{0.05}$ & $T_{0.1}$
							& $T_{0.5}$ & $T_{1}$ & $T_{2}$ & $T_{5}$ & $T_{10}$ & $T_{20}$ & $T_{50}$ & $T_{100}$ \\
		\midrule
            2 & ${\mathcal{S}}_{1.8}$ & 1.8 &  50 &   5.5 &   5.6 &   5.8 &   5.2 &   5.1 &   4.3 &   5.2 &   5.1 &   6.7 &   6.8 &   7.8 &   7.8 &   8.3 \\
 &  & 1.8 &  75 &   6.6 &   5.2 &   5.7 &   5.5 &   5.0 &   6.7 &   7.5 &   8.1 &  10.2 &  10.1 &   9.7 &   9.2 &   9.8 \\
 &  & 1.8 & 100 &   5.5 &   4.7 &   4.8 &   4.7 &   5.2 &   6.5 &   7.7 &   8.0 &   8.2 &   8.8 &   9.5 &  10.0 &  11.0 \\
 &  & 1.8 & 150 &   5.6 &   6.9 &   7.2 &   7.4 &   7.1 &   7.8 &   9.0 &   9.0 &   9.2 &   9.2 &   9.5 &  11.1 &   9.7 \\
 & $SN_{1.5}$ & 1.5 &  50 &   2.9 &   2.6 &   2.9 &   4.0 &   4.4 &  10.2 &  12.5 &  12.8 &  13.2 &  11.1 &  10.7 &  10.8 &  10.3 \\
 &  & 1.5 &  75 &   4.3 &   4.2 &   5.0 &   6.2 &   7.0 &  12.8 &  16.6 &  16.4 &  14.4 &  13.9 &  13.8 &  11.8 &  13.0 \\
 &  & 1.5 & 100 &   3.2 &   4.2 &   4.7 &   6.2 &   7.5 &  16.4 &  19.4 &  19.5 &  17.3 &  15.4 &  12.6 &   9.8 &   8.6 \\
 &  & 1.5 & 150 &   4.0 &   5.1 &   5.7 &   7.4 &   9.7 &  18.0 &  24.6 &  23.5 &  19.4 &  16.8 &  14.5 &  12.3 &  11.5 \\
 % &  & 1.5 & 200 &   6.4 &   6.6 &   8.1 &   9.9 &  12.8 &  25.8 &  30.2 &  29.4 &  24.4 &  21.4 &  17.0 &  15.0 &  13.8 \\
% \cmidrule(lr){2-17}
 & $SN_{2.0}$ & 2.0 &  50 &   2.9 &   3.5 &   4.2 &   3.1 &   5.2 &  10.6 &  14.5 &  13.9 &  13.9 &  12.3 &  12.0 &  12.2 &  12.0 \\
 &  & 2.0 &  75 &   3.1 &   3.8 &   3.5 &   4.6 &   4.8 &  12.5 &  16.1 &  17.1 &  15.8 &  14.2 &  13.5 &  12.2 &  10.6 \\
 &  & 2.0 & 100 &   2.3 &   3.0 &   3.2 &   6.1 &   9.5 &  19.3 &  23.2 &  24.5 &  21.7 &  18.3 &  16.3 &  12.1 &  11.6 \\
 &  & 2.0 & 150 &   4.1 &   4.6 &   4.2 &   5.3 &   8.5 &  23.1 &  33.8 &  34.1 &  30.9 &  28.2 &  21.5 &  17.5 &  14.4 \\
 % &  & 2.0 & 200 &   4.4 &   6.4 &   7.7 &  10.4 &  15.7 &  34.2 &  42.8 &  42.3 &  34.0 &  31.4 &  23.6 &  16.4 &  12.7 \\
% \cmidrule(lr){2-17}
%  & Exponential & 0.5 &  20 &  31.3 &  33.0 &  34.3 &  36.4 &  40.3 &  46.7 &  43.3 &  39.4 &  29.7 &  28.1 &  27.6 &  28.6 &  27.3 \\
%  &  & 0.5 &  50 &  87.5 &  87.5 &  87.9 &  87.3 &  87.8 &  86.8 &  87.3 &  87.4 &  87.3 &  86.7 &  86.8 &  87.9 &  88.6 \\
%  &  & 0.5 &  75 &  97.9 &  98.8 &  98.9 &  99.0 &  99.0 &  99.1 &  99.2 &  99.3 &  99.0 &  99.0 &  99.0 &  99.0 &  98.8 \\
 % &  & 0.5 & 100 &  99.4 &  99.6 &  99.7 &  99.7 &  99.7 &  99.7 &  99.9 &  99.9 &  99.9 &  99.9 &  99.9 &  99.9 &  99.9 \\
% \cmidrule(lr){2-17}
 & Laplace & 0.0 &  50 &   3.9 &   6.3 &   6.2 &   7.2 &   7.3 &   9.0 &  10.8 &  13.2 &  16.3 &  16.7 &  15.9 &  14.6 &  16.0 \\
 &  & 0.0 &  75 &   2.2 &   3.2 &   3.9 &   6.7 &   7.3 &  15.2 &  22.8 &  30.4 &  36.3 &  35.5 &  34.3 &  27.0 &  22.1 \\
 &  & 0.0 & 100 &   3.0 &   5.4 &   6.5 &  12.5 &  14.9 &  31.3 &  42.6 &  50.4 &  51.8 &  48.5 &  41.3 &  31.7 &  25.3 \\
 &  & 0.0 & 150 &   2.9 &   6.0 &   7.0 &  16.0 &  19.1 &  34.7 &  53.1 &  65.2 &  64.7 &  64.4 &  57.0 &  44.4 &  36.8 \\
 % &  & 0.0 & 200 &   4.4 &   8.7 &   9.6 &  17.1 &  20.5 &  48.7 &  73.4 &  82.0 &  83.3 &  79.3 &  70.2 &  55.8 &  42.4 \\
% \cmidrule(lr){2-17}
%  & Logistic & 0.0 &  50 &   3.3 &   5.7 &   7.4 &  12.2 &  14.8 &  14.9 &  14.0 &  11.6 &  10.7 &  10.7 &  10.4 &   9.4 &   8.1 \\
%  &  & 0.0 &  75 &   3.2 &   5.7 &   7.5 &  16.7 &  18.5 &  13.9 &  12.5 &  13.1 &  12.3 &  10.2 &   9.7 &   8.0 &   9.5 \\
%  &  & 0.0 & 100 &   3.8 &   7.1 &   9.8 &  17.9 &  20.6 &  19.0 &  16.1 &  13.8 &  12.6 &  11.1 &  11.0 &  11.6 &  13.2 \\
%  &  & 0.0 & 150 &   6.3 &   9.6 &  14.0 &  25.9 &  27.9 &  23.4 &  20.7 &  18.1 &  13.2 &  12.4 &  12.4 &   9.9 &  11.4 \\
%  &  & 0.0 & 200 &   6.4 &  12.3 &  20.2 &  36.0 &  34.8 &  31.0 &  27.6 &  22.5 &  18.4 &  17.2 &  15.1 &  11.9 &  12.0 \\
% \cmidrule(lr){2-17}
 & $GG_1$ & 1.0 &  50 &   2.4 &   2.6 &   3.5 &   4.9 &   6.4 &   8.3 &  11.0 &  14.4 &  17.5 &  19.3 &  16.8 &  13.2 &  12.3 \\
 &  & 1.0 &  75 &   2.3 &   2.5 &   2.9 &   4.8 &   6.4 &  12.8 &  16.5 &  21.1 &  26.9 &  23.6 &  22.4 &  16.6 &  14.4 \\
 &  & 1.0 & 100 &   0.7 &   1.8 &   2.6 &   5.4 &   8.7 &  15.5 &  21.0 &  32.0 &  33.4 &  31.6 &  28.7 &  21.9 &  18.6 \\
 &  & 1.0 & 150 &   1.9 &   2.7 &   3.4 &   7.4 &  12.2 &  21.5 &  33.6 &  47.4 &  50.4 &  46.3 &  38.8 &  28.1 &  22.7 \\
 % &  & 1.0 & 200 &   0.8 &   4.5 &   5.5 &   9.9 &  13.1 &  25.1 &  43.7 &  57.3 &  64.1 &  59.1 &  48.7 &  37.5 &  30.0 \\
		\bottomrule
	\end{tabular}
	\end{threeparttable}
\end{table}

\subsection{The high-dimensional case} \label{high_dim}

 It should be clear that so far, and despite the fact that the new test applies to any dimension $p$, the underlying setting is not that of high dimension  ($p>n$). In this connection we point out that the extension of goodness-of-fit methods specifically tailored for the classical ``small $p$--large $n$" regime to high or infinite dimension is not straightforward, and therefore it is beyond the scope of the present article. This is not restricted to our setting alone but applies more generally, and the collections of \cite{goia2016} and \cite{kokoszka2017} reflect the need for statistical methods specifically tailored for non-classical settings. If we restrict it to our context, in the latter settings such methods have so far been mostly confined to testing for normality and the interested reader is referred to \cite{bugni2009}, \cite{nieto2014}, \cite{barcenas2017}, \cite{kellner2019}, \cite{yamada2019}, \cite{jiang2019}, \cite{gorecki2020}, \cite{henze2021}, \cite{chenxia2023} and \cite{chengenton2023}. If however the setting is that of  regression, with or without conditional heteroskedasticity, the number of parameters rapidly increases with the dimension $p$. This is one extra reason that corresponding specification methods in high/infinite dimension need to be treated separately, and in this connection the methods of \cite{cuesta2019} and  \cite{rice2020} appear to be of the few available for testing the (auto)regression function.   

%Cuesta-Albertos et al. (2019), Ann.Statist. 47, 439-467
%Kokoszka et al. (2017), Econometrics and Statistics, Vol. 1, 99--100
%Goia and Vieu, JMVA 146 (2016), 1-6.
%Rice, G., Wirjanto, T., Zhao, Y. (2020). Tests for conditional heteroscedasticity of functional data. Journal of Time Series Analysis. 41(6), 733-758. <doi:10.1111/jtsa.12532>.

%Bugni, FA and Hall, P and Horowitz, JL and Neumann, GR, Econometr. J., 12 (2009),  pages={S1--S18},

%Henze, N and Jim\'enez--Gamero, MD, Scand. J. Statist., 48 (2021), 406--428

%Jiang, Q and Hu\v{s}kov\'a, M and Meintanis, SG and Zhu, LX , J. Multivar. Anal., 170 (2019), 202--220

%Kellner and Celisse, Bernoulli 25 (2019), 1816-1837. 

%Gorecki et al., Inter.Statist.Rev. (2020) 

%Chen ang Genton, Inter.Statist.Rev. 91 (2023), 114-139

%Chen ang Xia, JASA 118 (2023), 719-731

%Yamada and Himeno , Comput. Statist. 34 (2019), 911-941. 

%Write something about approx. of eigenvalues and chi-squared approxim. See Canad.J.Statist.

%Barcenas et al. (2017), TEST 26, 503-526

%Cuesta-Albertos et al. (2014), CSDA 75, 124-141.

For an illustration of the special circumstances that are brought forward as the dimension grows, consider the test statistic  
in \eqref{STS} for $\alpha_0=2$ and without standardization, i.e.\ replace  $\bm Y_j$ by $\bm X_j\sim{\cal S}_2(\vec 0, {\mat{ I}}_p)$, $j=1,...,n$. Then from \eqref{lambda}, and by using the density of the  $p$-variate normal distribution with mean zero and covariance matrix $2{\mat {I}}_p$, we obtain 
\[
\Lambda_r(\bm x;2)= \left(\frac{\pi}{r}\right)^{\frac{p}{2}} e^{-\|\bm x\|^2/(4r)}, 
\]
and consequently  
 \begin{equation*}\begin{split} 
T_{n,r}
={}&{} \frac{1}{n} \left (\sum_{j=1}^n \Lambda_r(\bm 0;2)+\sum_{j\neq k}^n \Lambda_r(\bm X_j-\bm X_k;2)\right) +n\Lambda_{r+2}(\bm 0;2) \\
{}&{}-2 \sum_{j=1}^n  \Lambda_{r+1}(\bm X_j;2) \\  
={}&{} \frac{1}{n} \left (\sum_{j=1}^n\left(\frac{\pi}{r}\right)^{\frac{p}{2}} +\sum_{j\neq k}^n \left(\frac{\pi}{r}\right)^{\frac{p}{2}}  e^{-\|\bm X_j-\bm X_k\|^2/(4r)}\right)+n \left(\frac{\pi}{r+2}\right)^{\frac{p}{2}}\\
{}&{}- 2 \sum_{j=1}^n   \left(\frac{\pi}{r+1}\right)^{\frac{p}{2}} e^{-\|\bm X_j\|^2/(4(r+1))}.
\end{split}\end{equation*}
% Now the quantity $\|x\|^2$  behaves like $\{p+{\rm{O}}_{\mathbb P}(\sqrt{p})\}^2=p^2(1+{\rm{O}}_{\mathbb P}(1/\sqrt{p}))$, and thus our test statistic contains sums of terms each of which is of the order $e^{-p^2}$, as $p\to\infty$. Hence
Now the quantity $\|\vec X_j\|^2=\sum_{k=1}^pX_{jk}^2$  behaves like $2p+{\rm{O}}_{\mathbb P}(\sqrt{p})=2p(1+{\rm{O}}_{\mathbb P}(1/\sqrt{p}))$, and thus our test statistic contains sums of terms, each of which is of the order $e^{-2p}$, as $p\to\infty$. (For simplicity we suppress the terms $4r$ and $4(r+1)$ which occur as denominators in these sums, as they are anyway irrelevant to our argument).

In order to get a feeling of this result, write $\|\vec X_j\|^2=2\sum_{k=1}^p(X_{jk}/\sqrt{2})^2=:2 S_p$, where, due to Gaussianity and independence, $S_p$ is distributed as chi-squared with $p$ degrees of freedom. Thus the expectation of the quantity $e^{-\|\vec X_j\|^2}$    figuring in $T_{n,r}$ coincides with the value of the CF of this chi-squared distribution computed at the point $t=2{\rm{i}}$. To proceed further, notice that the CF of $S_p/p$ at a point $t$ is given by the CF corresponding to $S_p$ computed at $t/p$, and recall that the CF  of the chi--squared distribution with $p$ degrees of freedom is given by $\varphi_{S_p}(t)=(1-2{\rm{i}}t)^{-p/2}$. Hence, in obvious notation,   
\[
\lim_{p\to\infty} \varphi_{\frac{S_p}{p}}(t)=\lim_{p\to\infty} \left(1-\frac{2{\rm{i}}t}{p}\right)^{-p/2}=e^{{\rm{i}}t}, 
\]
and reverting back to the CF of $S_p$ we get $\varphi_{S_p}(t)\approx e^{{\rm{i}}pt}$, and thus $\varphi_{S_p}(2{\rm{i}})\approx e^{-2p}$, as $p\to \infty$.  A similar reasoning applies to $\|\vec X_j-\vec X_k\|^2$ implying that   
the expectation of  $e^{-\|\vec X_j-\vec X_k\|^2}$ (also occurring in $T_{n,r}$) being approximated by $e^{-4p}$, as $p\to\infty$.   

Hence, by using these approximations, we obtain      
\[
T_{n,r} \approx \left(\frac{\pi}{r}\right)^{\frac{p}{2}} +n \left(\frac{\pi}{r+2}\right)^{\frac{p}{2}}, \ p\to\infty, 
\]
from which it follows that 
\[
	T_{n,r} \to
	\begin{cases}
		\infty & \text{for } r<\pi, \\
		1, & \text{for } r=\pi, \\
        0, & \text{for } r>\pi. 
	\end{cases}
\]
as $p\to\infty$. 
Consequently our test statistic degenerates in high dimension, a fact that calls for proper
high-dimensional modifications of the test criterion, which is definitely a worthwhile subject for future research.

Nevertheless, we have obtained some initial Monte Carlo results that show a reasonable performance for the test criterion in cases where the dimension is much higher than the maximum $p=6$ considered so far. These results may be found in Tables~S18 and S19 of the Supplement.

\section{Application to financial data} \label{sec_8}
We consider daily log returns from 4 January 2010 to 30 June 2017
of stocks of two mining companies listed on the London Stock Exchange:
Anglo American (AAL) and Rio Tinto (RIO). The complete data set (available from
Yahoo!\ Finance) consists of 1,891 log returns.
We model the log returns using the CCC-GARCH$(1,1)$ model (with intercept) given by
\begin{equation}\label{CCC_intercept}
    {\vec Y}_j = {\vec\omega} + {\mat Q}^{1/2}_j{\vec\varepsilon}_j,
    \quad
    j=1,\ldots,n,
\end{equation}
where ${\mat Q}_j$ is as in \eqref{CCC} with $\kappa_x=\kappa_q=1$
and ${\mat{B}}_1$ assumed diagonal.

We are interested in determining whether the innovations ${\vec\varepsilon}_j$
have a bivariate stable distribution.
Since the innovations are unobserved, we apply our test to the
residuals $\hat{\vec\varepsilon}_j=\hat{\mat Q}^{-1/2}_j({\vec Y}_j - \hat{\vec\omega})$,
where the estimates $\hat{\mat Q}_j$ and $\hat{\vec\omega}$ are obtained using
EbE estimation as discussed earlier.
To obtain critical values of the tests, we apply the bootstrap algorithm
of Section~\ref{sec:GARCHcase} with $B=1,000$.

Table~\ref{tab:VaR} show that, when a stability index of $\alpha\in\{1.75,1.8\}$
is assumed, the tests based on $T_{r}$ do not reject the
null hypothesis that the CCC-GARCH-$(1,1)$ innovations have a ${\cal{S}}_\alpha$
distribution (at a 10\% level of significance). 
On the other hand, the null hypothesis of stable innovations is rejected
when a stability index of $\alpha=\{1.7,1.85,1.9,2\}$ is assumed.

The correct choice of the innovation distribution has important implications
for value-at-risk (VaR) forecasts. For the considered stable distributions,
we fit the model in \eqref{CCC_intercept} to the first 1,000 observations
and calculated one-step-ahead 5\% and 1\% portfolio VaR forecasts for long and short
positions for the remaining time period (i.e.\ 891 forecasts for each position).
The portfolio is assumed to consist of 50\% AAL shares and 50\% RIO shares.

Table~\ref{tab:VaR} shows the empirical coverage rates of the forecasted VaR bounds,
that is, the proportion of times that the value of the portfolio exceeded the bounds.
For the cases where our test supports the null hypothesis, i.e.\ when
$\alpha\in\{1.75,1.8\}$, the empirical coverage rates of the value-at-risk bounds
are quite close to the nominal rates. In addition, the p-values of
Christofferson's (\citeyear{christoffersen1998}) $LR_{cc}$ test
(given in brackets in Table~\ref{tab:VaR}), indicate that,
if the stability index of the innovation distribution is chosen
either too high or too low, the true conditional coverage rates of the forecasted
VaR bounds are significantly different from the nominal rates.

\begin{table}[tbp]
	\centering\footnotesize
	\begin{threeparttable}
	\caption{p-values of tests that the CCC-GARCH innovations have a
             ${\cal{S}}_\alpha$ distribution for several choices of $\alpha$,
             along with empirical coverage rates of forecasted VaR
             bounds calculated under ${\cal{H}}_0$.             
             p-values of the $LR_{cc}$ test are given in brackets.
             All p-values less than 10\% are underlined.}
	\label{tab:VaR}
	\begin{tabular}{lrrrrrrrrrrrrrrrrrrr}
		\toprule
        & \multicolumn{2}{c}{}
            & \multicolumn{2}{c}{5\% VaR coverage rate}
            & \multicolumn{2}{c}{1\% VaR coverage rate} \\
        \cmidrule(lr){4-5}\cmidrule(lr){6-7}
		${\cal{H}}_0$ & $T_{0.1}$ & $T_{0.5}$
          & \multicolumn{1}{c}{Long pos.}
		 & \multicolumn{1}{c}{Short pos.}
		 & \multicolumn{1}{c}{Long pos.}
		 & \multicolumn{1}{c}{Short pos.} \\
		% \midrule
  %       ${\cal{S}}_{2.0}$
  %           && 0.0418 & 0.0260 & 0.0136 & 0.0068 \\
  %           && 0.4694 & 0.0014 & 0.5075 & 0.571 \\
  %       ${\cal{S}}_{1.9}$
  %           && 0.0531 & 0.0362 & 0.0147 & 0.0056 \\
  %           && 0.8545 & 0.1402 & 0.3473 & 0.3579 \\
  %       ${\cal{S}}_{1.8}$
  %           && 0.0520 & 0.0350 & 0.0090 & 0.0023 \\
  %           && 0.8751 & 0.0992 & 0.8916 & 0.0203 \\
  %       ${\cal{S}}_{1.7}$
  %           && 0.0508 & 0.0328 & 0.0034 & 0.0023 \\
  %           && 0.8710 & 0.0443 & 0.0722 & 0.0203 \\
		\midrule
        ${\cal{S}}_{2.00}$ & \underline{0.004} & \underline{0.006} &
            0.064 (\underline{0.041}) & 0.065 (0.126) & 0.025 (\underline{0.001}) & 0.028 (\underline{0.000}) \\
        ${\cal{S}}_{1.90}$ & \underline{0.013} & \underline{0.025} &
            0.064 (\underline{0.041}) & 0.065 (0.126) & 0.018 (\underline{0.074}) & 0.018 (\underline{0.074}) \\
        ${\cal{S}}_{1.85}$ & \underline{0.048} & \underline{0.081} &
            0.063 (\underline{0.059}) & 0.065 (0.126) & 0.016 (0.227) & 0.013 (0.519) \\
        ${\cal{S}}_{1.80}$ & 0.214 & 0.286 &
            0.063 (\underline{0.059}) & 0.063 (0.224) & 0.010 (0.912) & 0.009 (0.887) \\
        ${\cal{S}}_{1.75}$ & 0.393 & 0.496 &
            0.059 (0.148) & 0.059 (0.444) & 0.008 (0.758) & 0.007 (0.560) \\
        ${\cal{S}}_{1.70}$ & \underline{0.020} & \underline{0.018} &
            0.054 (\underline{0.056}) & 0.059 (0.444) & 0.004 (0.177) & 0.002 (\underline{0.019}) \\
        % \midrule
        % ${\cal{S}}_{\hat\alpha}$ &   &   &
        %     0.065 (0.028) & 0.066 (0.091) & 0.022 (0.004) & 0.026 (0.000) \\
		\bottomrule
	\end{tabular}
	\end{threeparttable}
\end{table}

\section{Conclusion} \label{sec_9}
We have studied goodness-of-fit tests with data involving multivariate SP laws. Our tests include the case of i.i.d. observations as well as the one with observations from GARCH models,  and cover both elliptical and skewed distributions. Moreover they refer to hypotheses whereby some parameters are assumed known, as well as to the full composite hypothesis with all parameters estimated from the data at hand. The new procedures are shown to perform well in finite samples and to be competitive against other methods, whenever such methods are available. An application  illustrates the usefulness of the new procedure for modeling stock returns and  explores the subsequent forecasting implications.

\backmatter

\bmhead{Supplementary information}

Some technical material and additional Monte Carlo results are provided in
an accompanying online supplement. The {\textsf {R}} codes for the Monte Carlo implementations of Sections \ref{sec_6} and \ref{sec_7} are available from the authors upon request.

\bmhead{Acknowledgments}

The authors wish to sincerely thank Christian Francq (Univesrity of Lille, CREST) for providing the code for EbE estimation.

\section*{Statements and Declarations}

No funds, grants, or other support was received for conducting this study.
The authors have no relevant financial or non-financial interests to disclose. The authors have no conflicts of interest to declare.

\end{document}

%%%%%%%%%%%%%%%%%%%%%%%%%%%%%%%%%%%%%%%%%%%%%%%%%%%%%%%%%%%%%%%%%%%%%%%%%%%

\title[Specification procedures for multivariate stable-Paretian laws]{Online supplement to: Specification procedures for multivariate stable-Paretian laws for independent and for conditionally heteroskedastic data}


\author*[1,2]{\fnm{Simos G.} \sur{Meintanis}}\email{simosmei@econ.uoa.gr}

\author[3]{\fnm{John P.} \sur{Nolan}}

\author[4]{\fnm{Charl} \sur{Pretorius}}

%\author[5]{\fnm{Zhou} \sur{Zhou}}

\affil*[1]{\orgdiv{Department of Economics}, \orgname{National and Kapodistrian University of Athens}, \orgaddress{\city{Athens}, \country{Greece}}}

\affil[2]{\orgdiv{Pure and Applied Analytics}, \orgname{North-West University}, \orgaddress{\city{Potchefstroom}, \country{South Africa}}}

\affil[3]{\orgdiv{Department of Statistics}, \orgname{American University}, \orgaddress{\city{Washington DC}, \country{USA}}}

\affil[4]{\orgdiv{Centre for Business Mathematics and Informatics}, \orgname{North-West University}, \orgaddress{\city{Potchefstroom}, \country{South Africa}}}

%\affil[5]{\orgdiv{Department of Statistical Sciences}, \orgname{University of Toronto}, \orgaddress{\state{Ontario}, \country{Canada}}}

%\keywords{Empirical characteristic function, Goodness-of-fit, Heavy-tailed distribution}

\abstract{This supplement contains some technical material (Sections \ref{sec:Kotzweight}--\ref{sec:estGamma}) and additional Monte Carlo results (Tables \ref{tab:alpha2_vs_stable_estimateFULL}--\ref{tab:skew_iid_estimate_alpha_FULL}).}

\maketitle

\section{Computation using the CF of the Kotz--type law}\label{sec:Kotzweight}

The integral 
\begin{equation*}
I_{\nu,r}(\bm x;s)=\int_{\mathbb R^p} \cos(\bm t^\top \bm x)
\left(\|\bm t\|^2\right)^{\nu}e^{-r (\|\bm t\|^2)^s} {\rm{d}} \bm t
\end{equation*}
required for the test statistic $T_{n,w}$ can be computed  by making use of the CF of the Kotz-type distribution. Specifically for $\bm x=\bm 0$, we have 
\begin{equation*}
I_{\nu,r}(\bm 0;s)=\int_{\mathbb R^p} 
\left(\|\bm t\|^2\right)^{\nu}e^{-r (\|\bm t\|^2)^s} {\rm{d}} \bm t=\frac{{\pi^{\frac{p}{2}}\Gamma(\frac{2\nu+p}{2s})}}
{{s\Gamma(\frac{p}{2})r^{\frac{2\nu+p}{2s}}}
}.  
\end{equation*}
On the other hand, if $\bm x\neq\bm 0$, the value of $I_{\nu,r}(\bm x;s)$ may be computed from absolutely convergent series  
which for $1/2<s<1$ is given by   \citet{streit1991} as 
\begin{equation*} \label{int6}
I_{\nu,r}(\bm x;s)=\frac{ \pi^{p/2}}{s r^{\frac{2\nu+p}{2s}}}\sum_{k=0}^\infty \frac{1}{k!} \left(-\frac{\|\bm x\|^2}{4 r^{1/s}}\right)^k \frac{\Gamma(\frac{2\nu+p+2k}{2s})}{\Gamma(\frac{p+2k}{2})},  
\end{equation*}
while for  $0<s<1/2$ we have
\begin{multline*} \label{int8} \scriptsize
I_{\nu,r}(\bm x;s)=\frac{\pi^{\frac{p}{2}-1}}{r^{\frac{2\nu+p}{2 s}}}\sum_{k=0}^\infty \frac{(-1)^{k+1}}{k!} 2^{2 sk+2\nu+p} \: \Gamma\left(sk+\nu+\frac{p}{2}\right) \\\times\Gamma\left(sk+\nu+1\right)\sin\left(\pi(sk+\nu\right))\left(\frac{\|\bm x\|^2}{r^{1/s}}\right)^{-s k-\nu-\frac{p}{2}};  
\end{multline*}
see \citet{kotz1994}. 
%\end{equation}

In the special case $s=1$ the computation simplifies to
\begin{equation*} \label{int7}
I_{\nu,r}(\bm x;1)=\frac{ \pi^{\frac{p}{2}} \Gamma(\frac{2\nu+p}{2})
}{ r^{\frac{2\nu+p}{2}}}e^{-\frac{\|\bm x\|^2}{4r}}\sum_{k=0}^\nu \begin{pmatrix} \nu  \\
 k    \end{pmatrix} \frac{\left(-\frac{\|\bm x\|^2}{4 r}\right)^k} {\Gamma(\frac{p+2k}{2})},  
\end{equation*}
see \citet{Iyengar1989}, and for $s=1/2$ the bivariate case was treated by \cite{NK01} leading to 
\begin{multline*} \label{int10} \scriptsize
I_{\nu,r}\left(\bm x;\frac{1}{2}\right)=\frac{2\pi \Gamma(2(\nu+1)) (\nu+1)!}{4^\nu r^{2(\nu+1)}(\|\bm x\|^2+r^2)^{2\nu+\frac{3}{2}}}\sum_{k=0}^\nu 
\frac{(2(\nu-k))! (2k)!}{(\nu-k)! (k!)^2} (\|\bm x\|^2+r^2)^{k}\\
\times\sum_{\ell=0}^{\nu-k}\frac{(-1)^\ell 4^\ell r^{4\nu+3-2k-2\ell} \|\bm x\|^{2\ell} }{(\nu+1-\ell)!(\nu-k-\ell)! (2\ell)!},
\end{multline*}
which for $\nu=0$ simplifies to
\[
I_{0,r}\left(\bm x;\frac{1}{2}\right)=\frac{2\pi r}{(\|\bm x\|^2+r^2)^{3/2}}. 
\]

%\begin{rem*}
%Note that for $0<s<1/2$ we use the representation of the CF of the spherical Kotz type distribution derived by \citet{kotz1994}. We note in passing that this representation  is incorrectly reported in Theorem 3 of that paper and reproduced in this incorrect form in \citet{nadarajah2003}. The correct formula may be obtained from  \citet[p.~177]{kotz1994}.  
%\end{rem*}

\section{Affine invariance}\label{sec:affine.invariance}

%\subsection{Affine invariance}
A desirable feature of potential estimators $\widehat{\vec\delta}_n$ and $\widehat{\bf Q}_n$
of the location vector $\vec\delta$ and the dispersion matrix $\bf Q$ are
the following equivariance/invariance properties 
\begin{equation*}\label{del}
\widehat {\bm \delta}_n(\bm A\bm X_1+\bm b,...,\bm A \bm X_n+\bm b)=\bm  A\widehat {\bm \delta}_n(\bm X_1,...,\bm X_n)  +\bm b
\end{equation*}
and
\begin{equation*}\label{q} 
\widehat {\bf Q}_n(\bm A\bm X_1+\bm b,...,\bm A \bm X_n+\bm b)=\bm  A\widehat {\bf Q}_n(\bm X_1,...,\bm X_n)\bm A^{\top},
\end{equation*}
for each $\bm b \in \mathbb R^p$ and each non-singular $(p\times p)$ matrix $\bm A$. As a consequence the test statistic $T_{n,w}:=T_{n,w}(\bm X_1,...,\bm X_n)$ satisfies 
\begin{equation*}\label{AI}
T_{n,w}(\bm A\bm X_1+\bm b,...,\bm A\bm X_n+\bm b)=T_{n,w}(\bm X_1,...,\bm X_n), 
\end{equation*}
i.e. it is affine invariant. We note that this property is  in line with the fact that if $\bm X \sim {\cal{S}}_\alpha(\bm \delta, \bf Q)$, then $\bm A \bm  X+\bm b \sim {\cal{S}}_\alpha(\bm A \bm \delta+\bm b, \bm A\bf Q \bm A^\top)$, meaning that  the SP family of distributions is itself invariant
with respect to affine transformations $\bm X\mapsto \bm A\bm X+\bm b$.  
%In turn, if a test statistic is affine invariant then its finite-sample as well as its limit null distribution is independent of the true parameters $\bm \delta$ and $\bf Q$, and thus they may both  be approximated by simple Monte Carlo by setting the
%location vector and the dispersion matrix equal to their standard values zero and identity matrix,  respectively.       

\section{EbE estimator for GARCH parameters}\label{sec:EbE}
We outline the procedure for computing the EbE estimator of ${\vec\vartheta}$ of \cite{francq2014}. 
Note in this connection that, under the CCC--GARCH model, ${\vec X}_j={\mat D}_j{\bf R}^{1/2}{\vec\varepsilon}_j=:{\mat D}_j{\vec\varepsilon}_j'$. Suppose the conditional dispersion of $X_{k,j}$ 
(the $k^{\rm{th}}$ component of ${\vec X}_j$), is parameterized by ${\vec\vartheta}_n^{(k)}$,
i.e.\ $\tilde q_{k,j}=\tilde q_{k,j}({\vec\vartheta}_n^{(k)})$.
The EbE estimator of ${\vec\vartheta}_n^{(k)}$ is given by
\[
	\hat{\vec\vartheta}_n^{(k)}
	= \argmax_{{\vec\vartheta}_n^{(k)}} \;
			\sum_{j=1}^n \left[- \log\tilde{q}_{k,j}({\vec\vartheta}_n^{(k)})
										+ \log \bar{f}_\alpha\left(\frac{\tilde\varepsilon_{k,j}}{\tilde{q}_{k,j}({\vec\vartheta}_n^{(k)})}\right)\right],
\]
where $\bar{f}_\alpha$ denotes the density of a symmetric univariate SP law
with stability index equal to $\alpha$, with location zero and dispersion equal to one.
Letting
${\widetilde{\mat D}}_j={\rm{diag}}(\tilde q_{1,j}(\hat{\vec\vartheta}_n^{(1)}),\ldots, \tilde q_{p,j}(\hat{\vec\vartheta}_n^{(p)}))$,
we then obtain the correlated residuals $\tilde{\vec\varepsilon}_j'={\widetilde{\mat D}}_j^{-1}{\vec X}_j$
from which ${\bf R}$ is calculated using maximum likelihood to obtain
$\tilde{\vec\varepsilon}_j={\widetilde{\bf R}}_j^{-1/2}\tilde{\vec\varepsilon}_j'$.

\section{Calculation of the stable density}\label{sec:calc.stabledens}

%\subsection{Calculation of the multivariate stable density}
Implementation of our test procedure  relies on the evaluation
of the density of a $p$-variate spherical SP law.
Efficient evaluation of the density is also needed for maximum likelihood estimation of the parameters.

%Following the idea of \citet{mittnik1999} for univariate densities,
%\citet{bonato2012} uses a discrete fast Fourier transform over a grid of equally-spaced points in order to evaluate this density by means of the inversion theorem for CFs.  
% This approach yields a smooth density suitable for maximum likelihood
%estimation, but is computationally infeasible in higher dimensions.

In our numerical work, we utilize the fact  that if $\vec X$ follows a spherical SP law with CF $\phi_{\alpha}(\cdot)$, then  the density of $\vec X$ can be expressed as
\[
	f_{\alpha}(\vec x)
	= \begin{cases}
			\frac{\Gamma(p/2)}{2\pi^{p/2}}\|\vec x\|^{1-p}f_R(\|\vec x\|) & \text{if } \|\vec x\| \ne 0, \\
			\frac{\Gamma(p/\alpha)}{\alpha 2^{p-1}\pi^{p/2}\Gamma(p/2)} & \text{if } \|\vec x\| = 0,
		\end{cases}
\]
where $f_R(\cdot)$ is the density of $\|\vec X\|$, the amplitude of $\vec X$.
This reduces the problem of calculating $f_{\alpha}$ to calculating the univariate density $f_R$.
Various integral expressions for $f_R$ are given in \citet{nolan2013} and an implementation to numerically evaluate $f_R$
is the function \texttt{damplitude} in the $\textsf{R}$ package \texttt{stable} \citep[provided by][]{stablepackage}.

To speed up calculations, we pre--calculate $f_R(u)$ for each $u\in\{k/(N-k), k=0,\ldots,N-1\}$,
for some large value of $N$ ($=10,000$ in our simulations). Intermediary points are approximated
using cubic spline interpolation and, for $u\ge N$, we set $f_R(u)=0$.

\section{Estimation of the discrete spectral measure}\label{sec:estGamma}
Below we outline the projection estimation procedure of \citet{nolan2001}
as implemented in our work. The estimation procedure assumes that the stability
index $\alpha$ is known and that the data have been centered. As this is usually
not the case, we first estimate the one-dimensional parameters 
$(\alpha_j,\beta_j,\sigma_j,\delta_j)$, $j=1,\ldots,p$, for each
of the coordinates of the $p$-dimensional data set and center the data set
using the location estimate ${\vec\delta}=(\delta_1,\ldots,\delta_p)$.
Furthermore, if $\alpha$ is unknown, we estimate it by 
$\hat\alpha=p^{-1}\sum_{j=1}^p\hat\alpha_j$.

Motivated by a Theorem~1 of \citet{byczkowski1993}, we assume that $\Gamma(\cdot)$ can be
approximated by the discrete spectral measure
$
	\widetilde\Gamma(\cdot) = \sum_{k=1}^K \gamma_k {\rm{I}}_{\vec s_k}(\cdot)
$
as defined in the main paper, where we take
\[%\label{eq:masspoints.est}
    \vec s_k = \left[\cos\left(\frac{2\pi(k-1)}{K}\right), \sin\left(\frac{2\pi(k-1)}{K}\right)\right]^\top,
    \qquad
     k=1,\ldots,K.
\]
For all estimates calculated in our simulations, we took the number of projections as $K=10$.

For each $k=1,\ldots,K$, the projections
${\vec X}_1^\top{\vec s}_k,\ldots,{\vec X}_n^\top{\vec s}_k$ are i.i.d.\ with
a centered, univariate $\alpha$-stable distribution with skewness and dispersion parameters
$\beta_k$ and $\sigma_k$.
These parameters are estimated using maximum likelihood, after which the estimates
of the $K$ projections are combined using the method described in Section~2.2
of \citet{nolan2001} to recover an estimate of the spectral measure $\tilde\Gamma$ above.

\begin{table}[p]
	\centering\footnotesize
	\setlength{\tabcolsep}{3pt}
	\begin{threeparttable}
	\caption{Percentage of rejection of $\mathcal{H}_0$ with $\alpha=2$
					 against stable alternatives.
					 Tests done at the 10\% significance level.
                      Due to the computational complexity of the HJM test,
                      it is only included in selected cases to reduce run time.}
	\label{tab:alpha2_vs_stable_estimateFULL}
\begin{tabular}{cccHHHHHHrrrrHHH|HrrrrrHr}
		\toprule
		$p$ & $n$ & $\mathcal{H}_1$ & $T_{0.001}$ & $T_{0.005}$ & $T_{0.01}$ & $T_{0.05}$ & $T_{0.1}$
							& $T_{0.5}$ & $T_{1}$ & $T_{2}$ & $T_{5}$ & $T_{10}$ & $T_{20}$ & $T_{50}$ & $T_{100}$
						& $M_{2}$ & $M_{4}$ & $M_{6}$ & $M_{10}$ & $M_{15}$ & $M_{20}$ & $M_{25}$ & HJM \\
\midrule      &       &       & T0.001 & T0.005 & T0.01 & T0.05 & T0.1  & T0.5  & T1    & T2    & T5    & T10   & T20   & T50   & T100  & M2    & M4    & M6    & M10   & M15   & M20   & M25   & HJM \\
\midrule
\multicolumn{1}{r}{4} & 50    &  $S_{1.9}$  & 3.1   & 13.9  & 15.2  & 17.8  & 21.3  & 34.1  & 41.7  & 52.5  & 59.9  & 60.9  & 61.8  & 62.0  & 62.0  & 6.9   & 12.3  & 15.8  & 15.1  & 13.4  & 14.2  & 15.0  & 62.6 \\
      &       &  $S_{1.8}$  & 3.0   & 17.6  & 21.8  & 32.4  & 37.9  & 57.8  & 68.1  & 76.6  & 82.9  & 82.8  & 83.2  & 83.2  & 83.2  & 5.4   & 20.0  & 24.2  & 27.7  & 26.1  & 24.5  & 21.7  & 84.0 \\
      &       &  $S_{1.7}$  & 4.9   & 26.2  & 30.1  & 48.7  & 60.0  & 82.0  & 87.9  & 92.1  & 93.0  & 93.2  & 93.5  & 93.5  & 93.4  & 9.0   & 37.1  & 43.8  & 44.3  & 41.6  & 35.0  & 29.4  & 94.6 \\
      &       &  $S_{1.5}$  & 13.1  & 42.8  & 52.6  & 75.9  & 88.0  & 98.5  & 99.1  & 99.4  & 99.6  & 99.4  & 99.4  & 99.3  & 99.3  & 28.4  & 74.8  & 80.3  & 81.8  & 77.0  & 69.3  & 62.6  & 99.8 \\
\cmidrule{2-24}      & 75    &  $S_{1.9}$  & 4.2   & 15.1  & 14.9  & 22.8  & 28.0  & 42.0  & 50.7  & 61.4  & 68.8  & 71.0  & 71.5  & 71.4  & 71.4  & 5.6   & 12.0  & 15.9  & 15.2  & 16.6  & 14.2  & 12.6  &  \\
      &       &  $S_{1.8}$  & 5.7   & 20.8  & 22.9  & 43.2  & 52.5  & 74.0  & 82.3  & 87.7  & 90.8  & 91.0  & 91.2  & 91.2  & 91.2  & 4.3   & 22.4  & 32.0  & 32.3  & 32.4  & 31.3  & 23.8  &  \\
\bottomrule
	\end{tabular}
	\end{threeparttable}
\end{table}

\begin{table}[p]
	\centering\footnotesize
	\begin{threeparttable}
	\setlength{\tabcolsep}{3pt}
	\caption{Percentage of rejection of $\mathcal{H}_0$ with $\alpha=2$
					 against $t$ alternatives.
					 Tests done at the 10\% significance level.
                      Due to the computational complexity of the HJM test,
                      it is only included in selected cases to reduce run time.}
	\label{tab:alpha2_vs_t_estimateFULL}
	\begin{tabular}{cccHHHHHHrrrrHHH|HrrrrrHr}
		\toprule
		$p$ & $n$ & $\mathcal{H}_1$ & $T_{0.001}$ & $T_{0.005}$ & $T_{0.01}$ & $T_{0.05}$ & $T_{0.1}$
							& $T_{0.5}$ & $T_{1}$ & $T_{2}$ & $T_{5}$ & $T_{10}$ & $T_{20}$ & $T_{50}$ & $T_{100}$
						& $M_{2}$ & $M_{4}$ & $M_{6}$ & $M_{10}$ & $M_{15}$ & $M_{20}$ & $M_{25}$ & HJM \\
		\midrule
  \multicolumn{1}{r}{4} & \multicolumn{1}{r}{50} &  $t_{10}$  & 1.2   & 10.5  & 12.6  & 17.5  & 17.4  & 29.6  & 39.3  & 50.4  & 57.2  & 56.0  & 56.6  & 56.6  & 56.6  & 3.5   & 12.9  & 14.9  & 12.5  & 11.3  & 14.0  & 9.5   & 59.1 \\
      &       &  $t_{ 5}$  & 1.3   & 16.4  & 20.1  & 29.1  & 40.3  & 71.9  & 80.6  & 86.4  & 89.6  & 89.4  & 89.1  & 88.5  & 88.4  & 4.0   & 28.2  & 33.3  & 28.9  & 22.3  & 16.4  & 13.2  & 90.8 \\
      &       &  $t_{ 4}$  & 2.8   & 17.1  & 20.6  & 42.6  & 61.0  & 87.0  & 92.8  & 94.6  & 95.1  & 94.5  & 94.3  & 94.2  & 94.4  & 6.3   & 46.4  & 51.7  & 43.6  & 34.9  & 26.0  & 21.1  & 97.2 \\
      &       &  $t_{ 2}$  & 14.7  & 60.0  & 74.0  & 94.5  & 97.9  & 99.9  & 100.0 & 100.0 & 99.9  & 99.9  & 99.9  & 99.9  & 99.9  & 58.7  & 94.3  & 96.5  & 93.9  & 91.1  & 82.6  & 78.5  & 100.0 \\
\cmidrule{2-24}      & \multicolumn{1}{r}{75} &  $t_{10}$  & 2.5   & 11.2  & 11.5  & 15.8  & 20.6  & 38.6  & 45.8  & 55.1  & 61.0  & 62.5  & 62.8  & 62.5  & 62.8  & 3.9   & 16.6  & 18.2  & 15.9  & 15.0  & 12.5  & 11.4  &  \\
      &       &  $t_{ 5}$  & 4.2   & 16.2  & 18.2  & 37.2  & 55.0  & 85.6  & 90.6  & 93.6  & 94.2  & 93.8  & 92.8  & 92.3  & 91.8  & 6.8   & 41.1  & 44.9  & 39.1  & 28.6  & 21.9  & 19.5  &  \\
      &       &  $t_{ 4}$  & 5.5   & 25.2  & 28.0  & 60.3  & 77.4  & 96.2  & 97.9  & 98.6  & 98.4  & 98.0  & 97.6  & 96.9  & 96.5  & 14.7  & 62.9  & 65.2  & 52.8  & 44.6  & 37.3  & 32.1  &  \\
\bottomrule
	\end{tabular}
	\end{threeparttable}
\end{table}

\begin{table}[p]
	\centering\footnotesize
	\begin{threeparttable}
	\setlength{\tabcolsep}{3pt}
	\caption{Percentage of rejection of $\mathcal{H}_0$ with $\alpha=2$
					 against skew normal alternatives.
					 Tests done at the 10\% significance level.
                      Due to the computational complexity of the HJM test,
                      it is only included in selected cases to reduce run time.}
	\label{tab:alpha2_vs_skew_estimateFULL}
	\begin{tabular}{cccHHHHHHrrrrHHH|HrrrrrHr}
		\toprule
		$p$ & $n$ & $\mathcal{H}_1$ & $T_{0.001}$ & $T_{0.005}$ & $T_{0.01}$ & $T_{0.05}$ & $T_{0.1}$
							& $T_{0.5}$ & $T_{1}$ & $T_{2}$ & $T_{5}$ & $T_{10}$ & $T_{20}$ & $T_{50}$ & $T_{100}$
						& $M_{2}$ & $M_{4}$ & $M_{6}$ & $M_{10}$ & $M_{15}$ & $M_{20}$ & $M_{25}$ & HJM \\
         \midrule\multicolumn{1}{r}{4} & \multicolumn{1}{r}{50} &  $SN_{2.0}$  & 8.0   & 9.7   & 10.3  & 13.1  & 16.7  & 24.2  & 32.3  & 34.2  & 35.5  & 35.6  & 33.7  & 33.5  & 33.6  & 11.8  & 14.0  & 14.5  & 15.6  & 9.7   & 11.8  & 13.6  & 17.2 \\
          &       &  $SN_{1.5}$  & 10.3  & 10.1  & 9.7   & 11.0  & 12.7  & 21.8  & 29.1  & 32.0  & 30.6  & 29.5  & 28.7  & 28.2  & 28.0  & 12.9  & 12.9  & 12.1  & 13.5  & 9.6   & 12.3  & 11.0  & 16.3 \\
          &       &  $SN_{1.0}$  & 10.3  & 11.2  & 11.3  & 10.6  & 12.7  & 16.9  & 20.2  & 22.0  & 24.1  & 23.3  & 23.3  & 23.4  & 23.3  & 11.6  & 11.5  & 11.9  & 11.8  & 9.2   & 12.7  & 9.1   & 14.4 \\
          &       &  $SN_{0.5}$  & 9.8   & 9.8   & 9.6   & 9.7   & 10.3  & 11.4  & 11.9  & 13.6  & 15.0  & 14.1  & 14.3  & 13.9  & 14.1  & 9.5   & 7.8   & 7.1   & 10.2  & 12.1  & 10.5  & 9.0   & 12.6 \\
\cmidrule{2-24}		    & \multicolumn{1}{r}{75} &  $SN_{2.0}$  & 9.2   & 10.0  & 10.3  & 10.7  & 14.6  & 27.1  & 35.2  & 42.9  & 47.0  & 46.9  & 46.0  & 44.7  & 43.6  & 12.5  & 13.0  & 11.4  & 12.5  & 11.7  & 9.2   & 11.6  &  \\
		    &       &  $SN_{1.5}$  & 8.7   & 9.8   & 10.2  & 8.8   & 14.0  & 21.7  & 27.6  & 32.8  & 37.8  & 37.5  & 36.9  & 36.0  & 35.4  & 11.5  & 10.9  & 10.0  & 12.6  & 15.0  & 13.9  & 10.0  &  \\
		    &       &  $SN_{1.0}$  & 9.9   & 9.3   & 9.3   & 8.4   & 11.8  & 16.1  & 17.5  & 21.2  & 23.4  & 23.9  & 24.4  & 23.9  & 23.6  & 11.2  & 10.6  & 9.5   & 10.0  & 9.2   & 8.8   & 13.9  &  \\
\bottomrule
	\end{tabular}
	\end{threeparttable}
\end{table}

\begin{table}[p]
	\centering\footnotesize
	\begin{threeparttable}
	\setlength{\tabcolsep}{3pt}
	\caption{Percentage of rejection of $\mathcal{H}_0$ with $\alpha=1.8$
					 against stable alternatives.
					 Tests done at the 10\% significance level.}
	\label{tab:alpha1.8_vs_stable_estimate_FULL}
	\begin{tabular}{cccHHHHHHrrrrHHH|HrrrrrH}
		\toprule
		$p$ & $n$ & $\mathcal{H}_1$ & $T_{0.001}$ & $T_{0.005}$ & $T_{0.01}$ & $T_{0.05}$ & $T_{0.1}$
							& $T_{0.5}$ & $T_{1}$ & $T_{2}$ & $T_{5}$ & $T_{10}$ & $T_{20}$ & $T_{50}$ & $T_{100}$
						& $M_{2}$ & $M_{4}$ & $M_{6}$ & $M_{10}$ & $M_{15}$ & $M_{20}$ & $M_{25}$ \\
		\midrule
\multicolumn{1}{r}{4} & \multicolumn{1}{r}{50} &  $S_{2.0}$  & 31.4  & 19.3  & 10.6  & 10.4  & 10.6  & 11.8  & 9.1   & 3.1   & 0.1   & 0.0   & 0.0   & 0.0   & 0.0   & 14.7  & 12.3  & 12.2  & 12.1  & 12.5  & 12.4  & 10.6 \\
      &       &  $S_{1.9}$  & 19.5  & 13.2  & 9.8   & 10.3  & 9.6   & 9.6   & 9.1   & 6.0   & 2.0   & 1.8   & 1.9   & 2.0   & 3.0   & 10.8  & 10.0  & 9.3   & 11.1  & 8.3   & 10.1  & 9.4 \\
      &       &  $S_{1.7}$  & 6.4   & 8.4   & 10.1  & 10.9  & 10.5  & 11.3  & 13.2  & 18.8  & 31.2  & 32.4  & 31.1  & 27.5  & 24.6  & 11.9  & 12.7  & 11.8  & 11.5  & 9.4   & 10.9  & 11.9 \\
      &       &  $S_{1.5}$  & 3.5   & 9.2   & 12.6  & 13.4  & 13.4  & 22.9  & 33.2  & 54.7  & 79.0  & 78.8  & 75.8  & 69.2  & 62.8  & 29.1  & 30.0  & 26.4  & 16.7  & 12.7  & 11.1  & 11.3 \\
\cmidrule{2-23}      & \multicolumn{1}{r}{75} &  $S_{2.0}$  & 32.9  & 11.3  & 9.8   & 9.2   & 10.2  & 10.2  & 8.7   & 4.5   & 0.1   & 1.5   & 0.0   & 0.0   & 0.0   & 13.8  & 11.9  & 12.2  & 12.4  & 11.0  & 12.2  & 13.0 \\
      &       &  $S_{1.9}$  & 21.1  & 9.6   & 9.7   & 10.2  & 10.8  & 11.2  & 9.1   & 5.8   & 1.5   & 0.9   & 0.9   & 1.3   & 1.4   & 10.4  & 11.7  & 11.4  & 9.3   & 9.0   & 9.0   & 9.2 \\
      &       &  $S_{1.7}$  & 5.1   & 9.7   & 10.8  & 11.2  & 11.9  & 13.1  & 15.7  & 21.3  & 33.3  & 33.5  & 32.3  & 29.4  & 28.9  & 11.3  & 13.8  & 14.2  & 12.2  & 9.1   & 10.7  & 10.6 \\
\bottomrule
	\end{tabular}
	\end{threeparttable}
\end{table}

\begin{table}[p]
	\centering\footnotesize
	\setlength{\tabcolsep}{3pt}
	\begin{threeparttable}
	\caption{Percentage of rejection of $\mathcal{H}_0$ with $\alpha=1$
					 against stable alternatives.
					 Tests done at the 10\% significance level.}
	\label{tab:alpha1_vs_stable_estimate_FULL}
	\begin{tabular}{cccHHHHHHrrrrHHH|HrrrrrH}
		\toprule
		$p$ & $n$ & $\mathcal{H}_1$ & $T_{0.001}$ & $T_{0.005}$ & $T_{0.01}$ & $T_{0.05}$ & $T_{0.1}$
							& $T_{0.5}$ & $T_{1}$ & $T_{2}$ & $T_{5}$ & $T_{10}$ & $T_{20}$ & $T_{50}$ & $T_{100}$
						& $M_{2}$ & $M_{4}$ & $M_{6}$ & $M_{10}$ & $M_{15}$ & $M_{20}$ & $M_{25}$ \\
		\midrule
 \multicolumn{1}{r}{4} & \multicolumn{1}{r}{50} &  $S_{1.7}$  & 93.6  & 92.2  & 90.7  & 58.6  & 13.9  & 38.6  & 71.6  & 79.3  & 98.4  & 99.6  & 99.6  & 85.2  & 0.0   & 92.8  & 98.2  & 95.1  & 75.8  & 45.3  & 25.6  & 17.9 \\
      &       &  $S_{1.5}$  & 83.3  & 80.6  & 77.6  & 40.1  & 12.3  & 26.9  & 48.8  & 48.5  & 77.7  & 92.0  & 90.5  & 54.5  & 0.0   & 69.2  & 88.2  & 77.2  & 50.8  & 26.9  & 18.3  & 15.4 \\
      &       &  $S_{1.2}$  & 40.2  & 36.8  & 33.6  & 14.6  & 8.3   & 12.8  & 17.6  & 13.6  & 15.6  & 26.3  & 24.2  & 8.0   & 0.4   & 17.8  & 31.2  & 22.9  & 16.9  & 13.3  & 14.2  & 9.6 \\
      &       &  $S_{0.8}$  & 1.6   & 3.0   & 4.0   & 11.7  & 14.9  & 20.6  & 25.0  & 36.7  & 57.5  & 61.9  & 63.3  & 64.2  & 60.0  & 32.7  & 18.0  & 13.4  & 12.7  & 9.5   & 10.0  & 10.7 \\
      &       &  $S_{0.5}$  & 8.9   & 20.5  & 30.0  & 60.4  & 70.8  & 88.8  & 94.8  & 99.0  & 100.0 & 100.0 & 100.0 & 100.0 & 99.9  & 95.7  & 87.8  & 75.8  & 26.6  & 15.1  & 13.7  & 10.1 \\
\cmidrule{2-23}      & \multicolumn{1}{r}{75} &  $S_{1.7}$  & 98.8  & 98.2  & 96.7  & 37.8  & 10.7  & 62.9  & 90.8  & 97.7  & 100.0 & 100.0 & 100.0 & 99.4  & 76.2  & 99.6  & 99.9  & 99.4  & 90.3  & 59.2  & 39.0  & 29.1 \\
      &       &  $S_{1.5}$  & 94.6  & 93.1  & 88.8  & 21.8  & 8.8   & 42.9  & 71.8  & 79.8  & 96.9  & 99.4  & 99.0  & 89.7  & 34.1  & 89.8  & 96.4  & 91.3  & 61.8  & 33.3  & 18.1  & 15.7 \\
\bottomrule
	\end{tabular}
	\end{threeparttable}
\end{table}

\begin{table}[p]
	\centering\footnotesize
	\begin{threeparttable}
	\setlength{\tabcolsep}{3pt}
	\caption{Percentage of rejection of $\mathcal{H}_0$ with $\alpha=1$
					 against Student $t$ alternatives.
					 Tests done at the 10\% significance level.}
	\label{tab:alpha1_vs_t_estimate}
	\begin{tabular}{cccHHHHHHrrrrHHH|HrrrrrH}
		\toprule
		$p$ & $n$ & $\mathcal{H}_1$ & $T_{0.001}$ & $T_{0.005}$ & $T_{0.01}$ & $T_{0.05}$ & $T_{0.1}$
							& $T_{0.5}$ & $T_{1}$ & $T_{2}$ & $T_{5}$ & $T_{10}$ & $T_{20}$ & $T_{50}$ & $T_{100}$
						& $M_{2}$ & $M_{4}$ & $M_{6}$ & $M_{10}$ & $M_{15}$ & $M_{20}$ & $M_{25}$\\
		\midrule
\multicolumn{1}{r}{4} & \multicolumn{1}{r}{50} &  $t_{10}$  & 94.7  & 94.0  & 91.9  & 60.4  & 10.7  & 42.7  & 80.0  & 91.0  & 100.0 & 100.0 & 100.0 & 100.0 & 0.0   & 97.3  & 99.6  & 99.3  & 91.5  & 62.9  & 40.5  & 27.9 \\
      &       &  $t_{ 5}$  & 87.8  & 87.1  & 84.6  & 48.5  & 13.8  & 35.3  & 66.3  & 74.0  & 96.8  & 99.8  & 99.8  & 99.0  & 0.0   & 86.5  & 97.4  & 94.0  & 82.2  & 51.5  & 33.6  & 22.8 \\
      &       &  $t_{ 4}$  & 85.4  & 83.0  & 80.8  & 42.6  & 13.0  & 28.1  & 55.5  & 60.4  & 91.3  & 99.2  & 99.8  & 97.6  & 0.0   & 78.0  & 93.3  & 89.6  & 67.9  & 40.9  & 26.0  & 20.7 \\
      &       &  $t_{ 2}$  & 52.8  & 48.5  & 42.6  & 18.0  & 8.2   & 14.6  & 24.8  & 22.4  & 35.0  & 62.2  & 67.2  & 37.8  & 0.0   & 32.0  & 51.0  & 43.2  & 28.9  & 22.3  & 16.1  & 10.8 \\
\cmidrule{2-23}      & \multicolumn{1}{r}{75} &  $t_{10}$  & 99.8  & 99.4  & 98.9  & 46.5  & 9.2   & 67.8  & 95.3  & 99.7  & 100.0 & 100.0 & 100.0 & 100.0 & 100.0 & 100.0 & 100.0 & 100.0 & 98.8  & 80.6  & 52.8  & 37.3 \\
      &       &  $t_{ 5}$  & 97.0  & 96.8  & 95.6  & 54.2  & 13.0  & 51.1  & 82.8  & 94.0  & 100.0 & 100.0 & 100.0 & 100.0 & 82.8  & 97.6  & 99.4  & 99.2  & 93.4  & 64.9  & 44.3  & 31.9 \\
      &       &  $t_{ 4}$  & 95.8  & 95.0  & 91.8  & 28.4  & 8.9   & 43.7  & 77.0  & 90.2  & 99.8  & 100.0 & 100.0 & 100.0 & 90.6  & 94.6  & 99.1  & 97.8  & 87.4  & 55.7  & 37.9  & 33.5 \\
\bottomrule
	\end{tabular}
	\end{threeparttable}
\end{table}

\begin{table}[p]
	\centering\footnotesize
	\begin{threeparttable}
	\setlength{\tabcolsep}{3pt}
	\caption{Percentage of rejection of $\mathcal{H}_0$ with $\alpha=1$
					 against skew Cauchy alternatives.
					 Tests done at the 10\% significance level.}
	\label{tab:alpha1_vs_skewC}
	\begin{tabular}{cccHHHHHHrrrrHHH|HrrrrrH}
		\toprule
		$p$ & $n$ & $\mathcal{H}_1$ & $T_{0.001}$ & $T_{0.005}$ & $T_{0.01}$ & $T_{0.05}$ & $T_{0.1}$
							& $T_{0.5}$ & $T_{1}$ & $T_{2}$ & $T_{5}$ & $T_{10}$ & $T_{20}$ & $T_{50}$ & $T_{100}$
						& $M_{2}$ & $M_{4}$ & $M_{6}$ & $M_{10}$ & $M_{15}$ & $M_{20}$ & $M_{25}$ \\
		\midrule
		\multicolumn{1}{r}{4} & \multicolumn{1}{r}{50} &  $SC_{2.0}$  & 12.1  & 12.9  & 13.1  & 14.6  & 15.2  & 25.1  & 35.1  & 39.8  & 30.1  & 47.2  & 60.1  & 58.5  & 44.8  & 32.0  & 36.0  & 23.6  & 11.7  & 13.2  & 9.1   & 12.3 \\
      &       &  $SC_{1.5}$  & 11.2  & 12.1  & 12.7  & 14.2  & 15.4  & 23.1  & 31.0  & 35.2  & 25.9  & 42.7  & 55.5  & 53.5  & 41.0  & 29.4  & 32.0  & 22.4  & 13.1  & 10.0  & 11.7  & 10.4 \\
      &       &  $SC_{1.0}$  & 11.0  & 11.5  & 11.3  & 12.3  & 13.2  & 19.6  & 24.5  & 27.5  & 21.2  & 33.8  & 43.0  & 42.6  & 33.5  & 22.6  & 25.7  & 18.8  & 14.3  & 13.4  & 11.6  & 9.4 \\
      &       &  $SC_{0.5}$  & 11.8  & 11.1  & 10.9  & 11.6  & 11.9  & 15.1  & 15.8  & 17.9  & 15.4  & 18.5  & 23.0  & 23.9  & 20.9  & 14.6  & 16.6  & 14.9  & 10.7  & 10.5  & 8.0   & 11.6 \\
\cmidrule{2-23}      & \multicolumn{1}{r}{75} &  $SC_{2.0}$  & 13.3  & 14.0  & 14.7  & 14.2  & 14.2  & 33.2  & 52.8  & 63.7  & 57.7  & 77.6  & 88.7  & 88.2  & 80.1  & 46.0  & 48.5  & 30.0  & 13.0  & 10.5  & 10.3  & 7.7 \\
      &       &  $SC_{1.5}$  & 12.4  & 13.6  & 13.7  & 13.9  & 13.4  & 28.7  & 46.5  & 56.8  & 50.1  & 71.4  & 83.4  & 83.7  & 74.9  & 40.2  & 44.1  & 26.3  & 14.1  & 12.3  & 10.3  & 9.1 \\
      &       &  $SC_{1.0}$  & 12.2  & 12.6  & 13.4  & 13.4  & 13.3  & 22.5  & 35.2  & 43.2  & 37.0  & 54.1  & 69.8  & 70.6  & 61.8  & 31.1  & 32.8  & 21.6  & 12.1  & 8.4   & 8.7   & 8.2 \\
\bottomrule
	\end{tabular}
	\end{threeparttable}
\end{table}

\begin{table}[p]
	\centering\footnotesize
	\begin{threeparttable}
	\caption{Percentage of rejection of $\mathcal{H}_0$ with $\alpha=2$
					 against stable alternatives when using the test in
                      (8) with the Kotz-type weight function.
					 Tests done at the 10\% significance level.}
	\label{tab:alpha2_vs_stable_estimate_Kotz}
	\begin{tabular}{cccrHHHHHrrrrrrHHHH}
		\toprule
		$p$ & $n$ & $\mathcal{H}_1$ &
            $\nu$ & $T_{\nu,0.001}$ & $T_{\nu,0.005}$ & $T_{\nu,0.01}$ & $T_{\nu,0.05}$ & $T_{\nu,0.1}$ & $T_{\nu,0.5}$ & $T_{\nu,1}$ & $T_{\nu,2}$ & $T_{\nu,5}$ & $T_{\nu,10}$ & $T_{\nu,20}$ & $T_{\nu,50}$ & $T_{\nu,100}$ \\
		\midrule
        % Directory: C:/Documents/Papers/Working/StableLaws/Kode.v2.2.1/Output 
        % Pattern: ^test_Kotz_2.2.1_p[[:digit:]]_null2_altdist1.*txt$
		2 &  50 & $S_{1.9}$ & 0 & 12.7 & 15.8 & 16.6 & 23.0 & 25.9 & 32.1 & 35.0 & 37.6 & 40.6 & 42.1 & 42.7 & 43.0 & 42.9 \\
 &  &  & 1 & 11.4 & 11.6 & 13.9 & 15.3 & 14.7 & 21.8 & 28.5 & 35.0 & 41.5 & 43.0 & 43.5 & 44.2 & 44.2 \\
 &  &  & 2 & 13.0 & 10.5 & 11.9 & 14.5 & 13.9 & 19.4 & 26.8 & 34.8 & 42.4 & 44.7 & 46.5 & 46.1 & 46.3 \\
 &  &  & 3 &  9.4 & 13.8 & 10.6 & 10.9 & 13.5 & 17.2 & 19.1 & 26.5 & 38.9 & 41.1 & 43.0 & 43.8 & 44.2 \\
\cmidrule(lr){3-17}
 &  & $S_{1.8}$ & 0 & 17.8 & 22.9 & 27.9 & 34.6 & 41.0 & 54.2 & 60.7 & 65.3 & 66.7 & 67.0 & 67.5 & 67.6 & 67.3 \\
 &  &  & 1 & 12.6 & 13.4 & 15.6 & 19.8 & 23.5 & 40.2 & 47.8 & 57.7 & 65.1 & 67.1 & 67.5 & 67.5 & 67.3 \\
 &  &  & 2 & 15.3 & 14.4 & 16.5 & 15.4 & 18.6 & 34.5 & 45.2 & 56.4 & 66.4 & 68.4 & 69.0 & 69.1 & 68.9 \\
\bottomrule
	\end{tabular}
	\end{threeparttable}
\end{table}

\begin{table}[p]
	\centering\footnotesize
	\begin{threeparttable}
	\caption{Percentage of rejection of $\mathcal{H}_0$ with $\alpha=1$
					 against stable alternatives when using the test in
                      (8) with the Kotz-type weight function.
					 Tests done at the 10\% significance level.}
	\label{tab:alpha1_vs_stable_estimate_Kotz}
	\begin{tabular}{cccrHHHHHrrrrrrHHHH}
		\toprule
		$p$ & $n$ & $\mathcal{H}_1$ &
            $\nu$ & $T_{\nu,0.001}$ & $T_{\nu,0.005}$ & $T_{\nu,0.01}$ & $T_{\nu,0.05}$ & $T_{\nu,0.1}$ & $T_{\nu,0.5}$ & $T_{\nu,1}$ & $T_{\nu,2}$ & $T_{\nu,5}$ & $T_{\nu,10}$ & $T_{\nu,20}$ & $T_{\nu,50}$ & $T_{\nu,100}$ \\
		\midrule
		% Directory: C:/Documents/Papers/Working/StableLaws/Kode.v2.2.1/Output 
% Pattern: ^test_Kotz_2.2.1_p[[:digit:]]_null1_altdist1.*txt$
2 &  50 & $S_{1.4}$ & 0 & 29.4 & 11.9 & 11.7 & 11.4 & 12.0 & 26.3 & 33.9 & 37.2 & 54.9 & 54.3 & 37.3 &  0.1 &  0.0 \\
 &  &  & 1 & 13.2 & 12.8 & 12.0 &  9.8 & 10.4 & 13.4 & 18.9 & 26.5 & 31.4 & 56.0 & 52.1 & 20.2 &  0.1 \\
 &  &  & 2 & 11.3 & 10.3 &  9.0 & 11.1 & 11.1 & 10.1 & 15.1 & 29.1 & 16.5 & 42.7 & 57.9 & 41.2 &  1.4 \\
 &  &  & 3 & 72.8 & 11.4 & 12.3 &  8.9 & 11.1 & 10.7 &  9.9 & 23.6 & 21.8 & 22.3 & 59.7 & 49.3 & 13.0 \\
\cmidrule(lr){3-17}
 &  & $S_{1.3}$ & 0 & 12.5 &  8.1 &  7.8 &  8.5 & 11.0 & 21.9 & 23.3 & 20.8 & 34.5 & 34.1 & 21.8 &  0.3 &  0.5 \\
 &  &  & 1 & 16.5 & 13.3 & 13.2 &  9.8 &  9.3 & 10.8 & 15.8 & 18.4 & 18.7 & 40.9 & 42.1 & 12.4 &  0.2 \\
 &  &  & 2 &  8.9 &  9.9 & 10.3 & 10.4 &  9.0 &  8.7 & 13.6 & 18.6 & 10.0 & 27.1 & 39.9 & 21.7 &  0.3 \\
\bottomrule
	\end{tabular}
	\end{threeparttable}
\end{table}

\begin{table}[tbp]
	\centering\footnotesize
	\begin{threeparttable}
	\setlength{\tabcolsep}{3pt}
	\caption{Percentage of rejection of $\mathcal{H}_0$ with $\alpha=2$
             against stable alternatives in the CCC-GARCH$(1,1)$ case.
             Tests done at the 10\% significance level.}
	\label{tab:alpha2_vs_stable_garch_FULL}
	\begin{tabular}{cccHHHHHHrrrrHHH|HrrrrrHr}
		\toprule
		$p$ & $n$ & $\mathcal{H}_1$ & $T_{0.001}$ & $T_{0.005}$ & $T_{0.01}$ & $T_{0.05}$ & $T_{0.1}$
							& $T_{0.5}$ & $T_{1}$ & $T_{2}$ & $T_{5}$ & $T_{10}$ & $T_{20}$ & $T_{50}$ & $T_{100}$
						& $M_{2}$ & $M_{4}$ & $M_{6}$
						& $M_{10}$ & $M_{15}$ & $M_{20}$
                        & $M_{25}$ & HJM \\
\midrule
\multicolumn{1}{r}{4} & \multicolumn{1}{r}{100} &  $S_{2.0}$  & 10.9  & 10.4  & 9.4   & 8.1   & 8.2   & 8.9   & 9.5   & 9.4   & 8.5   & 7.9   & 7.6   & 7.9   & 8.1   & 9.2   & 9.6   & 8.7   & 10.0  & 9.3   & 10.1  & 11.1  & 10.2 \\
      &       &  $S_{1.9}$  & 10.7  & 11.3  & 16.3  & 18.0  & 23.0  & 34.0  & 33.8  & 29.3  & 19.2  & 11.2  & 7.5   & 7.2   & 7.2   & 5.5   & 13.2  & 15.6  & 18.9  & 17.6  & 14.8  & 14.8  & 73.9 \\
      &       &  $S_{1.8}$  & 23.0  & 25.0  & 31.3  & 38.7  & 51.1  & 68.9  & 69.7  & 62.9  & 37.0  & 17.0  & 8.2   & 7.2   & 8.0   & 5.1   & 28.5  & 36.7  & 32.4  & 32.5  & 29.6  & 25.7  & 96.0 \\
      &       &  $S_{1.7}$  & 26.9  & 31.6  & 50.0  & 65.0  & 76.7  & 90.8  & 92.0  & 85.9  & 52.8  & 21.8  & 8.0   & 7.6   & 8.2   & 9.0   & 53.9  & 61.8  & 58.5  & 52.2  & 46.9  & 38.8  & 99.9 \\
\cmidrule{2-24}      & \multicolumn{1}{r}{150} &  $S_{2.0}$  & 10.4  & 10.8  & 12.3  & 10.4  & 11.1  & 9.1   & 8.8   & 7.4   & 7.3   & 7.6   & 7.3   & 7.6   & 7.0   & 15.0  & 13.1  & 11.3  & 11.2  & 11.3  & 9.8   & 11.1  & 12.5 \\
      &       &  $S_{1.9}$  & 13.8  & 16.7  & 23.9  & 29.3  & 36.0  & 48.0  & 47.7  & 42.4  & 27.6  & 14.8  & 8.1   & 6.3   & 7.4   & 4.3   & 15.6  & 20.0  & 20.5  & 18.3  & 19.8  & 14.6  & 89.0 \\
      &       &  $S_{1.8}$  & 23.6  & 28.2  & 45.6  & 57.1  & 67.3  & 84.8  & 85.5  & 82.3  & 59.2  & 30.3  & 10.9  & 6.2   & 6.7   & 4.8   & 37.0  & 45.3  & 44.6  & 39.1  & 34.6  & 32.4  & 98.6 \\
\bottomrule
	\end{tabular}
	\end{threeparttable}
\end{table}

\begin{table}[p]
	\centering\footnotesize
	\begin{threeparttable}
	\caption{Percentage of rejection of $\mathcal{H}_0$ with $\alpha=2$
					 against $t$ alternatives (GARCH model errors).
					 Tests done at the 10\% significance level.}
	\label{tab:alpha2_vs_t_garch}
	\begin{tabular}{cccHHHHrrrrrHHHH|HrrH}
		\toprule
		$p$ & $n$ & $\mathcal{H}_1$ & $T_{0.001}$ & $T_{0.005}$ & $T_{0.01}$ & $T_{0.05}$ & $T_{0.1}$
							& $T_{0.5}$ & $T_{1}$ & $T_{2}$ & $T_{5}$ & $T_{10}$ & $T_{20}$ & $T_{50}$ & $T_{100}$
						& $M_{2}$ & $M_{4}$ & $M_{6}$ & HJM \\
		\midrule
4 & 100 & $t_{10}$ & 10.7 & 10.8 &  10.8 &  14.1 &  16.5 &  28.3 &  24.9 &  17.4 &  10.1 &  9.3 &  8.9 &  9.7 &  9.9 &   4.5 &  18.2 &  21.8 \\
 &  & $t_{ 5}$ & 18.5 & 19.3 &  26.3 &  38.7 &  57.3 &  81.1 &  79.4 &  58.2 &  15.4 &  8.9 &  7.6 &  7.9 &  8.8 &  10.1 &  51.1 &  54.5 \\
 &  & $t_{ 4}$ & 22.9 & 25.6 &  40.5 &  60.5 &  79.3 &  95.8 &  95.4 &  84.7 &  30.7 &  9.8 &  8.1 &  8.8 &  9.0 &  23.1 &  72.2 &  73.4 \\
 &  & $t_{ 2}$ & 72.9 & 83.8 &  96.5 &  99.2 & 100.0 & 100.0 & 100.0 &  99.9 &  90.9 & 37.7 &  5.1 &  4.6 &  5.5 &  91.9 &  99.9 & 100.0 \\
\cmidrule(lr){2-19}
 & 150 & $t_{10}$ & 15.0 & 15.6 &  15.3 &  18.3 &  24.6 &  44.5 &  40.0 &  25.6 &   9.1 &  7.0 &  6.7 &  6.9 &  7.7 &   7.0 &  23.9 &  26.0 \\
 &  & $t_{ 5}$ & 23.1 & 25.4 &  42.9 &  59.4 &  76.7 &  94.9 &  94.7 &  85.8 &  28.4 &  9.5 &  7.8 &  8.4 &  8.8 &  29.8 &  71.5 &  71.0 \\
 &  & $t_{ 4}$ & 21.3 & 32.5 &  62.3 &  82.4 &  94.8 &  99.7 &  99.4 &  96.6 &  53.8 & 12.9 &  7.5 &  8.4 &  9.1 &  52.2 &  90.6 &  90.7 \\
\bottomrule
	\end{tabular}
	\end{threeparttable}
\end{table}

\begin{table}[p]
	\centering\footnotesize
	\begin{threeparttable}
	\caption{Percentage of rejection of $\mathcal{H}_0$ with $\alpha=2$
             against skew normal alternatives in the CCC-GARCH$(1,1)$ case.
             Tests done at the 10\% significance level.
             Due to the computational complexity of the HJM test,
             it is only included in selected cases to reduce run time.}
	\label{tab:alpha2_vs_skewN_garch}
	\begin{tabular}{cccHHHHrrrrrHHHH|Hrrr}
		\toprule
		$p$ & $n$ & $\mathcal{H}_1$ & $T_{0.001}$ & $T_{0.005}$ & $T_{0.01}$ & $T_{0.05}$ & $T_{0.1}$
							& $T_{0.5}$ & $T_{1}$ & $T_{2}$ & $T_{5}$ & $T_{10}$ & $T_{20}$ & $T_{50}$ & $T_{100}$
						& $M_{2}$ & $M_{4}$ & $M_{6}$ & HJM \\
		\midrule
  4 & 50 & $SN_{0.4}$ &  9.3 & 10.4 & 10.6 & 12.6 & 18.8 & 54.5 & 71.3 & 81.3 & 87.1 & 87.7 & 88.1 & 88.0 & 88.2 & 18.9 & 47.7 & 52.5 & 12.6 \\
 &  & $SN_{0.3}$ & 11.2 & 13.0 & 13.3 & 15.0 & 15.7 & 38.0 & 51.9 & 62.4 & 67.6 & 67.7 & 67.7 & 67.6 & 67.9 & 14.3 & 31.9 & 34.2 & 11.1 \\
 &  & $SN_{0.2}$ &  9.7 &  9.2 &  9.0 &  9.5 & 10.2 & 21.1 & 28.8 & 34.8 & 35.9 & 36.6 & 36.8 & 36.6 & 35.9 & 10.3 & 16.0 & 18.3 & 11.5 \\
 &  & $SN_{0.1}$ & 10.5 &  9.8 & 10.2 & 11.8 & 13.6 & 12.1 & 14.9 & 17.1 & 18.5 & 18.4 & 18.5 & 17.8 & 17.3 &  9.5 & 11.9 & 13.4 &  9.6 \\
 % &  & $SN_{0.0}$ & 12.5 & 12.0 & 12.0 & 12.1 &  9.8 & 10.8 &  10.2 &  12.0 &  12.0 &  12.9 &  12.4 &  12.7 &  12.6 & 10.2 &  9.9 &  9.6 \\
\cmidrule(lr){2-20}
  & 100 & $SN_{0.4}$ & 12.6 & 13.2 & 13.8 & 19.9 & 32.8 & 85.9 & 97.5 & 99.2 & 99.8 & 99.9 &  99.9 & 100.0 & 100.0 & 38.1 & 75.8 & 81.7 & 21.4 \\
 &  & $SN_{0.3}$ &  9.8 & 11.1 & 15.7 & 16.9 & 24.0 & 65.7 & 84.4 & 93.5 & 95.5 & 96.1 &  96.8 &  97.0 &  97.0 & 21.8 & 51.8 & 61.6 & 21.9 \\
 &  & $SN_{0.2}$ & 10.8 & 11.0 & 13.1 & 13.3 & 17.2 & 43.1 & 54.3 & 64.4 & 71.1 & 74.0 &  74.7 &  74.5 &  74.7 & 13.0 & 28.9 & 36.2 & 13.9 \\
\bottomrule
	\end{tabular}
	\end{threeparttable}
\end{table}

\begin{table}[p]
	\centering\footnotesize
	\begin{threeparttable}
	\caption{Percentage of rejection of $\mathcal{H}_0$ with $\alpha=1.8$
             against stable alternatives in the CCC-GARCH$(1,1)$ case.
					 Tests done at the 10\% significance level.}
	\label{tab:alpha1.8_vs_stable_garch}
	\begin{tabular}{cccHHHrrrrrHHHHH|Hrrr}
		\toprule
		$p$ & $n$ & $\mathcal{H}_1$ & $T_{0.001}$ & $T_{0.005}$ & $T_{0.01}$ & $T_{0.05}$ & $T_{0.1}$
							& $T_{0.5}$ & $T_{1}$ & $T_{2}$ & $T_{5}$ & $T_{10}$ & $T_{20}$ & $T_{50}$ & $T_{100}$
						& $M_{2}$ & $M_{4}$ & $M_{6}$ \\
		\midrule
  4 & 100 & $S_{2.0}$ & 12.0 & 11.0 &  9.9 &  8.8 & 11.3 & 16.4 & 15.7 & 12.2 & 11.5 & 12.7 & 11.9 &  7.9 &  4.9 & 19.8 & 13.8 & 14.5 \\
 &  & $S_{1.9}$ & 10.0 & 10.0 & 10.5 & 11.3 & 10.2 & 14.4 & 14.1 & 12.5 & 11.2 & 10.9 &  9.8 &  8.4 &  6.4 & 13.7 & 10.5 &  9.8 \\
 &  & $S_{1.8}$ &  6.0 &  6.2 &  8.1 &  9.1 &  9.1 & 10.5 &  9.0 &  8.3 &  8.4 &  9.7 &  9.8 &  9.0 &  9.5 &  9.7 &  7.9 &  8.1 \\
 &  & $S_{1.7}$ &  8.2 &  8.5 &  8.2 &  8.5 &  9.1 & 12.4 & 11.7 & 10.9 &  9.0 &  9.6 & 11.2 & 14.2 & 16.8 &  8.6 & 12.1 & 10.0 \\
 &  & $S_{1.6}$ & 11.2 & 11.4 & 12.7 & 11.2 & 11.9 & 18.5 & 17.2 & 14.0 &  9.7 & 11.1 & 18.2 & 24.7 & 28.8 & 15.1 & 24.3 & 22.3 \\
 &  & $S_{1.5}$ & 13.7 & 16.7 & 17.0 & 19.6 & 26.0 & 37.8 & 32.8 & 19.6 & 12.6 & 18.8 & 30.0 & 42.9 & 48.0 & 30.0 & 46.2 & 44.5 \\
\cmidrule(lr){2-19}
 & 150 & $S_{2.0}$ &  8.7 &  8.3 &  9.7 &  9.2 & 11.3 & 25.1 & 24.1 & 17.7 & 12.9 & 11.7 & 10.7 &  6.8 &  4.1 & 22.9 & 20.3 & 20.3 \\
\bottomrule
	\end{tabular}
	\end{threeparttable}
\end{table}

\begin{table}[p]
	\centering\footnotesize
	\begin{threeparttable}
	\caption{Percentage of rejection of $\mathcal{H}_0$ with $\alpha=1.5$
					 against stable alternatives (GARCH model errors).
					 Tests done at the 10\% significance level.}
	\label{tab:alpha1.5_vs_stable_garch}
	\begin{tabular}{cccHHHHrrrrrHHHH|Hrrr}
		\toprule
		$p$ & $n$ & $\mathcal{H}_1$ & $T_{0.001}$ & $T_{0.005}$ & $T_{0.01}$ & $T_{0.05}$ & $T_{0.1}$
							& $T_{0.5}$ & $T_{1}$ & $T_{2}$ & $T_{5}$ & $T_{10}$ & $T_{20}$ & $T_{50}$ & $T_{100}$
						& $M_{2}$ & $M_{4}$ & $M_{6}$ \\
		\midrule
4 & 100 & $S_{1.8}$ & 55.7 &  9.7 &  7.0 &  7.9 &  8.9 & 34.2 & 36.5 & 24.2 & 17.9 & 19.2 & 21.4 &  9.8 &  2.8 & 32.4 & 41.0 & 40.9 \\
 &  & $S_{1.7}$ & 24.4 &  9.6 &  8.7 &  8.6 &  9.8 & 24.9 & 24.3 & 18.7 & 17.1 & 17.4 & 16.5 &  8.8 &  2.8 & 20.5 & 25.5 & 21.7 \\
 &  & $S_{1.6}$ & 13.3 &  8.4 &  7.7 &  7.5 &  8.4 & 10.9 & 13.3 & 13.7 & 13.4 & 11.9 &  9.9 &  6.2 &  6.3 &  9.8 & 13.3 & 11.6 \\
 &  & $S_{1.5}$ & 10.9 & 11.4 & 11.4 & 11.5 & 10.9 &  9.0 &  9.7 & 10.5 &  8.7 &  8.0 &  8.1 &  8.7 &  7.9 &  9.3 & 11.1 & 12.9 \\
 &  & $S_{1.4}$ &  8.7 & 10.1 & 10.3 & 10.9 & 10.4 & 10.6 & 13.4 & 11.3 & 11.7 & 14.8 & 19.3 & 23.6 & 24.8 & 15.1 & 14.7 & 14.6 \\
 &  & $S_{1.3}$ &  7.4 & 10.3 & 12.3 & 11.3 & 15.4 & 21.7 & 21.8 & 15.2 & 13.3 & 21.2 & 33.3 & 42.2 & 41.1 & 29.9 & 31.4 & 25.5 \\
 &  & $S_{1.2}$ &  8.0 & 12.0 & 15.8 & 17.0 & 24.6 & 44.8 & 41.7 & 25.6 & 15.7 & 32.0 & 52.3 & 63.5 & 65.8 & 56.5 & 52.3 & 44.5 \\
\bottomrule
	\end{tabular}
	\end{threeparttable}
\end{table}

\begin{table}[p]
	\centering\footnotesize
	\begin{threeparttable}
	\caption{Percentage of rejection of $\tilde{\mathcal{H}}_0$ with $\alpha=1.5$
					 against stable alternatives.
					 Tests done at the 10\% significance level.}
	\label{tab:skew_alpha1.5_FULL}
	\begin{tabular}{HccHHrrrrHHHHHHHHHrrrrHHHHHHHHrrrrr}
		\toprule
            &&&\multicolumn{13}{c}{Skew SP distributions}&\multicolumn{13}{c}{Symmetric SP distributions} \\
            \cmidrule(lr){4-16}\cmidrule(lr){17-29}
		$p$ & $n$ & $\mathcal{H}_1$ &
            $\tilde T_{0.001}$ & $\tilde T_{0.005}$ & $\tilde T_{0.01}$ & $\tilde T_{0.05}$ & $\tilde T_{0.1}$
							& $\tilde T_{0.5}$ & $\tilde T_{1}$ & $\tilde T_{2}$ & $\tilde T_{5}$ & $\tilde T_{10}$ & $\tilde T_{20}$ & $\tilde T_{50}$ & $\tilde T_{100}$
             & $\tilde T_{0.001}$ & $\tilde T_{0.005}$ & $\tilde T_{0.01}$ & $\tilde T_{0.05}$ & $\tilde T_{0.1}$
							& $\tilde T_{0.5}$ & $\tilde T_{1}$ & $\tilde T_{2}$ & $\tilde T_{5}$ & $\tilde T_{10}$ & $\tilde T_{20}$ & $\tilde T_{50}$ & $\tilde T_{100}$ \\
		\midrule
2     & 100   & $\tilde{\cal{S}}_{2.0}$   & 33.7  & 52.5  & 59.0  & 73.7  & 71.8  & 59.6  & 51.2  & 46.4  & 34.4  & 23.6  & 17.3  & 12.6  & 33.9  & 57.2  & 66.3  & 73.6  & 68.1  & 56.2  & 50.0  & 40.6  & 24.1  & 16.7  & 11.9  & 9.9 \\
      &       & $\tilde{\cal{S}}_{1.7}$   & 7.0   & 10.7  & 13.4  & 17.8  & 19.4  & 17.2  & 15.6  & 15.3  & 12.2  & 9.1   & 8.6   & 8.9   & 5.5   & 9.7   & 13.2  & 18.7  & 18.5  & 16.8  & 14.6  & 12.4  & 9.3   & 9.1   & 9.7   & 9.7 \\
      &       & $\tilde{\cal{S}}_{1.5}$   & 15.0  & 14.2  & 13.3  & 9.8   & 8.9   & 7.4   & 8.0   & 8.4   & 9.9   & 8.9   & 8.0   & 9.4   & 9.2   & 8.4   & 9.0   & 8.5   & 8.4   & 8.3   & 7.9   & 7.1   & 9.4   & 8.2   & 7.9   & 8.8 \\
      &       & $\tilde{\cal{S}}_{1.4}$  & 23.6  & 21.8  & 20.6  & 16.5  & 15.2  & 13.9  & 15.6  & 15.4  & 15.4  & 14.0  & 13.5  & 13.5  & 18.6  & 17.8  & 17.1  & 13.8  & 11.2  & 11.4  & 12.0  & 12.6  & 12.1  & 11.4  & 11.6  & 11.3 \\
      &       & $\tilde{\cal{S}}_{1.3}$   & 40.4  & 38.4  & 37.5  & 30.3  & 27.8  & 28.7  & 29.6  & 30.5  & 30.3  & 28.3  & 25.1  & 22.1  & 34.6  & 31.8  & 31.2  & 24.3  & 21.1  & 24.1  & 26.7  & 27.1  & 23.4  & 20.7  & 18.8  & 16.5 \\
\bottomrule
	\end{tabular}
	\end{threeparttable}
\end{table}

\begin{table}[p]
	\centering\footnotesize
	\begin{threeparttable}
	\caption{Percentage of rejection of $\tilde{\mathcal{H}}_0$ with $\alpha=1.8$
					 against stable alternatives.
					 Tests done at the 10\% significance level.}
	\label{tab:skew_alpha1.8}
	\begin{tabular}{HccHHrrrrHHHHHHHHHrrrrHHHHHHHHrrrrr}
		\toprule
            &&&\multicolumn{13}{c}{Skew SP distributions}&\multicolumn{13}{c}{Symmetric SP distributions} \\
            \cmidrule(lr){4-16}\cmidrule(lr){17-29}
		$p$ & $n$ & $\mathcal{H}_1$ &
            $\tilde T_{0.001}$ & $\tilde T_{0.005}$ & $\tilde T_{0.01}$ & $\tilde T_{0.05}$ & $\tilde T_{0.1}$
							& $\tilde T_{0.5}$ & $\tilde T_{1}$ & $\tilde T_{2}$ & $\tilde T_{5}$ & $\tilde T_{10}$ & $\tilde T_{20}$ & $\tilde T_{50}$ & $\tilde T_{100}$
             & $\tilde T_{0.001}$ & $\tilde T_{0.005}$ & $\tilde T_{0.01}$ & $\tilde T_{0.05}$ & $\tilde T_{0.1}$
							& $\tilde T_{0.5}$ & $\tilde T_{1}$ & $\tilde T_{2}$ & $\tilde T_{5}$ & $\tilde T_{10}$ & $\tilde T_{20}$ & $\tilde T_{50}$ & $\tilde T_{100}$ \\
		\midrule
2     & 100   & $\tilde{\cal{S}}_{2.0}$   & 4.8   & 8.7   & 12.1  & 17.1  & 16.0  & 12.9  & 13.2  & 10.9  & 11.1  & 10.0  & 8.6   & 9.3   & 4.0   & 10.1  & 14.2  & 16.3  & 13.4  & 13.0  & 13.5  & 11.2  & 9.0   & 8.4   & 7.5   & 8.5 \\
      &       & $\tilde{\cal{S}}_{1.9}$   & 6.8   & 7.2   & 7.2   & 9.2   & 8.7   & 7.6   & 8.5   & 8.4   & 9.7   & 8.9   & 9.5   & 8.0   & 5.1   & 5.8   & 6.8   & 7.9   & 6.5   & 9.0   & 9.8   & 9.7   & 9.0   & 8.6   & 8.3   & 10.3 \\
      &       & $\tilde{\cal{S}}_{1.8}$  & 14.1  & 14.5  & 12.6  & 8.6   & 7.6   & 8.3   & 8.3   & 7.3   & 8.5   & 8.9   & 10.1  & 10.4  & 8.7   & 7.5   & 7.6   & 5.5   & 5.1   & 7.8   & 8.2   & 7.8   & 8.6   & 8.9   & 10.1  & 8.6 \\
      &       & $\tilde{\cal{S}}_{1.7}$   & 21.8  & 23.5  & 21.1  & 18.4  & 15.2  & 14.4  & 12.1  & 13.2  & 10.8  & 10.0  & 8.6   & 10.0  & 15.9  & 15.7  & 15.6  & 10.8  & 9.8   & 10.0  & 10.3  & 8.7   & 6.4   & 6.0   & 8.5   & 10.0 \\
      &       & $\tilde{\cal{S}}_{1.6}$   & 34.8  & 39.8  & 37.6  & 31.0  & 26.2  & 26.8  & 24.8  & 25.1  & 19.7  & 16.3  & 13.5  & 11.5  & 30.1  & 33.9  & 33.4  & 23.0  & 22.5  & 25.5  & 22.0  & 21.0  & 16.4  & 11.2  & 11.5  & 11.5 \\
\bottomrule
	\end{tabular}
	\end{threeparttable}
\end{table}

\begin{table}[p]
	\centering\footnotesize
	\begin{threeparttable}
	\caption{Percentage of rejection of $\tilde{\mathcal{H}}_0'$ using tests based on
                  the statistic in (21) with $\Psi(\xi)=e^{-r\xi}$.
				 Tests done at the 10\% significance level.}
	\label{tab:skew_iid_estimate_alpha_FULL}
	\begin{tabular}{HlHrHHHHHrrrrrHHHHHrr}
		\toprule
		$p$ & $\mathcal{H}_1$ & & $n$ & $T_{0.001}$ & $T_{0.005}$ & $T_{0.01}$ & $T_{0.05}$ & $T_{0.1}$
							& $T_{0.5}$ & $T_{1}$ & $T_{2}$ & $T_{5}$ & $T_{10}$ & $T_{20}$ & $T_{50}$ & $T_{100}$ \\
		\midrule
            2 & $\tilde{{S}}_{1.8}$ & 1.8 &  50 &   5.5 &   5.6 &   5.8 &   5.2 &   5.1 &   4.3 &   5.2 &   5.1 &   6.7 &   6.8 &   7.8 &   7.8 &   8.3 \\
 &  & 1.8 &  75 &   6.6 &   5.2 &   5.7 &   5.5 &   5.0 &   6.7 &   7.5 &   8.1 &  10.2 &  10.1 &   9.7 &   9.2 &   9.8 \\
 &  & 1.8 & 100 &   5.5 &   4.7 &   4.8 &   4.7 &   5.2 &   6.5 &   7.7 &   8.0 &   8.2 &   8.8 &   9.5 &  10.0 &  11.0 \\
 &  & 1.8 & 150 &   5.6 &   6.9 &   7.2 &   7.4 &   7.1 &   7.8 &   9.0 &   9.0 &   9.2 &   9.2 &   9.5 &  11.1 &   9.7 \\
 &  & 1.8 & 200 &  10.0 &   7.7 &   7.1 &   7.2 &   6.6 &   9.2 &   9.6 &   8.7 &  10.3 &  11.7 &  12.6 &  12.3 &  11.2 \\
\bottomrule
	\end{tabular}
	\end{threeparttable}
\end{table}

\begin{table}[p]
	\centering\footnotesize
	\begin{threeparttable}
	\caption{Percentage of rejection of $\mathcal{H}_0$ with $\alpha=2$
                against stable alternatives.
                For these results, $\mat Q$ and $\vec \delta$
                were assumed known and not estimated.
                Tests done at the 10\% significance level.}
	\label{tab:alpha2_vs_stable_highdim}
	\begin{tabular}{cccHHHHHHrrrrrHH}
		\toprule
		$p$ & $n$ & $\mathcal{H}_1$ & $T_{0.001}$ & $T_{0.005}$ & $T_{0.01}$ & $T_{0.05}$ & $T_{0.1}$
							& $T_{0.5}$ & $T_{1}$ & $T_{2}$ & $T_{5}$ & $T_{10}$ & $T_{20}$ & $T_{50}$ & $T_{100}$ \\
  \midrule
50 &  250 & $S_{2.00}$ & 100.0 & 100.0 & 100.0 & 100.0 & 100.0 &  16.3 &  10.0 &  10.4 &  10.1 &   8.1 &   8.5 &  9.7 & 10.0 \\
 &  & $S_{1.97}$ & 100.0 & 100.0 & 100.0 & 100.0 & 100.0 &   3.6 &  15.7 &  18.2 &  21.6 &  13.1 &  11.0 & 11.4 & 12.6 \\
 &  & $S_{1.95}$ & 100.0 & 100.0 & 100.0 & 100.0 & 100.0 &   1.7 &  15.5 &  28.0 &  40.3 &  20.3 &  10.0 & 13.9 & 18.5 \\
\cmidrule(lr){2-16}
 &  500 & $S_{2.00}$ & 100.0 & 100.0 & 100.0 & 100.0 & 100.0 &  12.1 &  10.9 &   9.8 &   7.2 &   9.8 &  10.0 & 10.5 & 10.0 \\
 &  & $S_{1.97}$ & 100.0 & 100.0 & 100.0 & 100.0 & 100.0 &   3.5 &  14.9 &  21.4 &  36.9 &  20.3 &  11.8 & 11.8 & 14.2 \\
 &  & $S_{1.95}$ & 100.0 & 100.0 & 100.0 & 100.0 & 100.0 &   3.4 &  20.2 &  44.6 &  68.3 &  35.6 &  13.9 & 18.7 & 27.7 \\
\cmidrule(lr){2-16}
 &  750 & $S_{2.00}$ & 100.0 & 100.0 & 100.0 & 100.0 & 100.0 &  11.6 &  10.6 &   9.3 &   9.4 &  10.8 &  11.1 & 11.7 & 11.2 \\
\bottomrule
	\end{tabular}
	\end{threeparttable}
\end{table}

\begin{table}[p]
	\centering\footnotesize
	\begin{threeparttable}
	\caption{Percentage of rejection of $\mathcal{H}_0$ with $\alpha=1.5$
                against stable alternatives.
                For these results, $\mat Q$ and $\vec \delta$
                were assumed known and not estimated.
                Tests done at the 10\% significance level.}
	\label{tab:alpha1.5_vs_stable_highdim}
	\begin{tabular}{cccHHHHHHrrrrrHH}
		\toprule
		$p$ & $n$ & $\mathcal{H}_1$ & $T_{0.001}$ & $T_{0.005}$ & $T_{0.01}$ & $T_{0.05}$ & $T_{0.1}$
							& $T_{0.5}$ & $T_{1}$ & $T_{2}$ & $T_{5}$ & $T_{10}$ & $T_{20}$ & $T_{50}$ & $T_{100}$ \\
  \midrule
50 &  250 & $S_{1.60}$ & 100.0 & 100.0 & 100.0 & 71.4 & 71.6 & 75.3 & 78.5 & 60.5 &  89.4 &  75.0 &  23.2 & 17.0 & 25.6 \\
 &  & $S_{1.55}$ & 100.0 & 100.0 & 100.0 & 43.7 & 42.5 & 44.5 & 44.2 &  7.7 &  36.6 &  27.8 &  11.8 & 10.5 & 12.6 \\
 &  & $S_{1.50}$ & 100.0 & 100.0 & 100.0 & 11.7 & 10.6 & 11.2 &  9.8 & 10.5 &   9.0 &   9.6 &   9.4 &  9.2 &  9.1 \\
 &  & $S_{1.45}$ & 100.0 & 100.0 & 100.0 &  1.4 &  1.0 &  1.3 &  5.5 & 23.9 &  44.7 &  31.3 &  12.3 & 12.1 & 17.2 \\
 &  & $S_{1.40}$ & 100.0 & 100.0 & 100.0 &  0.0 &  0.0 &  0.8 & 15.7 & 57.1 &  84.8 &  70.1 &  24.2 & 16.3 & 32.0 \\
\cmidrule(lr){2-16}
 &  500 & $S_{1.60}$ & 100.0 & 100.0 & 100.0 & 84.2 & 84.4 & 87.2 & 90.7 & 55.0 & 100.0 &  95.3 &  34.8 & 16.7 & 38.9 \\
 &  & $S_{1.55}$ & 100.0 & 100.0 & 100.0 & 46.5 & 46.8 & 50.4 & 49.2 &  4.8 &  61.3 &  46.6 &  13.6 &  8.9 & 15.7 \\
\bottomrule
	\end{tabular}
	\end{threeparttable}
\end{table}

% \bibliographystyle{dcu}
%\bibliography{refs}